\documentclass[a4paper,10pt]{amsart}
\usepackage{amssymb}
\usepackage{amsthm}
\usepackage{zmacros}
\usepackage{url}
\usepackage[ps,dvips,2cell,matrix,arrow,line,tips]{xy} 
\usepackage{srcltx}
\usepackage{leftidx}
\usepackage{xspace}
\input{zmathmacros.cfg}
\input{abbreviations.cfg}
\title{Yetter-Drinfeld modules for Turaev crossed structures}     
\author{Marco Zunino}
\date{Friday, 13 September 2002}
\address{Marco Zunino -- 
        Institut de Recherche Math\'ematiques Avanc\'ee (IRMA) -- 
        Universit\'e Louis Pasteur et CNRS -- 
        7,~rue Ren\'e Descartes, 67084 Strasbourg \textsc{cedex}, France. }
\email{zunino@math.u-strasbg.fr}

\begin{document}

\bibliographystyle{amsplain}

\begin{abstract}
 We provide an analog of the Joyal-Street center construction
 and of the Kassel-Turaev categorical quantum double
 in the context of the crossed categories introduced by Turaev.
 Then, we focus or attention to the case of categories of representation.
 In particular, we introduce the notion of a Yetter-Drinfeld module over
 a crossed group coalgebra $H$ and we prove that both the category
 of Yetter-Drinfeld modules over $H$ and the center of the category
 of representations of $H$ as well as the category of representations
 of the quantum double of $H$ are isomorphic as braided crossed 
 categories.
\end{abstract}

\maketitle

\tableofcontents

\section*{Introduction}
 Recently, Turaev~\cite{Tur-pi,Tur-CPC} (see also~\cite{LeTur} 
 and~\cite{Virelizi2}) generalized quantum invatiants of
 $3$\nobreakdash-\hspace{0pt}manifold to the case of a 
 $3$\nobreakdash-\hspace{0pt}manifold $M$ endowed with 
 a homotopy class of maps $M\to K(\pi,1)$, where $\pi$ is a group
 (such homotopy classes of maps $M\to K(\pi,1)$ classify principal flat 
 $\pi$\nobreakdash-\hspace{0pt}bundles  over $M\/$).
 
 One of the key points of the theory~\cite{Tur-CPC} is a generalization 
 of the definition of a tensor category to the notion of a 
 \textit{crossed \picat,} here called a \textit{Turaev category} 
 or, briefly, a \textit{\Tcat.} 
 The algebraic counterpart is the generalization of the definition 
 of Hopf algebra to the notion of a \textit{crossed Hopf \picoalg,} 
 here called a \textit{Turaev coalgebra} or, briefly, a \textit{\Tcoalg.\/}  
 As the category of representations of
 a Hopf algebra has a structure of a tensor category, 
 the category of representations of a \Tcoalg has a structure of a 
 \Tcat. Concepts like braiding (universal R-matrix), ribbon, and
 modularity can be extended to the crossed case. Ribbon and modular
 \Tcats play a central role in the construction of the new invariants.
 
  Roughly speaking, a \textit{\Tcoalg $H$} is a family $\{H_{\alpha}\}_{\alpha\in\pi}$ of
 algebras endowed with a \textit{comultiplication}
 $\Delta_{\alpha,\beta}\colon H_{\alpha\beta}\to H_{\alpha}\otimes H_{\beta}$,
 a \textit{counit} $\varepsilon\colon\Bbbk\to H_{1}$ (where $1$ is the neutral 
 element of $\pi\/$),  and an \textit{antipode}
 $s_{\alpha}\colon H_{\alpha}\to H_{\alpha^{-1}}$.
 It is required that $H$ satisfies axioms that generalize 
 those of a Hopf algebra. It is also required that $H$
 is endowed with a family of algebra isomorphisms 
 $\varphi^{\alpha}_{\beta}=\varphi_{\beta}\colon H_{\alpha}\to H_{\beta\alpha\beta^{-1}}$, 
 the \textit{conjugation,} compatible with the above structures 
 and such that $\varphi_{\beta\gamma}=\varphi_{\beta}\circ\varphi_{\gamma}$.
  In particular, when $\pi=1$, we recover the usual definition 
  of a Hopf algebra.
 A \textit{\Tcat} is a tensor category $\mathcal{T}$ disjoint 
 union of a family of categories $\{\mathcal{T}_{\alpha}\}_{\alpha\in\pi}$
 such that, if $U\in\mathcal{T}_{\alpha}$ and $V\in\mathcal{T}_{\beta}$, then 
 $U\otimes V\in\mathcal{T}_{\alpha\beta}$.
 It is required that $\mathcal{T}$ is endowed with a group 
 homomorphism $\varphi\colon\pi\to\Aut(\mathcal{T})$, the \textit{conjugation,}
 where $\Aut(\mathcal{T})$ is the group of strict tensor 
 automorphisms of $\mathcal{T}$.
 Given $\alpha\in\pi$ and $U\in\mathcal{T}_{\alpha}$, the functor $\varphi_{\alpha}$ 
 is also denoted $\lidx{U}{(\cdot)}$.
 Notice that the component $\mathcal{T}_{1}$ is a tensor category. 
 In particular, when $\pi=1$, we recover the usual definition 
 of a tensor category.
 If $H$ is a \Tcoalg, then the disjoint 
 union $\Rep(H)=\coprod_{\alpha\in\pi}\Rep_{\alpha}(H)$ 
 of the categories of representations $\Rep_{\alpha}(H)=\Rep(H_{\alpha})$
 has a structure of a \Tcat.
 
 In~\cite{Zunino-CHC} we provide an analog of the Drinfeld 
 quantum double (see~\cite{Drn}) and of the ribbon extension 
 (see~\cite{RT}), that is we showed how to obtain a quasitriangular
 and a ribbon \Tcoalg starting from any finite-type \Tcoalg.
 In this article we focus our attention on the corresponding 
 categorical constructions and the relations between 
 these and the algebraic one.
 
 In the first part of the article, we study how to obtain braided
 or ribbon \Tcats starting from a \Tcat that is not braided.
 Firstly, given any \Tcat $\mathcal{T}$, we obtain a braided \Tcat 
 $\mathcal{Z}(\mathcal{T})$, the \textit{center of $\mathcal{T}$} 
 ({Theorem~\ref{thm:ZZ}}). When $\pi=1$, we recover the 
 definition of the center of a tensor category, see~\cite{JS}.
 Given a braided \Tcat $\mathcal{T}'$, by generalizing the constructions
 in~\cite{Street-double} and~\cite{KasTur}, we obtain a ribbon \Tcat
 $\mathcal{N}\bigl((\mathcal{T}')^{Z}\bigl)$. We can apply this 
 construction when
 $\mathcal{T}'=\mathcal{Z}(\mathcal{T})$ so that for any \Tcat
 $\mathcal{T}$ we obtain a ribbon \Tcat 
 $\mathcal{D}(\mathcal{T})=
 \mathcal{N}\Bigl(\bigl(\mathcal{Z}(\mathcal{T})\bigl)^{Z}\Bigr)$. 
 
 In the second part of the article, we consider the case of category
 of representations of \Tcoalgs and we discuss the relation between the
 above categorical constructions and the algebraic constructions 
 that we previously introduced in~\cite{Zunino-CHC}. 
 Firstly, we introduce the fundamental
 notion of a 
   \textit{Yetter-Drinfeld module} over a \Tcoalg $H$, or, 
   briefly, a \textit{\Yetter,} as a module $V$ 
  over a component $H_{\alpha}$ of $H$ endowed with a family of \klin morphisms
 $\Delta_{\beta}\colon V\to V\otimes H_{\beta}$
 (for any $\beta\in\pi\/$) satisfying axioms that generalize the usual 
 definition of a Yetter-Drinfeld module over a Hopf
 algebra. We state now the main result relating our
 double constructions for crossed categories and for crossed
 Hopf group coalgebras.\medskip
 
 \textit{\textsc{\bfseries Theorems \ref{thm:C-moll-Messe} and 
 \ref{thm:verpastwt}.}
 Let $H$ be a \Tcoalg of finite type. 
 We have the isomorphisms of braided \Tcats
 \begin{equation*}
 \mathcal{Z}\bigl(\Rep(H)\bigr)=\mathcal{YD}(H)=
 \Rep\bigl(\overline{D}(H)\bigr)
 \text{,}
\end{equation*}
 where $\YD\,(H)$ is the category of Yetter-Drinfeld modules over $H$ and 
 $\overline{D}(H)$ is the mirror of the quantum double $D(H)$
 of $H$ \textup{(\/}see~\cite{Zunino-CHC}\textup{).}}\medskip
 
 We can describe Theorems~\ref{thm:C-moll-Messe}
 and~\ref{thm:verpastwt} via the commutative diagram
 \begin{equation*}
 \vcenter{\xymatrix{\Rep(H)\ar[r]^{\mathcal{Z}\ \ \ \ } & \mathcal{Z}\bigl(\Rep(H)\bigr)\ar@{=}[r]
               & \YD\,(H)\ar@{=}[r] 
               & \Rep\bigl(\overline{D}(H)\bigr)\\ 
           H\ar[u]^{\Rep}\ar[rrr]_{\overline{D}} & \ & \ & \overline{D}(H)\ar[u]_{\Rep}}}\text{.}
\end{equation*}
 
Similar results are obtained for ribbon structures.\medskip

\paragraph{\scshape Acknowledgments} 
The author wants to thank his adviser, V.\@~Turaev for his
stimulating and constructive direction in the research. 
The author also wants to thank A.\@~Brugui\`eres, 
B.\@~Enriquez, Ch.\@~Kassel, and H.J.\@~Schneider
for the useful remarks.
The author was partially supported by the 
\textsc{INdAM, Istituto Nazionale di Alta Matematica,} Rome.

A special thanks to the little blond angel that invited me to
return to work with the kindness of her smile and the persuasion
of a well pointed trident. 

\section{Tensor categories}\label{tensorcat}
To fix our notations, let us recall few basic definitions.
A \textit{tensor category} $\mathcal{C}=(\mathcal{C},\otimes,a,l,r)$ 
(see~\cite{ML,ML-coherence}), also called a \textit{monoidal category,}
is a category $\mathcal{C}$ endowed with a functor 
$\otimes\colon\mathcal{C}\times\mathcal{C}\to\mathcal{C}$ 
(the \textit{tensor product\/}), an object 
$\tunit\in\mathcal{C}$ (the \textit{tensor unit\/}) 
and natural isomorphisms
$a=a_{U,V,W}\colon (U\otimes V)\otimes W\to U\otimes (V\otimes W)$
for all $U,V,W\in\mathcal{C}$ (the \textit{associativity constraint\/}) and
$l=l_{U}\colon \tunit\otimes U\to U$, $r=r_{U}\colon U\otimes\tunit\to U$,
for any $U\in\mathcal{C}$ (the \textit{left unit constraint} and the 
\textit{right unit constraint,} respectively) such that, for all 
$U, V, W, X\in\mathcal{C}$, the two identities
$a_{U,V,W\otimes X}\circ a_{U\otimes V,W,X}= 
(U\otimes a_{V,W,X}) \circ a_{U,V\otimes W,X} \circ (a_{U,V,W}\otimes X)$
(called the \textit{associativity pentagon}) and
$(U\otimes l_{V}) \circ (r_{U}\otimes V) = a_{U,\tunit,V}$
are satisfied.
A tensor category $\mathcal{C}$ is \textit{strict} when all the 
constraints are identities.

Given two tensor categories $\mathcal{C}$ and $\mathcal{D}$, a 
\textit{tensor functor} $F=(F,F_{2},F_{0})\colon\mathcal{C}\to\mathcal{D}$ 
consists of the following items.
\begin{itemize}
 \item A functor $F\colon \mathcal{C}\to\mathcal{D}$.
 \item A natural family $\bigl\{F_{2}(U,V)\colon F(U)\otimes F(V)\to F(U\otimes V)\bigr\}_{U,V\in\mathcal{C}}$ of isomorphisms in $\mathcal{D}$
       such that
       $F(a_{U,V,W}) \circ F_{2}(U\otimes V, W) \circ F_{2}(U,V)\otimes F(W)
        =F_{2}(U,V\otimes W) \circ \bigl(F(U)\otimes F_{2}(V,W)\bigr) \circ a_{U,V,W}$, 
        for any $U, V, W\in\mathcal{C}$.
 \item An isomorphism $F_{0}\colon\tunit\to F(\tunit)$ in $\mathcal{D}$
      such that
      $F(r_{U}) \circ F_{2}(U,\tunit) \circ \bigl(F(U)\otimes F_{0}\bigr) =
       r_{F(U)}$ and
       $F(l_{U}) \circ F_{2}(\tunit,U) \circ \bigl(F_{0}\otimes F(U)\bigr) =
       l_{F(U)}$, for any $U\in\mathcal{C}$.
\end{itemize}
\noindent $F$ is said \textit{strict} when $F_{0}$ 
and all the $F_{2}(U,V)$ are identities.

\begin{rmk}\label{rmk:strict}
 Let $\mathcal{C}$ be a tensor category. We recall~\cite{ML-coherence,ML}
 that $\mathcal{C}$ is equivalent to a strict tensor category 
 $\mathcal{S}(\mathcal{C})$ via a tensor functor 
 $F\colon\mathcal{S}(\mathcal{C})\to\mathcal{C}$ and a tensor functor
 $G\colon\mathcal{C}\to\mathcal{S}(\mathcal{C})$. More precisely, the category
 $\mathcal{S}(\mathcal{C})$ and the functors $F$ and $G$ can be obtained
 as follows.
 \begin{itemize}
  \item The objects of $\mathcal{S}(\mathcal{C})$ are all the finite 
        sequences $u=(U_{1},\ldots, U_{n})$ of objects 
        $U_{1},..., U_{n}\in\mathcal{C}$. Also
        the empty sequence, denoted $u_{0}$, is an object in
        $\mathcal{S}(\mathcal{C})$.
  \item For any $u\in\mathcal{S}(\mathcal{C})$, the object $F(u)$ is given by
        \begin{equation*} 
          F(u) = \begin{cases}
           \tunit & \text{if $u=u_{0}$}\text{,}\\
           \Bigl(\cdots \bigl((U_{1}\otimes U_{2})\otimes U_{3}\bigr)\otimes\cdots\Bigr)\otimes U_{n} &
            \text{if $u=(U_{1},\ldots,U_{n})$, with $n\in\N\setminus\{0\}$}\text{,}
        \end{cases}\end{equation*}
        where on the right all pairs of parenthesis begin if front.
        For any $u,v\in\mathcal{S}(\mathcal{C})$, the arrows from $u$ to $v$
        in $\mathcal{C}$ are given by
        $
          \mathcal{S}(\mathcal{C})(u,v)\eqdef\mathcal{C}
          \bigl(F(v),F(v)\bigr)$.
        In that way, with the composition induced by $\mathcal{C}$, we 
        obtain the category $\mathcal{S}(\mathcal{C})$.
  \item $\mathcal{S}(\mathcal{C})$ becomes a tensor category with the
        tensor product of objects given by the concatenation product
        and the tensor product of two arrows
        $f\in\mathcal{S}(\mathcal{C})(u,v)$ and
        $g\in\mathcal{S}(\mathcal{C})(u',v')$ given by the composite
        $F_{2}(v,v') \circ (f \otimes g) \circ F_{2}^{-1}(u,u')$,
        where, for any $w,w'\in\mathcal{S}(\mathcal{C})$, the arrow 
        $F_{2}(w,w')$ is the canonical isomorphism in $\mathcal{C}$
        from $F(w)\otimes F(w')$ to $F(w\otimes w')$ obtained iterating the
        associativity constraint $a$
        (well defined by the coherence theorem in~\cite{ML-coherence}).
  \item The definition of the functor $F$ is completed by setting
        $F(f)=f$ for any arrow 
        $f\in\mathcal{S}(\mathcal{C})(u,v)=\mathcal{C}\bigl(F(u),F(v)\bigr)$,
        with $u,v\in\mathcal{S}(\mathcal{C})$.
        $F$ becomes a tensor functor by defining $F_{2}(\cdot,\cdot)$ as above and 
        $F_{0}=\Id_{\tunit}$.
  \item The category $\mathcal{C}$ can be embedded in
        $\mathcal{S}(\mathcal{C})$ by identifying $\mathcal{C}$ with the
        full subcategory of $\mathcal{S}(\mathcal{C})$ given by the
        sequences of length one. The functor $G$ is given by the immersion
        of $\mathcal{C}$ in $\mathcal{S}(\mathcal{C})$. $G$ becomes a
        tensor functor by setting
        $G_{2}(U,V)=a_{U,V}$ for any $U,V\in\mathcal{C}$, and
        $G_{0}=\Id_{\tunit}\in\mathcal{S}(\mathcal{C})(u_{0},\tunit)=
          \mathcal{C}(\tunit,\tunit)$.
 \end{itemize}
\end{rmk}

\paragraph{\scshape Dualities}
Let $\mathcal{C}$ be a tensor category. For simplicity, allowed by
Remark~\ref{rmk:strict}, we suppose that $\mathcal{C}$ is strict.
Given $U,V\in\mathcal{C}$, a \textit{pairing between $V$ and $U$} is 
an arrow $d\colon V\otimes U\longrightarrow\tunit\text{.}$ in $\mathcal{C}$.
If, for any arrow $f\colon X\to U\otimes Y$
in $\mathcal{C}$, we set
$d^{\sharp}(f)=\biggl(V\otimes X\xrightarrow{V\otimes f} V\otimes U\otimes Y
  \xrightarrow{d\otimes Y} Y\biggr)$,
then we obtain an application
$d^{\sharp}\colon\mathcal{C}(X,U\otimes Y)\longrightarrow\mathcal{C}(V\otimes X, Y)$.
The pairing $d$ is \textit{exact}{\label{p:exact}} when $d^{\sharp}$ is
bijective for any $X,Y\in\mathcal{C}$, i.e., if we have an adjunction
of functors $\adjunction{V\otimes\_}{U\otimes\_}$. 
It follows that $d$ is exact  if and only if there exists a map
$b\colon\tunit\to U \otimes V$ (that is, $b=\bigl(d^{\sharp}\bigr)^{-1}(\Id_{V})\/$)  such that the relations
(called \textit{adjunction triangles\/} or \textit{duality relations\/})
$(U \otimes d) \circ (b \otimes U)= U$ and $(d \otimes V) \circ (V \otimes b) = V$
hold. When the pairing is exact, we say that the pair $(b,d)$ is an 
\textit{adjunction} or a \textit{duality} between $V$ and $U$. 
We also say that $V$ is \textit{left adjoint} or \textit{left dual} 
to $U$, that $U$ is \textit{right adjoint} or \textit{right dual} to 
$V$, and we write $(b,d)\colon\adjunction{V}{U}$.
We call $b$ the \textit{unit} and $d$ the \textit{counit} of the
adjunction.
 
Given two adjunction $(b_{1},d_{1})\colon\adjunction{V_{1}}{U_{1}}$ and 
$(b_{2},d_{2})\colon\adjunction{V_{2}}{U_{2}}$ in $\mathcal{C}$, 
we have a bijection 
$(\,\hat{ }\,)\colon\mathcal{C}(V_{1},V_{2})\to\mathcal{C}(U_{2},U_{1})$
with inverse 
$(\,\check{ }\,)\colon\mathcal{C}(U_{2},U_{1})\to\mathcal{C}(V_{1},V_{2})$,
obtained by setting, for any $f\in\mathcal{C}(V_{1},V_{2})$ and 
$g\in\mathcal{C}(U_{2},U_{1})$,
$\hat{g} = (V_{2}\otimes b_{1})\circ(V_{2}\otimes g\otimes V_{1})\circ(d_{2}\otimes V_{1})$, and
$\check{f}  = (b_{2}\otimes U_{1})\circ(U_{2}\otimes f\otimes U_{1})\circ(U_{2}\otimes d_{1})$.
When $g=\check{f}$ we write
$
  \adjunction{f}{g}
$.
We say that $\mathcal{C}$ is \textit{left} 
(respectively, \textit{right\/}) \textit{autonomous} when any object has
a left (respectively, right) dual. We say that $\mathcal{C}$ is 
\textit{autonomous} if it is both left and right autonomous.
 
When $\mathcal{C}$ is left autonomous, choosen an adjunction 
$(b_{V},d_{V})\colon\adjunction{V^{\ast}}{V}$ for any $V\in\mathcal{C}$, we get 
a functor $(\cdot)^{\ast}\colon\mathcal{C}\to\mathcal{C}^{\ast}$,
defined on $f\in\mathcal{C}(V,U)$ by the condition
$\adjunction{f^{\ast}}{f}\text{.}$.

Fix two adjunctions
$(b_{U},d_{U})\colon\adjunction{U^{\ast}}{U}$ and
$(b_{V},d_{V})\colon\adjunction{V^{\ast}}{V}$
in a tensor category $\mathcal{C}$.
Given $f\in\mathcal{C}(X\otimes U,V\otimes Y)$, the \textit{mate $f^{\text{@}}$ of $f$} 
is the arrow
\begin{multline}\label{e:mate}
  f^{\text{@}}=\biggl(V^{\ast}\otimes X\xrightarrow{V^{\ast}\otimes X\otimes b_{U}}
  V^{\ast}\otimes X\otimes U\otimes U^{\ast}\xrightarrow{V^{\ast}\otimes f\otimes U^{\ast}} V^{\ast}\otimes V\otimes Y\otimes U^{\ast}\\
  V^{\ast}\otimes V\otimes Y\otimes U^{\ast}\xrightarrow{d_{V}\otimes Y\otimes U^{\ast}}Y\otimes U^{\ast}\biggr)\text{.}
\end{multline}

\section{\protect\Tcats}\label{sec:cpc}
Let $\pi$ be a group.
A \textit{\Tcat $\mathcal{T}$ \textup{(}over $\pi\/$\textup{)}} 
is given by the following data.
\begin{itemize}
 \item A tensor category $\mathcal{T}$.
 \item A family of subcategories $\{\mathcal{T}_{\alpha}\}_{\alpha\in\pi}$
       such that $\mathcal{T}$ is disjoint union of this 
       family and that
       $U\otimes V\in\mathcal{T}_{\alpha\beta}$, for any $\alpha,\beta\in\pi$, $U\in\mathcal{T}_{\alpha}$, and 
         $V\in\mathcal{T}_{\beta}$.
 \item Denoted $\aut(\mathcal{T})$ the group of the invertible strict 
       tensor functors from $\mathcal{T}$ to itself, a group homomorphism
       $\map{\varphi}{\pi}{\aut(\mathcal{T})}{\beta}{\varphi_{\beta}}$,
       the \textit{conjugation,} such that 
       $\varphi_{\beta}(\mathcal{T}_{\alpha})=\mathcal{T}_{\beta\alpha\beta^{-1}}$ for any $\alpha,\beta\in\pi$.
\end{itemize}
In the terminology of~\cite{Tur-CPC}, a \Tcat is called a 
\textit{crossed group category.} Differently 
from~\cite{Tur-CPC}, we do not require that a \Tcat is a linear category.
Given $\alpha\in\pi$, the subcategory $\mathcal{T}_{\alpha}$ is called the 
\textit{\alphath component} of $\mathcal{T}$ while the functors $\varphi_{\beta}$ 
are called \textit{conjugation isomorphisms.} 
$\mathcal{T}$ is called \textit{strict} when it is strict 
as a tensor category. When $\pi=1$, 
then $\mathcal{T}$ is nothing but a tensor category.

Given two \Tcats $\mathcal{T}$ and $\mathcal{T}'$, 
a \textit{\Tfun}$F\colon\mathcal{T}\to\mathcal{T}'$ is a tensor functor from
$\mathcal{T}$ to $\mathcal{T}'$ that satisfies the following 
two conditions.
\begin{enumerate}
 \item For any $\alpha\in\pi$, $F(\mathcal{T}_{\alpha})\subset\mathcal{T}'_{\alpha}$.
 \item $F$ commutes with the conjugation isomorphisms.
\end{enumerate}
Two \Tcats $\mathcal{T}$ and 
$\mathcal{T}'$ are \textit{equivalent as \Tcats} if they are equivalent 
as categories via a \Tfun $F\colon\mathcal{T}\to\mathcal{T}'$ and a 
\Tfun $G\colon\mathcal{T}\to\mathcal{T}'$.
 
\paragraph{\scshape Left index notation}
Given $\beta\in\pi$ and an object 
$V\in\mathcal{T}_{\beta}$, the functor $\varphi_{\beta}$ will be denoted $\lidx{V}{(\cdot)}$, 
as in~\cite{Tur-CPC}, or also $\lidx{\beta}{\,(\cdot)}$. We introduce the 
notation $\lidx{\overline{V}}{(\cdot)}$ for $\lidx{\beta^{-1}}{(\cdot)}$.
Since $\lidx{V}{(\cdot)}$ is a functor, for any object $U\in\mathcal{T}$ and 
for any couple of composable arrows
$\cdot\,\xrightarrow{\ f \ }\,\cdot\,\xrightarrow{\ g \ }\,\cdot$ in $\mathcal{T}$,
we obtain
$\leftidx{^{V}}{\!\Id}{_{U}}=\Id_{\lidx{V}{U}}$ and
$\lidx{V}{(g\circ f)}=\lidx{V}{g}\circ\lidx{V}{f}$.
Since the conjugation $\varphi\colon\pi\to\aut(\mathcal{T})$ is a group homomorphism, 
for any $V,W\in\mathcal{T}$, we have
$\lidx{V\otimes W}{(\cdot)}=\lidx{V}{\Bigl(\lidx{W}{(\cdot)}\Bigr)}$
and
$\lidx{\tunit}{\,(\cdot)} = \lidx{V}{\Bigl(\lidx{\overline{V}}{(\cdot)}\Bigr)}=
   \lidx{\overline{V}}{\Bigl(\lidx{V}{(\cdot)}\Bigr)}= 
   \Id_{\mathcal{T}}$.
Since, for any $V\in\mathcal{C}$, the functor $\lidx{V}{(\cdot)}$ is
strict, we have 
$\lidx{V}{(f\otimes g)} =\lidx{V}{f}\otimes\lidx{V}{g}$,
for any arrow $f$ and $g$ in $\mathcal{T}$, and
$\lidx{V}{\,\tunit} =  \tunit$.
  
\paragraph{\scshape Strict equivalence for \protect\Tcats}
 
\begin{thm}\label{lemma:strict}
 Let $\mathcal{T}$ be a \Tcat. $\mathcal{T}$ is equivalent as a 
 \Tcat to a strict \Tcat $\mathcal{S}(\mathcal{T})$.
\end{thm}
 
\begin{proof}[Proof \textup{(sketch)}]
 Define the category $\mathcal{S}(\mathcal{T})$ and the functors $F$ 
 and $G$ as in Remark~\ref{rmk:strict}.  
 \begin{itemize}
  \item Let $u=(U_{1},...,U_{n})$ be in $\mathcal{S}(\mathcal{T}_{\alpha})$,
        with $n\geq 1$ and let $U_{1}\in\mathcal{T}_{\alpha_{1}},
        U_{2}\in\mathcal{T}_{\alpha_{2}}, \ldots, U_{n}\in\mathcal{T}_{\alpha_{n}}$. We set
        $m(u)=\alpha_{1}\alpha_{2}\cdots\alpha_{n}$, and $m(u_{0})=1$,
        where $u_{0}$ is the empty sequence. Any \alphath 
        component of $\mathcal{S}(\mathcal{T})$ is defined as the full 
        subcategory $\mathcal{S}_{\alpha}(\mathcal{T})$ 
        whose objects are the objects $u$ of 
        $\mathcal{S}(\mathcal{T})$ such that $m(u)=\alpha$.
  \item The conjugation $\varphi^{\mathrm{str}}$ of $\mathcal{S}(\mathcal{T})$ is
        obtained by setting, for any $\alpha\in\pi$,
        $\varphi^{\mathrm{str}}_{\alpha}(u)=\varphi^{\mathrm{str}}_{\alpha}(U_{1},\ldots,U_{n})=
          \bigl(\varphi_{\alpha}(U_{1}),\ldots,\varphi_{\alpha}(U_{n})\bigr)$,
        for any $u=(U_{1},\ldots,U_{n})\in\mathcal{S}(\mathcal{T})$, and
        $\varphi_{\alpha}^{\mathrm{str}}(u_{0})=u_{0}$.
        The definition is completed by setting
        $\varphi^{\mathrm{str}}_{\alpha}(f)=\varphi_{\alpha}(f)$, for any arrow
        $f\in\mathcal{S}(\mathcal{T})$.
 \end{itemize}

 It is easy to prove that, in that way, $\mathcal{S}(\mathcal{T})$ 
 becomes a \Tcat and $F$ and $G$ become \Tfuns.
 Notice that the hypothesis that the functor $\varphi_{\alpha}$ ($\alpha\in\pi\/$) is strict 
 is essential to obtain the functor $\varphi^{\mathrm{str}}_{\alpha}$.
\end{proof}
 
In virtue of Theorem~\ref{lemma:strict}, often we will
consider only strict \Tcats.
 
\paragraph{\scshape Adjunctions in a \protect\Tcat}
A \textit{left autonomous \Tcat}$\mathcal{T}=(\mathcal{T},(\cdot)^{\ast})$ is a 
\Tcat $\mathcal{T}$ endowed with a choice of left dualities $(\cdot)^{\ast}$
satisfying the following two conditions.
\begin{itemize}
 \item If $U$ is an object in $\mathcal{T}_{\alpha}$ (with $\alpha\in\pi\/$), then 
       $U^{\ast}$ is an object in $\mathcal{T}_{\alpha^{-1}}$.
 \item The conjugation preserve the chosen dualities,
       i.e., $\varphi_{\beta}(b_{U}) = b_{\varphi_{\beta}(U)}$
       and $\varphi_{\beta}(d_{U}) = d_{\varphi_{\beta}(U)}$ for any $\beta\in\pi$ and $U\in\mathcal{T}$.
\end{itemize}
\noindent Similarly, it possible to introduce the notion of a 
right autonomous \Tcat. An \textit{autonomous \Tcat}\ is a \Tcat that 
is both left and right autonomous.
  
Given two left autonomous \Tcats 
$\mathcal{T}$ and $\mathcal{T}'$, a \textit{left autonomous \Tfun 
$F\colon\mathcal{T}\to\mathcal{T}'$} is a \Tfun from $\mathcal{T}$ to 
$\mathcal{T}'$ that preserves the dualities, that is $F(b_{U})=b_{F(U)}$
and $F(d_{U})=d_{F(U)}$, for any $U\in\mathcal{T}$.
Two left autonomous \Tcats $\mathcal{T}$ and $\mathcal{T}'$ are 
\textit{equivalent as left autonomous \Tcats}if they are equivalent 
as categories via a left autonomous \Tfun $F\colon\mathcal{T}\to\mathcal{T}'$ 
and a left autonomous \Tfun $G\colon\mathcal{T}\to\mathcal{T}'$.
Similarly, it is possible to introduce the notions of a right 
autonomous \Tfun and of an autonomous \Tfun and the notions of equivalence
of right autonomous \Tcats and of autonomous \Tcats.
 
\begin{rmk}\label{rmk:Prometeus}
 Let $\mathcal{T}$ be a left autonomous T-category. Define the \Tcat 
 $\mathcal{S}(\mathcal{T})$ and the \Tfuns $F$ and $G$ as in 
 Theorem~\ref{lemma:strict}. Given $u\in\mathcal{S}(\mathcal{T})$, if we set
 $ 
   u^{\ast}=G\bigl(F(u)^{\ast}\bigr)\text{,}
 $
 then the exact pairing
 $
   F(u^{\ast}\otimes u)\xrightarrow{\,\cong\,} F(u)^{\ast}\otimes F(u)\xrightarrow{\,d_{F(U)}\,}\tunit
 $
 in $\mathcal{T}$ gives also an exact pairing $u^{\ast}\otimes u\to\tunit$ under the identification
 $
   \mathcal{S}(\mathcal{T})(u^{\ast}\otimes u,u_{0})=
   \mathcal{T}\bigl(F(u^{\ast}\otimes u),\tunit\bigr)\text{.}
$
 It is easy to check that $\mathcal{S}(\mathcal{T})$
 becomes a left autonomous \Tcat and that $\mathcal{T}$ 
 is equivalent to $\mathcal{S}(\mathcal{T})$ as a left autonomous \Tcat 
 via $F$ and $G$.
\end{rmk}
  
\paragraph{\scshape Stable left duals}
Let{\label{par:stable-dual}} 
$\mathcal{T}$ be a \Tcat and let $U\in\mathcal{T}_{\alpha}$ (with $\alpha\in\pi\/$) be an object 
endowed with an adjunction $(b_{U},d_{U})\colon\adjunction{U^{\ast}}{U}$. We say 
that $(b_{U},d_{U})$ is a \textit{stable adjunction} and $U^{\ast}$ a 
\textit{stable left dual of $U$} when, 
for any $\beta\in\pi$ that commutes with $\alpha$,   
if $\varphi_{\beta}(U)=U$ then $\bigl(\varphi_{\beta}(b_{U}),\varphi_{\beta}(d_{U})\bigr)=(b_{U},d_{U})$.
If we set
$
  \Phi(U)\eqdef\bigl\{\varphi_{\beta}(U)\bigr\}_{\beta\in\pi}\text{,}
$
then, given $V\in\Phi(U)$ and $\beta\in\pi$ such that $V=\varphi_{\beta}(U)$, the stable 
adjunction $(b_{U},d_{U})$ gives rise to another stable adjunction 
$\bigl(\varphi_{\beta}(b_{U}),\varphi_{\beta}(d_{U})\bigr)\colon\adjunction{\varphi_{\beta}(U^{\ast})}{V}$
that does not depends on $\beta$.

\begin{lemma}\label{l:aut-op}
 A \Tcat $\mathcal{T}$ admits a structure of left autonomous \Tcat if
 and only if, for any $U\in\mathcal{T}$, there exists an object $U_{0}\in\Phi(U)$
 endowed with a stable adjunction $(b_{0},d_{0})\colon\adjunction{U^{\ast}_{0}}{U_{0}}$.
\end{lemma}
 
\begin{rmk}\label{rmk:choice}
 The terminology concerning a category with dualities is not completely 
 fixed. In particular, if some authors~\cite{JS0} only require that any
 object $V$ in a left autonomous category admits an exact pairing, other
 authors~\cite{Kas,Tur-QG} also require the choice of a pairing, i.e.,
 they consider a fixed adjunction for any object of the category.
 To be coherent with the definition in~\cite{Tur-CPC}, 
 we choose the second
 convention. This will also be useful in the second part of the article,
 since the considered categories will be endowed with 
 a natural
 choice of stable dualities. However, we will see that, starting
 from a \Tcat $\mathcal{T}$ endowed with a twist $\theta$, it is possible to
 obtain a ribbon subcategory $\mathcal{N}(\mathcal{T})$ of $\mathcal{T}$, 
 i.e., a subcategory of $\mathcal{T}$ endowed with stable dualities 
 compatible with the twist. With the exception of the trivial case in
 which we just know that $\mathcal{T}$ is ribbon (so that we have
 $\mathcal{N}(\mathcal{T})=\mathcal{T}\/$), there is no natural way to
 obtain a canonical duality $(\cdot)^{\ast}$ for $\mathcal{N}(\mathcal{T})$.
\end{rmk}

\paragraph{\scshape Braiding}
A \textit{braiding} for a \Tcat $\mathcal{T}$ is
a family of isomorphisms
\begin{equation*}
  c=\Biggl\{c_{U, V}\in\mathcal{T}\biggl(U\otimes V,
  \Bigl(\lidx{U}{V}\Bigr)\otimes U\biggr)\Biggr\}_{U,V\in\mathcal{T}}
\end{equation*}
satisfying the following conditions.
\begin{subequations}\label{e:braiding}
\begin{itemize}
 \item For any arrow $f\in\mathcal{T}_{\alpha}(U,U')$ (with $\alpha\in\pi\/$),
       $g\in\mathcal{T}(V,V')$ we have
       \begin{equation}\label{e:braiding-a}
        \Bigl(\bigl(\lidx{\alpha}{g}\bigr)\otimes f\Bigr) \circ c_{U,V} 
        = c_{U',V'} \circ (f \otimes g).
      \end{equation}
 \item For any $U,V,W\in\mathcal{T}$, we have
       \begin{equation}\label{e:braiding-b}
        c_{U\otimes V, W} = a_{\lidx{U\otimes V}{W},U,V} \circ (c_{U,\lidx{V}{W}}\otimes V) \circ 
        a^{-1}_{U,\lidx{V}{W},V} \circ (U\otimes c_{V,W})
       \end{equation}
       \noindent and
       \begin{equation}\label{e:braiding-c}
       c_{U, V\otimes W} = a^{-1}_{\lidx{U\otimes V}{W},U,V} \circ 
                 \Bigl(\bigl(\lidx{U}{V}\bigr)\otimes c_{U,W}\Bigr) \circ 
                 a_{\lidx{U}{V},U,W} \circ (c_{U,V}\otimes W) \circ a^{-1}_{U,V,W}\text{.}
       \end{equation}
   \item For any $U,V\in\mathcal{T}$ and $\beta\in\pi$, we have
         \begin{equation}\label{e:braiding-d}
            \varphi_{\beta}(c_{U,V}) = c_{\varphi_{\beta}(U),\varphi_{\beta}(V)}\text{.}
         \end{equation}
  \end{itemize}
 \end{subequations}
 
\noindent A \Tcat endowed with a braiding is called a 
\textit{braided \Tcat.} In particular,  when $\pi=1$, we recover 
the usual definition of a braided tensor 
category~\cite{JS}.

Given two braided \Tcat $\mathcal{T}$ and $\mathcal{T}'$, a 
\textit{braided \Tfun}$F\colon\mathcal{T}\to\mathcal{T}'$ is \Tfun from 
$\mathcal{T}$ to $\mathcal{T}'$ that preserves the braiding, 
i.e., such that $F(c_{U,V})=c_{F(U),F(V)}$ for any $U,V\in\mathcal{T}$.
Two braided \Tcats $\mathcal{T}$ and $\mathcal{T}'$ are 
\textit{equivalent as braided \Tcats} if they are equivalent as \Tcats 
via a braided \Tfun $F\colon\mathcal{T}\to\mathcal{T}'$ and a braided \Tfun 
$G\colon\mathcal{T}'\to\mathcal{T}$.
 
\begin{rmk}\label{rmk:rol}
 Let $\mathcal{T}$ be a braided \Tcat  and define 
 $\mathcal{S}(\mathcal{T})$,
 $F$, and $G$ as in Theorem~\ref{lemma:strict}. The family of
 arrows
 \begin{equation*}
   c_{u,v}=\Biggl(F(u\otimes v)\xrightarrow{\,\cong\,}
          F(u)\otimes F(v)\xrightarrow{\,c_{F(u),F(v)}\,}
          \Bigl(\lidx{m(u)}{F(v)}\Bigr)\otimes F(u)\xrightarrow{\,\cong\,}
          F\biggl(\Bigl(\lidx{u}{\,v}\Bigr)\otimes u\biggr)\Biggr)
 \end{equation*}
 (for any $u,v\in\mathcal{S}(\mathcal{T})\/$) is a braiding in 
 $\mathcal{S}(\mathcal{T})$. 
 With this structure of braided \Tcat on $\mathcal{S}(\mathcal{T})$, the 
 functors $F$ and $G$ become braided \Tfuns and so $\mathcal{T}$ is 
 equivalent to $\mathcal{S}(\mathcal{T})$ as a braided \Tcat.
\end{rmk}

\paragraph{\scshape The morphism $\omega$} Let $U$ be an object in a braided \Tcat $\mathcal{T}$ endowed with
a left dual $U^{\ast}$ via an adjunction $(b_{U},d_{U})$. We set
\begin{equation}\label{e:omega}
 \omega_{U}=(d_{\leftidx{^{U\otimes U}}{\! U}{^{\ast}}}\otimes U)\circ
 \biggl(\Bigl(\leftidx{^{U\otimes U}}{\! U}{^{\ast}}\Bigr)\otimes
 \tilde{c}_{\lidx{U}{U}, \lidx{U\otimes U}{U}}\biggr)\circ
 \Bigl((c_{\,\lidx{U}{U},\leftidx{^{U}}{\! U}{^{\ast}}}\circ
  b_{\lidx{U}{U}})\otimes\lidx{U\otimes U}{U}\Bigr)\text{.}
\end{equation}

\begin{lemma}\label{l:omega-invar}
  $\omega_{U}$ is independent from the choice of the stable left
  adjunction of $U$.
\end{lemma}

Since the proof is relatively long, it is omitted.

\paragraph{\scshape Twist}
A \textit{twist} for a braided \Tcat $\mathcal{T}$ 
is a family of isomorphisms
\begin{equation*}
  \theta=\Bigl\{\theta_{U}\colon U\to\lidx{U}{U}\Bigr\}_{U\in\mathcal{T}}
\end{equation*}
satisfying the following conditions.
\begin{subequations}\label{e:twist}
\begin{itemize}
 \item $\theta$ is \textit{natural}, i.e., for any $f\in\mathcal{T}_{\alpha}(U,V)$ 
       (with $\alpha\in\pi\/$),
       \begin{equation}\label{e:twist-natural}
        \theta_V \circ f = \bigl(\lidx{\alpha}{f}\bigr) \circ \theta_{U}\text{.}
        \end{equation}
 \item For any $U\in\mathcal{T}_{\alpha}$ and $V\in\mathcal{T}_{\beta}$ 
       (with $\alpha,\beta\in\pi\/$),
       \begin{equation}\label{e:twist-main}
        \theta_{U\otimes V}=c_{\lidx{U\otimes V}{V}, \lidx{U}{U}} 
        \circ c_{\lidx{U}{U},\lidx{V}{V}} \circ (\theta_{U}\otimes\theta_{V})\text{.}
       \end{equation}
  \item For any $U\in\mathcal{T}$ and $\alpha\in\pi$,
        \begin{equation}\label{e:twist-ultra}
         \varphi_{\alpha}(\theta_{U}) = \theta_{\varphi_{\alpha}(U)}\text{.}
        \end{equation}
\end{itemize}
\end{subequations}
A braided \Tcat endowed with a twist
is called a \textit{balanced \Tcat.} In particular, for $\pi=1$
we recover the usual definition of a balanced tensor
category~\cite{JS}.\index{tcat@\Tcat!balanced|)}
 
Given two balanced \Tcats $\mathcal{T}$ and 
$\mathcal{T}'$, a braided \Tfun $F\colon\mathcal{T}\to\mathcal{T}'$ is a 
\textit{balanced \Tfun}if it preserves the twist, i.e., if
$F(\theta_{U})=\theta_{F(U)}$ for any $U\in\mathcal{T}$.
Two \Tcats $\mathcal{T}$ and 
$\mathcal{T}'$ are \textit{equivalent as balanced \Tcats}if they are 
equivalent as \Tcat via balanced \Tfun  $F\colon\mathcal{T}\to\mathcal{T}'$ 
and a balanced \Tfun $G\colon\mathcal{T}'\to\mathcal{T}$.
 
A \textit{ribbon \Tcat}$\mathcal{T}$ is a balanced \Tcat that is also 
a left autonomous \Tcat such that for any 
$U\in\mathcal{T}_{\alpha}$ (with $\alpha\in\pi\/$), 
\begin{equation}\label{e:tortility} 
\Bigl(\bigl(\lidx{U}{U}\bigr)\otimes\theta_{\leftidx{^{U}}{\! U}{^{\ast}}}\Bigr) \circ 
 b_{\lidx{U}{U}}=
 (\theta_{U} \otimes U^{\ast}) \circ b_{U}\text{.}
\end{equation}
For $\pi=1$ we recover the usual definition of a ribbon 
category~\cite{RT,Tur-QG} also called 
\textit{tortile tensor category}~\cite{JS,JS-tortile,Shum}.
 
Given two ribbon \Tcats $\mathcal{T}$ and 
$\mathcal{T}'$, a \Tfun $F\colon\mathcal{T}\to\mathcal{T}'$ that is at the same
time a balanced \Tfun and a left autonomous \Tfun is called 
\textit{ribbon \Tfun.}
Two ribbon \Tcats 
$\mathcal{T}$ and $\mathcal{T}'$ are 
\textit{equivalent as ribbon \Tcats}if they are equivalent as
\Tcats via a ribbon \Tfun $F\colon\mathcal{T}\to\mathcal{T}'$ and a ribbon
\Tfun $G\colon\mathcal{T}'\to\mathcal{T}$.
 
\begin{rmk}\label{ref:sinus}
 Let $\mathcal{T}$ be a balanced \Tcat and define the equivalent braided 
 \Tcat $\mathcal{S}(\mathcal{T})$ and the functors $F$ and $G$ as in 
 Theorem~\ref{lemma:strict} and Remark~\ref{rmk:rol}. The family of 
 arrows $\theta_{u}=\theta_{F(u)}$
 (with $u\in\mathcal{S}(\mathcal{T})\/$) is a twist in 
 $\mathcal{S}(\mathcal{T})$ such that $\mathcal{T}$ is equivalent to 
 $\mathcal{S}(\mathcal{T})$ as a balanced \Tcat via $F$
 and $G$. If $\mathcal{T}$ is ribbon, then, with the 
 structure of left autonomous \Tcat on $\mathcal{S}(\mathcal{T})$ 
 provided in Remark~\ref{rmk:Prometeus}, $\mathcal{T}$ is equivalent to 
 $\mathcal{S}(\mathcal{T})$ as a ribbon \Tcat.
\end{rmk}
 
\paragraph{\scshape Mirror \protect\Tcat}\label{par:mirror-cat}
Let $\mathcal{T}$ be \Tcat. The 
\textit{mirror $\overline{\mathcal{T}}$ of a \Tcat $\mathcal{T}$}
(see~\cite{Tur-CPC} and the description of the corresponding
algebraic notion of mirror \Tcoalg below) is the \Tcat defined as follows.
\begin{itemize}
 \item For any $\alpha\in\pi$, we set 
       $\overline{\mathcal{T}}_{\!\!\alpha}=\mathcal{T}_{\alpha^{-1}}$ as a category. 
       So, as a category, $\overline{\mathcal{T}}=\mathcal{T}$.
 \item The tensor product $U\overline{\otimes}V$ in $\overline{\mathcal{T}}$ of 
       $U\in\overline{\mathcal{T}}_{\alpha}=\mathcal{T}_{\alpha^{-1}}$
       and $V\in\overline{\mathcal{T}}_{\beta}=\mathcal{T}_{\beta^{-1}}$ (with $\alpha,\beta\in\pi\/$)
       is given by
       $ U\overline{\otimes}V=\varphi_{\beta^{-1}}(U)\otimes V\in\overline{\mathcal{T}}_{\!\!\alpha\beta}
        $.
       Given an arrow $f$ in $\overline{\mathcal{T}}_{\alpha}$ and an arrow 
       $g$ in $\overline{\mathcal{T}}_{\beta}$ (with $\alpha,\beta\in\pi\/$), the tensor 
       product $f\overline{\otimes}g$ of $f$ and $g$ in $\overline{\mathcal{T}}$
       is given by
       $
          f\overline{\otimes}g=\varphi_{\beta^{-1}}(f)\otimes g
       $.
 \item The associativity constraint is given by
       $\overline{a}_{U,V,W}=a_{\varphi_{\gamma^{-1}\beta^{-1}}(U),\varphi_{\gamma^{-1}}(V),W}$,
       for any
       $U\in\mathcal{T}_{\alpha}$, $V\in\mathcal{T}_{\beta}$, and $W\in\mathcal{T}_{\gamma}$ 
       (with $\alpha,\beta,\gamma\in\pi\/$).
 \item The left unit constraint and the right unit constraint 
       are given by $\overline{l}_{U}=l_{U}$ and
       $\overline{r}_{U}=r_{U}$, for any $U\in\overline{\mathcal{T}}$.
 \item The conjugation is given by $\overline{\varphi}_{\alpha}=\varphi_{\alpha}$.
\end{itemize}
When $\mathcal{T}$ is a left autonomous, $\overline{\mathcal{T}}$ 
is left autonomous by setting, for any $U\in\overline{\mathcal{T}}_{\!\!\alpha}$ 
(with $\alpha\in\pi\/$), $\overline{b}_{U}=\varphi_{\alpha}(b_{U})$ and 
$\overline{d}_{U}=d_{U}$.
When $\mathcal{T}$ is braided, $\overline{\mathcal{T}}$ is braided 
via $\overline{c}_{U,V}=(c_{V,U})^{-1}$,
for any $U,V\in\mathcal{T}$.
When $\mathcal{T}$ is balanced (respectively, ribbon),  
$\overline{\mathcal{T}}$ is also balanced (respectively, ribbon)
via $\overline{\theta}_{U}=(\theta_{\varphi_{\alpha}(U)}\bigr)^{-1}$,
for any $U\in\overline{\mathcal{T}}_{\!\!\alpha}$ (with $\alpha\in\pi\/$).
The mirror construction is involutive, i.e., we have 
$\overline{\overline{\mathcal{T}}}=\mathcal{T}$.
  
Two [(left/right) autonomous, braided, balanced or ribbon] 
\Tcats $\mathcal{T}$ 
and $\mathcal{T}'$ are \textit{mirror equivalent 
\textup{[\/}as \textup{(\/}left/right\textup{)} autonomous, braided, balanced \textup{or} 
ribbon \Tcats\/\textup{]}} if $\mathcal{T}$ is equivalent to 
$\overline{\mathcal{T}'}$ as [(left/right) autonomous, braided, balanced
or ribbon] \Tcat.
 
\section{\protect\Tcats of representations}\label{ex:main-catex}

\paragraph{\scshape\protect\Tcoalgs}
 Let $\Bbbk$ be a commutative field and let $\pi$ be a group. 
 A \textit{\Tcoalg $H$ \textup{(}over $\pi$ 
 and $\Bbbk$\/\textup{)}}~\cite{Tur-CPC}
 is given by the following data.
 \begin{itemize}
  \item For any $\alpha\in\pi$, an associative \kalg $H_{\alpha}$, called the
        \textit{\alphath component of $H$.} The multiplication 
        is denoted $\mu_{\alpha}\colon H_{\alpha}\otimes H_{\alpha}\to H_{\alpha}$
        and the unit is denoted $\eta_{\alpha}\colon \Bbbk\to H_{\alpha}$, with 
        $1_{\alpha}\eqdef\eta_{\alpha}(1)$. 
  \item A family of algebra morphisms 
        $\Delta=\bigl\{\Delta_{\alpha,\beta}\colon H_{\alpha\beta}\to H_{\alpha}\otimes H_{\beta}\bigr\}_{\alpha,\beta\in\pi}$,
        called \textit{comultiplication,} that is \textit{coassociative} 
        in the sense that, for any $\alpha,\beta,\gamma\in\pi$, we have
        $(H_{\alpha} \otimes \Delta_{\beta, \gamma}) \circ \Delta_{\alpha,\beta\gamma} = (H_{\alpha\beta}\otimes H_{\gamma}) \circ \Delta_{\alpha\beta,\gamma}$.
  \item An algebra morphism $\varepsilon\colon H_{1}\to  \Bbbk$,
        called \textit{counit,} such that, for any $\alpha\in\pi$, we have
        $(\varepsilon\otimes H_{\alpha})\circ \Delta_{1,\alpha} = \Id$ and $(H_{\alpha}\otimes \varepsilon) \circ \Delta_{\alpha,1}=\Id$.
 \item A set of algebra isomorphisms 
       $\varphi= \bigl\{\varphi_{\beta}^{\alpha}\colon H_{\alpha}\to H_{\beta\alpha\beta^ {-1}} \bigr\}_{\alpha,\beta\in\pi}$,
       called \textit{conjugation.}
       When not strictly necessary, the upper index will be omitted. 
       We require that $\varphi$ satisfies the following conditions. 
       \begin{itemize}
         \item $\varphi$ is \textit{multiplicative,} i.e., for any $\beta,\gamma\in\pi$,
               we have $\varphi_{\beta}\circ\varphi_{\gamma}=\varphi_{\beta\gamma}$.
               It follows that, for any $\alpha\in\pi$, we have 
               $\varphi_{1}^{\alpha}=\Id_{H_{\alpha}}$.
         \item $\varphi$ is \textit{compatible with $\Delta$,} 
               i.e, for any $\alpha,\beta,\gamma\in\pi$, we have
               $(\varphi_{\gamma}\otimes\varphi_{\gamma}) \circ \Delta_{\alpha,\beta}= \Delta_{\gamma\alpha\gamma^{-1},\gamma\beta\gamma^{-1}} \circ \varphi_{\gamma}$.
         \item $\varphi$ is \textit{compatible with $\varepsilon$,} i.e., 
               for any $\gamma\in\pi$, we have $\varepsilon = \varphi_{\gamma}\circ\varepsilon$.
       \end{itemize}
  \item Finally, a set of \klin morphisms  
        $s=\{s_{\alpha} \colon H_{\alpha}\to H_{\alpha^{-1}}\}_{\alpha\in\pi}$,
        called \textit{antipode,} such that, for any $\alpha\in\pi$, we have
        $\mu_{\alpha}\circ (H_{\alpha}\otimes s_{\alpha^{-1}})\circ \Delta_{\alpha,\alpha^{-1}} =
         \eta_{\alpha}\circ \varepsilon = 
         \mu_{\alpha}\circ (s_{\alpha^{-1}}\otimes H_{\alpha})\circ \Delta_{\alpha^{-1},\alpha}$.
\end{itemize}
 We{\label{pag:tf}} say that $H$ is \textit{of finite-type} when any 
 component $H_{\alpha}$ (with $\alpha\in\pi\/$) is a finite-dimensional \kvector space. 
 We say that $H$ is \textit{totally-finite} when $\dim_{\Bbbk}\bigoplus_{\alpha\in\pi}H_{\alpha}<\infty$, 
 i.e., when $H$ is of finite-type and almost all the $H_{\alpha}$ are zero. 
 It was proved in~\cite{Virelizi} that the antipode of a finite-type 
 \Tcoalg is always bijective.
 We observe that the component $H_{1}$ of a \Tcoalg $H$ is a Hopf algebra
 in the usual sense. 

\paragraph{\scshape Coopposite \protect\Tcoalg}
 Let $H$ be a \Tcoalg with invertible antipode. 
 The \textit{coopposite \Tcoalg}$H^{\cop}$ 
 is the \Tcoalg defined as follows.
      \begin{itemize}
       \item For any $\alpha\in\pi$, we set $H^{\cop}_{\alpha}\eqdef H_{\alpha^{-1}}$ 
             as an algebra.
       \item The comultiplication $\Delta^{\cop}$ is obtained by setting, 
             for any $\alpha,\beta\in\pi$,
             \begin{equation*}
              \Delta^{\cop}_{\alpha,\beta}\eqdef\Bigl(H^{\cop}_{\alpha\beta}=H_{\beta^{-1}\alpha^{-1}}
              \xrightarrow{\Delta_{\beta^{-1},\alpha^{-1}}}H_{\beta^{-1}}\otimes H_{\alpha^{-1}} 
              \xrightarrow{\sigma}H_{\alpha^{-1}}\otimes H_{\beta^{-1}}
              =H^{\cop}_{\alpha}\otimes H^{\cop}_{\beta}\Bigr)
              \text{.}
             \end{equation*}
             The counit is given by $\varepsilon^{\cop}=\varepsilon$.
       \item The antipode $s^{\cop}$ is obtained by setting
             $s^{\cop}_{\alpha}=s^{-1}_{\alpha}$, for any $\alpha\in\pi$.
       \item The conjugation  $\varphi^{\cop}$ is obtained by setting
             $\varphi^{\cop}_{\beta}\eqdef\varphi_{\beta}$, for any $\beta\in\pi$.
      \end{itemize}

\paragraph{\scshape Heynemann-Sweedler notation}
 The coassociativity of $H$ allows us to introduce 
 an analog of the Heynemann\hspace{0pt}-\hspace{0pt}Sweedler 
 notation~\cite{Sweedler}. 
 Given $\alpha_{1},\dots, \alpha_{n}\in\pi$ and defined
 \begin{equation*}\begin{split}
  \Delta_{\alpha_{1},\alpha_{2},\dots,\alpha_{n}}\eqdef \Bigl(& H_{\alpha_{1}\alpha_{2}\cdots\alpha_{n}}
  \xrightarrow{\Delta_{\alpha_{1},\alpha_{2}\cdots\alpha_{n}}} H_{\alpha_{1}}\otimes H_{\alpha_{2}\cdots\alpha_{n}}
    \xrightarrow{H_{1}\otimes\Delta_{\alpha_{2},\alpha_{3}\cdots\alpha_{n}}}\\
    & H_{1}\otimes H_{2}\otimes H_{\alpha_{3}\cdots\alpha_{n}}
  \to\ \cdots\ \to H_{\alpha_{1}}\otimes H_{\alpha_{2}}\otimes \cdots\otimes H_{\alpha_{n}}\Bigr)\text{,}
 \end{split}\end{equation*}
 for any $h\in H_{\alpha_{1}\alpha_{2}\cdots\alpha_{n}}$, we set
 $h'_{(\alpha_{1})}\otimes h''_{(\alpha_{2})}\otimes\cdots\otimes {h'}{}^{n}_{(\alpha_{n})}\eqdef
  \Delta_{\alpha_{1},\alpha_{2},\dots,\alpha_{n}}(h)$.
 Let $H$ be a \Tcoalg over a group $\pi$ and a field $\Bbbk$
 and let $M$ be a \kvector space. Suppose that 
 $f\colon H_{\alpha_{1}}\times H_{\alpha_{2}}\times\cdots\times H_{\alpha_{n}}\longrightarrow M$
 is a \kmlin map. Denoted $\hat{f}$ the tensor lift of $f$, 
 we introduce the notation
 $f\Bigl(h'_{(\alpha_{1})}, h''_{(\alpha_{1})},\cdots, {h'}{}^{n}_{(\alpha_{n})}\Bigr)\eqdef
  \hat{f}\bigl(\Delta_{\alpha_{1},\alpha_{2},\dots,\alpha_{n}}(h)\bigr)$.
 For simplicity, we also suppress the subscript ``$(\alpha_{i})$'' 
 when $\alpha_{i}=1$.

 \paragraph{\scshape Quasitriangular \protect\Tcoalgs}
 A \textit{quasitriangular \Tcoalgs} (see~\cite{Tur-CPC})
 is a \Tcoalg $H$ endowed with a family
 $R=\{R_{\alpha,\beta}=\xi_{(\alpha).i}\otimes\zeta_{(\beta).i}\in H_{\alpha}\otimes H_{\beta}\}_{\alpha,\beta\in\pi}$,
 the \textit{universal \Rmatrix,} such that $R_{\alpha,\beta}$ 
 is invertible for any $\alpha,\beta\in\pi$ and the following conditions are satisfied.
 \begin{itemize}\begin{subequations}\label{e:R}
   \item For any $\alpha,\beta\in\pi$ and $h\in H_{\alpha\beta}$,
         \begin{equation}\label{e:R-a}     
          R_{\alpha,\beta}\Delta_{\alpha,\beta}(h)=\bigl(\sigma\circ(\varphi_{\alpha^{-1}}\otimes H_{\alpha})\circ
          \Delta_{\alpha\beta\alpha^{-1},\alpha}\bigr)(h)R_{\alpha,\beta}\text{.}
         \end{equation}
   \item For any $\alpha,\beta,\gamma\in\pi$,
         \begin{equation}\label{e:R-b}
          (H_{\alpha}\otimes\Delta_{\beta,\gamma})(R_{\alpha,\beta\gamma})=(R_{\alpha,\gamma})_{1\beta 3}
          (R_{\alpha,\beta})_{12\gamma}\text{,}
         \end{equation}
         where, given two vector spaces $P$ and $Q$ over $\Bbbk$,
         for any $x=p_{i}\otimes q_{i}\in P\otimes Q$ we set
         $x_{1\beta 3}=p_{i}\otimes 1_{\beta}\otimes q_{i}$
         and
         $x_{12\gamma}=p_{i}\otimes q_{i}\otimes 1_{\gamma}$.
   \item For any $\alpha,\beta,\gamma\in\pi$,
         \begin{equation}\label{e:R-c}
          (\Delta_{\alpha,\beta}\otimes H_{\gamma})(R_{\alpha\beta,\gamma})=
          \bigl((\varphi_{\beta}\otimes H_{\gamma})(R_{\beta^{-1}\alpha\beta,\gamma})\bigr)_{1\beta 3}(R_{\beta,\gamma})_{\alpha 23}\text{,}
         \end{equation}
         where, given two vector spaces $P$ and $Q$,
         for any $x=p_{i}\otimes q_{i}\in P\otimes Q$ we set
         $x_{\alpha 23}= 1_{\alpha}\otimes p_{i}\otimes q_{i}$.
   \item $R$ is \textit{compatible with $\varphi$,} in the sense that,
         for any $\alpha,\beta,\gamma\in\pi$, we have
         \begin{equation}\label{e:R-d}      
          (\varphi_{\alpha}\otimes\varphi_{\alpha})(R_{\beta,\gamma})=R_{\alpha\beta\alpha^{-1},\alpha\gamma\alpha^{-1}}\text{.}
         \end{equation}
  \end{subequations}\end{itemize}
 
\noindent Notice that $(H_{1},R_{1,1})$ is a quasitriangular
Hopf algebra in the usual sense. 
For any $\alpha, \beta\in\pi$, we introduce the notation
$\tilde{\xi}_{(\alpha).i}\otimes\tilde{\zeta}_{(\beta).i}=\tilde{R}_{\alpha,\beta}=R^{-1}_{\alpha,\beta}$.

Following~\cite{Virelizi}, we set
$u_{\alpha}\eqdef (s_{\alpha^{-1}}\circ\varphi_{\alpha})(\zeta_{(\alpha^{-1}).i})\xi_{(\alpha).i}\in H_{\alpha}$
 and $u\eqdef\{u_{\alpha}\}_{\alpha\in\pi}$. 
 The $u_{\alpha}$ are called \textit{Drinfeld elements} of $H$.
 When $\pi=1$, we recover the usual definition 
 of Drinfeld element of a quasitriangular Hopf algebra.
 The following properties of $u$ are proved in~\cite{Virelizi}. 
 Let $\alpha$ and $\beta$ be in $\pi$ and let $h$ be in $H_{\alpha}$.
 \begin{eqnenumerate}\label{Ovid}
  \item $u_{1}=s_{1}(\zeta_{(1).i})\xi_{(1).i}$.\label{Ovid-1}
  \item $u_{\alpha}$ is invertible with inverse 
        $u^{-1}_{\alpha}=s^{-1}_{\alpha}\Bigl(\tilde{\zeta}_{(\alpha^{-1}).i}\Bigr)\tilde{\xi}_{(\alpha).i}$.
        Moreover we have{\label{Ovid-2}}\
        $u_{\alpha}^{-1}=\xi_{(\alpha).i}(s_{\alpha^{-1}}\circ s_{\alpha})(\zeta_{(\alpha).i})$.
  \item $(u_{\alpha\beta})'_{(\alpha)}\otimes (u_{\alpha\beta})''_{(\beta)}=
         \tilde{\xi}_{(\alpha).i}\tilde{\zeta}_{(\alpha).j}u_{\alpha}\otimes\tilde{\zeta}_{(\beta).i}
         \varphi_{\alpha^{-1}}(\tilde{\xi}_{(\alpha\beta\alpha^{-1}).j})u_{\beta}$.\label{Ovid-3}
  \item $\varepsilon(u_{1})=1$.\label{Ovid-4}
  \item $s_{\alpha^{-1}}(u_{\alpha^{-1}})u_{\alpha}=u_{\alpha}s_{\alpha^{-1}}(u_{\alpha^{-1}})$.\label{Ovid-5}
  \item $\varphi_{\beta}(u_{\alpha})=u_{\beta\alpha\beta^{-1}}$.\label{Ovid-6}
  \item $(s_{\alpha^{-1}}\circ s_{\alpha}\circ\varphi_{\alpha})(h)=u_{\alpha}hu_{\alpha}^{-1}$.\label{Ovid-7}
  \item $u_{\alpha}s_{\alpha^{-1}}(u_{\alpha^{-1}})h=
         \varphi_{\alpha^{2}}(h)u_{\alpha}s_{\alpha^{-1}}(u_{\alpha^{-1}})$.\label{Ovid-8}
 \end{eqnenumerate}

\paragraph{\scshape The mirror \protect\Tcoalg}\label{par:mirror}
 Let $H=(H,R)$ be a quasitriangular Hopf algebra (with $R=\xi_{i}\otimes\zeta_{i}$ and 
 $R^{-1}=\tilde{R}=\tilde{\xi}_{i}\otimes\tilde{\zeta}_{i}\/$). By replacing $R$ with 
 $\overline{R}=\sigma(\tilde{R})=\tilde{\zeta}_{i}\otimes\tilde{\xi}_{i}$ we obtain another 
 quasitriangular structure $\overline{H}=(H,\overline{R})$. 
 This means that, in the category of 
 representations of $H$, we replace the braiding $c_{R}$ provided 
 by $R$ by the braiding $c_{R}^{-1}$ provided by $\overline{R}$. 
 When $H$ is a \Tcoalg, the family 
 $\{R^{-1}_{\alpha,\beta}=\tilde{\zeta}_{\alpha.i}\otimes\tilde{\xi}_{\beta.i}\}_{\alpha,\beta\in\pi}$ 
 is not a universal \Rmatrix for $H$. Nevertheless,
 it is still possible to generalize
 the definition of $\overline{H}$.
Let $H$ be a \Tcoalg. The \Tcoalg $\overline{H}$, called 
\textit{mirror of $H$}~\cite{Tur-CPC}, is defined as follows.
\begin{itemize}
 \item For any $\alpha\in\pi$, we set $\overline{H}_{\alpha}\eqdef H_{\alpha^{-1}}$.
  \item For any $\alpha,\beta\in\pi$, the component $\overline{\Delta}_{\alpha,\beta}$ of the 
        comultiplication $\overline{\Delta}$
        of $\overline{H}$ is
        \begin{equation}\label{e:mirror-temp}
         \overline{\Delta}_{\alpha,\beta}(h)\eqdef \bigl((\varphi_{\beta}\otimes H_{\beta^{-1}})\circ
         \Delta_{\beta^{-1}\alpha\beta,\beta^{-1}}\bigr)(h)\in H_{\alpha^{-1}}\otimes 
         H_{\beta^{-1}}=\overline{H}_{\alpha}\otimes\overline{H}_{\beta}\text{,}
        \end{equation}
        for any $h\in H_{\beta^{-1}\alpha^{-1}}=\overline{H}_{\alpha\beta}$. If,  we set
        $h'_{\overline{(\alpha)}}\otimes h''_{\overline{(\beta)}}=\overline{\Delta}_{\alpha,\beta}(h)$,
        then~\eqref{e:mirror-temp} can be written in the form
        $h'_{\overline{(\alpha)}}\otimes h''_{\overline{(\beta)}}\eqdef 
        \varphi_{\beta}(h'_{(\beta^{-1}\alpha^{-1}\beta)})\otimes 
         h''_{(\beta^{-1})}$.
        The counit of $\overline{H}$ is given by
        $\varepsilon\in H_{1}^{\ast}=\overline{H}_{1}^{\ast}$.
  \item For any $\alpha\in\pi$, the \alphath component 
        of the antipode $\overline{s}$ of $\overline{H}$ 
        is given by
        $\overline{s}_{\alpha}=\varphi_{\alpha}\circ s_{\alpha^{-1}}\colon\overline{H}_{\alpha}=H_{\alpha^{-1}}\to 
         H_{\alpha}=\overline{H}_{\alpha^{-1}}$.
  \item Finally, for any $\alpha\in\pi$, we set $\overline{\varphi}_{\alpha}=\varphi_{\alpha}$.
 \end{itemize}
 
 If $H$ is quasitriangular, then $\overline{H}$ 
 is also quasitriangular with universal \Rmatrix
 \begin{equation}\label{e:R-mirror}
  \overline{\xi}_{(\alpha).i}\otimes\overline{\zeta}_{(\beta).i}=
  \overline{R}_{\alpha,\beta}=\bigl(\sigma(R_{\beta^{-1},\alpha^{-1}})\bigr)^{-1}\in 
  H_{\alpha^{-1}}\otimes H_{\beta^{-1}}=\overline{H}_{\alpha}\otimes\overline{H}_{\beta}
 \end{equation}
 for any $\alpha,\beta\in\pi$. 
 
 \paragraph{\scshape Ribbon \protect\Tcoalgs}
We say that a quasitriangular \Tcoalg $H$ is \textit{ribbon} 
when it is endowed with a family
 $\theta=\{\theta_{\alpha}\vert\theta_{\alpha}\in H_{\alpha}\}_{\alpha\in\pi}$
 such that $\theta_{\alpha}$ is invertible for any $\alpha\in\pi$ and 
 the following conditions are satisfied
 for any $\alpha,\beta\in\pi$ and $h\in H_{\alpha}$.
 \begin{enumerate}
  \item $\varphi_{\alpha}(h) = \theta^{-1}_{\alpha}h\theta_{\alpha}$.
  \item $s_{\alpha}(\theta_{\alpha})=\theta_{\alpha^{-1}}$.
  \item $(\theta_{\alpha\beta})'_{(\alpha)}\otimes (\theta_{\alpha\beta})''_{(\beta)}=
        \theta_{\alpha}\zeta_{(\alpha).i}\xi_{(\alpha).j}\otimes \theta_{\beta}        
        \varphi_{\alpha^{-1}}(\xi_{(\alpha\beta\alpha^{-1}).i})\zeta_{(\beta).j}$.\label{p:bkdin3}
  \item $\varphi_{\beta}(\theta_{\alpha})=\theta_{\beta\alpha\beta^{-1}}$.
\end{enumerate}
\noindent Notice that $(H_{1},R_{1,1},\theta_{1})$ is a ribbon Hopf algebra 
in the usual sense.

\paragraph{\scshape\protect\Tcats of representations}
 Let $H$ be a \Tcoalg over a field $\Bbbk$. The \Tcat $\Rep(H)$ (see~\cite{Tur-CPC}) is defined as follows.
\begin{itemize}
 \item For any $\alpha\in\pi$, the \alphath component of $\Rep(H)$, denoted 
       $\Rep_{\alpha}(H)$, is the category of representations of the algebra 
       $H_{\alpha}$.
 \item The tensor product $U\otimes V$ of $U\in\Rep_{\alpha}(H)$ and $V\in\Rep_{\beta}(H)$ 
       (with $\alpha,\beta\in\pi\/$) is obtained by the tensor product of 
       \kvector spaces
       $U\otimes_{\Bbbk}V$ endowed with the action of $H_{\alpha\beta}$ given by
       $h(u\otimes v) = \Delta_{\alpha,\beta}(h)(u\otimes v)= h'_{(\alpha)}u\otimes h''_{(\beta)}v$,
       for any $h\in H_{\alpha\beta}$, $u\in U$, and $v\in V$.
 \item The tensor product of two arrows $f\in\Rep_{\alpha}(H)$ and $g\in\Rep_{\beta}(H)$
       is given by the tensor product of \klin morphisms, i.e,
       the forgetful functor from $\Rep(H)$ to the
       category of vector spaces over $\Bbbk$ is faithful.
 \item The unit $\tunit$ is the ground field $\Bbbk$ with the 
       action of $H_{1}$ provided by $\varepsilon$.
 \item Given $\beta\in\pi$, we need to define the functor $\lidx{\beta}{\,(\cdot)}$.
       To avoid confusion, in this context we reserve the notation 
       $\varphi_{\beta}$ for the isomorphism of algebras
       $\varphi_{\beta}\colon H_{\alpha}\to H_{\beta\alpha\beta^{-1}}$.
       Let $U$ be in $\Rep_{\alpha}(H)$, with $\alpha\in\pi$. The object $\lidx{\beta}{U}$
       has the same underling vector space of $U$. 
       Given $u\in U$, we denote $\lidx{\beta}{\, u}$ the corresponding
       element in $\lidx{\beta}{U}$. The action of $H_{\beta\alpha\beta^{-1}}$ on
       $\lidx{\beta}{U}$ is given by
       \begin{equation}\label{e:shift-prod}
           h\,\lidx{\beta}{\,u}=\lidx{\beta}{\,\bigl(\varphi_{\beta^{-1}}(h)u\bigr)}
       \end{equation}
       for any $u\in U$ and $h\in H_{\beta\alpha\beta^{-1}}$.
\end{itemize}
The objects of $\Rep(H)$ are called \textit{representations of $H$.}
   
When $H$ is quasitriangular, $\Rep(H)$ is braided by setting,
for any $u\in U$, $v\in V$, $U\in\Rep_{\alpha}(H)$, $V\in\Rep_{\beta}(H)$, and $\alpha,\beta\in\pi$,
\begin{equation*}
    \map{c_{U,V}}{U\otimes V}{\Bigl(\lidx{U}{V}\Bigr)\otimes U}{u\otimes v}%
    {\Bigl(\lidx{\alpha}{\,(\zeta_{(\beta).i}v)}\Bigr)\otimes\xi_{(\alpha).i}u}\text{.}
\end{equation*}
   
Let us consider the full subcategory $\Repf(H)$ of
\textit{finite-dimensional representations of $H$,} i.e., 
of representations $U$ of $H$ such that $\dim_{\Bbbk}U\in\N$.
The \Tcat $\Repf(H)$ has a structure left autonomous
\Tcat obtained in the following way. For any $U\in\mathcal{T}_{\alpha}$
we set $U^{\ast}=\Hom_{\Bbbk}(U,\Bbbk)$, with the action of $H_{\alpha^{-1}}$ on 
$U^{\ast}$ given by
\begin{subequations}\label{e:tardi}
\begin{equation}
    \langle hf, u\rangle = \bigl\langle f, s_{\alpha^{-1}}(h)u\bigr\rangle
\end{equation}
for any $h\in H_{\alpha^{-1}}$, $f\in U^{\ast}$ and $u\in U$. Then, $U^{\ast}$ is a 
left dual of $U$ via
\begin{equation}
    \map{b_{U}}{\tunit}{U\otimes U^{\ast}}{1}{e_{i}\otimes e^{i}}
\end{equation}
(where $\{e_{i}\}$ is a basis of $U$ as a \kvector space and $\{e^{i}\}$ 
its dual basis), and
\begin{equation}
  \map{d_{U}}{U^{\ast}\otimes U}{\Bbbk}{f\otimes u}{\langle f,u\rangle=f(u)}
\end{equation} 
(for any $f\in U^{\ast}$ and $u\in U\/$).
\end{subequations}
   
If $H$ is endowed with a twist $\{\theta_{\alpha}\in H_{\alpha}\}_{\alpha\in\pi}$, then $\Rep(H)$
is a balanced \Tcat, with the twist
$\map{\theta_{U}}{U}{\lidx{U}{U}}{u}{\lidx{\alpha}{\,\bigl(\theta_{\alpha}u\bigr)}}$,
for any $u\in U$, with $U\in\Rep_{\alpha}(H)$, and $\alpha\in\pi$. Similarly,
$\Repf(H)$ is a ribbon \Tcat.
   
Notice that $\overline{\Rep(H)}$ is isomorphic
to $\Rep\Bigl(\overline{H}\Bigr)$. 
Similarly, $\overline{\Repf(H)}=\Repf\Bigl(\overline{H}\Bigr)$.
  
\section{The center of a \protect\Tcat}\label{sec:Z}
We generalize the center construction of a tensor category 
described in~\cite{JS} to the case of a \Tcat $\mathcal{T}$
(that, for simplicity, we suppose strict), obtaining a
braided \Tcat $\mathcal{Z}(\mathcal{C})$ in the following way.
\begin{itemize}\begin{subequations}
 \item The objects of $\mathcal{Z}(\mathcal{T})$, called 
       \textit{half-braidings,} 
       are the pairs $(U, \halfbrd_{\_})$ satisfying the following
       conditions.
       \begin{itemize}
        \item $U$ is an object of $\mathcal{T}$.
        \item $\halfbrd_{\_}$ is a natural isomorphism from the functor 
              $U\otimes\_$ to the functor $\lidx{U}{(\_)}\otimes U$ such that for any
              $X, Y\in\mathcal{T}$, we have
              \begin{equation}\label{e:ax-c}
               \halfbrd_{X\otimes Y} = 
               \biggl(\Bigl(\lidx{U}{X}\Bigr)\otimes \halfbrd_{Y}\biggr)\circ
               (\halfbrd_{X}\otimes Y)\text{.}
              \end{equation}
       \end{itemize}
 \item The arrows in $\mathcal{Z}(\mathcal{T})$ from an object 
       $(U,\halfbrd_{\_})$ to an object  $(V,\halfbrdtwo_{\_})$ are 
       the arrows $f\in\mathcal{T}(U,V)$ such that, for any 
       $X\in\mathcal{T}$, we have
       \begin{equation}\label{e:ax-arr}
         \biggl(\Bigl(\lidx{U}{X}\Bigr)\otimes f\biggr)\circ \halfbrd_{X} 
         = \halfbrdtwo_{X}\circ (f\otimes X)\text{.}
       \end{equation}
       The composition of two arrows in $\mathcal{Z}(\mathcal{T})$ 
       is given by the composition in $\mathcal{T}$,
       i.e., by requiring that the forgetful 
       $\mathcal{Z}(\mathcal{T})\to\mathcal{T}\colon (U,\halfbrd_{\_})\mapsto U$
       is faithful.
 \item Given{\label{p:Dr-Moriarty}} 
       $Z=(U,\halfbrd_{\_}), Z'=(U',\halfbrd'_{\_})\in\mathcal{Z}(\mathcal{T})$, 
       their tensor product $Z\otimes Z'$ in $\mathcal{Z}(\mathcal{C})$ is the
       couple $\bigl(U\otimes U', (\halfbrd\tenbrd\halfbrd')_{\_}\bigr)$, where 
       $(\halfbrd\tenbrd\halfbrd')_{\_}$ is obtained by setting, 
       for any $X\in\mathcal{T}$, 
       \begin{equation}\label{e:ax-prod} 
         (\halfbrd\tenbrd\halfbrd')_{X}\eqdef 
         (\halfbrd_{\lidx{U'}{X}}\otimes U')\circ (U\otimes \halfbrd'_{X})\text{.}
       \end{equation}
 \item The tensor unit of $\mathcal{Z}(\mathcal{T})$ is the couple 
       $Z_{\tunit}=(\tunit,\Id_{\_})$, where $\tunit$ is the tensor unit of
       $\mathcal{T}$.
 \item For any $\alpha\in\pi$, the \alphath component of 
       $\mathcal{Z}(\mathcal{T})$, denoted $\mathcal{Z}_{\alpha}(\mathcal{T})$,
       is the full subcategory of $\mathcal{Z}(\mathcal{T})$ whose 
       objects are the pairs $(U,\halfbrd_{\_})$ with $U\in\mathcal{T}_{\alpha}$. 
 \item For any $\beta\in\pi$, the automorphism $\varphi_{\mathcal{Z}.\beta}$ is obtained by 
       setting, for any $(U,\halfbrd_{\_})\in\mathcal{Z}(\mathcal{T})$,
       \begin{equation}\label{e:ax-aut}
         \varphi_{\mathcal{Z}.\beta}(U,\halfbrd_{\_})=
         \bigl(\varphi_{\beta}(U), \varphi_{\mathcal{Z}.\beta}(\halfbrd_{\_})\bigr)\text{,}
       \end{equation}
       where $\varphi_{\mathcal{Z}.\beta}(\halfbrd)_{X}=\varphi_{\beta}
        \bigl(\halfbrd_{\varphi_{\beta}^{-1}(X)}\bigr)$, for any $X\in\mathcal{T}$.
       The definition of $\varphi_{\beta}$ is completed by requiring 
       that the forgetful $\mathcal{Z}(\mathcal{T})\to\mathcal{T}$
       is a \Tfun.
 \item The braiding $c$ in $\mathcal{Z}(\mathcal{T})$ is obtained by 
       setting, $c_{Z,Z'}= \halfbrd_{U'}$, for any
       $Z=(U,\halfbrd_{\_}), Z'=(U',\halfbrd'_{\_})
       \in\mathcal{Z}(\mathcal{T})$.
 \end{subequations}\end{itemize}

\begin{thm}\label{thm:Z}
 $\mathcal{Z}(\mathcal{T})$ is a braided \Tcat.
\end{thm}

 The proof of Theorem~\ref{T:7} and of most of the
 results in the next three sections can be obtained by
 introducing a graphical calculus as
 is~\cite{Tur-QG} and in~\cite{KasTur} and then to follow,
 \textit{mutatis mutandis,} the computation in~\cite{KasTur}, with
 the main difference that we does not obtain the categorical
 quantum double directly, but in more steps, as
 described before. Alternatively, it is possible to made every
 computation algebraically. This has the advantage that you can
 consider \Tcats that are not equivalent to strict \Tcats.
 In particular, you can consider \Tcats such that the conjugate
 automorphism $\varphi$ are not strict (see note at the end of the proof
 of Theorem~\ref{lemma:strict}). However, in that way, computations
 become very long.

\begin{rmk}
Brugui\`eres~\cite{Bruguieres} noticed that, if $\mathcal{C}$ is an
autonomous tensor category, then $\mathcal{Z}(\mathcal{T})$ is autonomous
too. He also noticed that this result is still true if we replace 
$\mathcal{C}$ with a \Tcat $\mathcal{T}$.
\end{rmk}
 
\begin{rmk}\label{rmk:si01}
 The definition of the center given here generalizes the most usual
 convention adopted also in~\cite{JS,KasTur}. 
 However, in~\cite{Kas}, the center of a tensor category is constructed
 in a similar way, but considering the natural transformation of the kind
 $\_\otimes U\to U\otimes\_$.  The choice in~\cite{Kas} seems more appropriate in
 some context, e.g., in the construction of the isomorphism between the
 center of the category of representations of a Hopf algebra $H$ and the
 category of representations of $D(H)$. 
\end{rmk}

\section{The twist extension of a braided \protect\Tcat}\label{sec:Rib}
Let $\mathcal{T}$ be a braided \Tcat
(again, for simplicity, we suppose $\mathcal{T}$ strict). 
Generalizing the construction described in~\cite{Street-double}, 
we obtain a balanced \Tcat $\mathcal{T}^{Z}$, the 
\textit{twist extension of $\mathcal{T}$.} Even if 
we do not have, in general, an embedding 
$\mathcal{T}\hookrightarrow\mathcal{T}^{Z}$, the name is justified by 
the observation that we still have an embedding 
$\mathcal{T}_{1}\hookrightarrow\mathcal{T}^{Z}_{1}$. We will see
that, if $H$ is a \Tcoalg and $\mathcal{T}=\Rep(H)$ or 
$\mathcal{T}=\Repf(H)$, then there is an embedding 
$\mathcal{T}\hookrightarrow\mathcal{T}^{Z}$.

\begin{itemize}\begin{subequations}
 \item The objects of $\mathcal{T}^{Z}$ are the pairs $T=(U,t)$, 
       where $U\in\mathcal{T}$ and
       $t\in\mathcal{T}\Bigl(U,\lidx{U}{U}\Bigr)$ is invertible.
 \item For any $T_{1}=(U_{1},t_{1}), T_{2}=(U_{2},t_{2})\in\mathcal{T}^{Z}$,
       the arrows from $T_{1}$ to $T_{2}$ in $\mathcal{T}^{Z}$ 
       are the arrows $f\in\mathcal{T}(U_{1},U_{2})$ such that 
       \begin{equation}\label{e:twistator-arrow}
        \Bigl(\lidx{U}{f}\Big)\circ t_{1}=t_{2}\circ f\text{.}
       \end{equation}
       The composition is given by the composition in $\mathcal{T}$, 
       i.e., we require that the forgetful functor
       from $\mathcal{T}^{Z}\to\mathcal{T}\colon(U,t)\mapsto U$ is faithful.
 \item The tensor product of $T_{1}=(U_{1},t_{1}), T_{2}=(U_{2},t_{2})
       \in\mathcal{T}^{Z}$ is the couple
       $T_{1}\twten T_{2}=(U_{1}\otimes U_{2},t_{1}\twten t_{2})$, where
       \begin{equation}\label{e:twistator}
         t_{1}\twten t_{2}\eqdef c_{\lidx{U\otimes U'}{U'},\lidx{U}{U}}\circ
         c_{\lidx{U}{U},\lidx{U'}{U'}}\circ(t_{1}\otimes t_{2})
       \end{equation}
 \item The tensor product of two arrows in $\mathcal{T}^{Z}$ is 
       given by the tensor product of arrows in $\mathcal{T}$.
 \item The tensor unit in $\mathcal{T}^{Z}$ is the couple 
       $T_{\tunit}=(\tunit,\Id_{\tunit})$, where $\tunit$ is the tensor
       unit of $\mathcal{T}$.
 \item For any $\alpha\in\pi$, the component 
       $\mathcal{T}_{\alpha}^{Z}=\bigl(\mathcal{T}^{Z}\bigr)_{\alpha}$ is
       the full subcategory of $\mathcal{T}^{Z}$ whose objects are 
       the pairs $(U,t)$ with $U\in\mathcal{T}_{\alpha}$.
 \item For any $\beta\in\pi$, the functor $\varphi^{Z}_{\beta}$ is obtained by setting,
       for any $(U,t)\in\mathcal{T}^{Z}$,
       \begin{equation}\label{e:twtw}
         \varphi^{Z}_{\beta}(U,t)=\bigl(\varphi_{\beta}(U),\varphi_{\beta}(t)\bigr)
       \end{equation}
       and $\varphi^{Z}(f)=\varphi(f)$, for any arrow $f$ in $\mathcal{T}^{Z}$.
 \item The braiding in $\mathcal{T}^{Z}$ is obtained by requiring that 
       the forgetful functor from $\mathcal{T}^{Z}$
       to $\mathcal{T}$ is braided.
 \item The twist $\theta$ of $\mathcal{T}^{Z}$ is obtained by 
       setting $\theta_{T}=t$, for any $T=(U,t)\in\mathcal{T}^{Z}$.
 \end{subequations}
\end{itemize}
 
\begin{thm}\label{thm:ZZ}
   $\mathcal{T}^{Z}$ is a balanced \Tcat.
\end{thm}

\paragraph{\scshape Dualities in $\mathcal{T}^{Z}$}
Even when $\mathcal{T}$ is left autonomous, 
an object in $\mathcal{T}^{Z}$ not necessarily admits a left dual. 
So, in particular, $\mathcal{T}^{Z}$ is not necessarily ribbon. 
The following lemma gives a characterization of the objects in 
$\mathcal{T}^{Z}$ endowed with a stable left dual.
  
\begin{lemma}\label{l:iounsa}
 Let $T=(U,t)$ and $T^{\ast}=(U^{\ast},\tau)$ be objects in $\mathcal{T}^{Z}$. 
 Then, $T^{\ast}$ is a stable left dual
 of $T$ with unit $b_{T}$ and counit $d_{T}$ if and only if
 \begin{itemize}
  \item $U^{\ast}$ is a stable left dual of $U$ in $\mathcal{T}$ via unit 
        $b_{U}=b_{T}$ and counit $d_{U}=d_{T}$ and
  \item $\tau=\leftidx{^{U^{\ast}}}{\!\!\hat{t}\,}{^{\ast}}$, where 
        $\hat{t}\in\mathcal{T}\Bigl(U,\lidx{U}{U}\Bigr)$ 
        satisfies $t^{-1}\circ\leftidx{^{U}}{\!\hat{t}\,}{^{-1}}=\omega_{U}$.
 \end{itemize} 
\end{lemma}

\section{Dualities in a balanced \protect\Tcat}\label{sec:functor-N}Let $\mathcal{T}$ be a balanced \Tcat.
Generalizing some results in~\cite{JS-tortile,JS,Tur-QG,KasTur} to the 
case of a \Tcat, we study the properties of dualities in 
$\mathcal{T}$. In particular, this will allow us to obtain a full
subcategory $\mathcal{N}(\mathcal{T})$ of $\mathcal{T}$ that will
be the biggest ribbon category included in $\mathcal{T}$. This is
the analog, in the case of a \Tcat, of the construction given 
in~\cite{Street-double} in the case of a tensor category.
 
\paragraph{\scshape Reflexive objects} 
Let $\mathcal{T}$ be a balanced \Tcat and $U\in\mathcal{T}$. We set
\begin{equation*}
 \theta^{2}_{U}\eqdef\Bigl(U\xrightarrow{\,\theta_{U}\,}\lidx{U}{U}
 \xrightarrow{\,\leftidx{^{U}}{\!\theta}{_{U}}\,} \lidx{U\otimes U}{U}\Bigr)
\end{equation*}
and $\theta^{-2}_{U}=(\theta^{2}_{U})^{-1}$.
We say that $U$ is \textit{reflexive} if it is endowed with a stable 
left dual $U^{\ast}$, such that $\theta^{-2}_{U} =\omega_{U}$.

\begin{lemma}\label{l:reflexivity}
 If $U\in\mathcal{T}$ has a stable left dual $U^{\ast}$ such that the ribbon 
 condition~\eqref{e:tortility} is satisfied,
 then $U$ is reflexive. In particular, any object in a ribbon
 \Tcat is reflexive.
\end{lemma}

\paragraph{\scshape Reversed duality}\label{p:r-d}
Let $U$ be a reflexive object in $\mathcal{T}$. We set
\begin{equation*}
   \begin{cases}  b'_{U} \!\!\!\! &= 
   \biggl(\tunit\xrightarrow{\,\leftidx{^{U^{\ast}}}{\! b}{_{U}}\,}
   \lidx{U^{\ast}}{U}\otimes\leftidx{^{U^{\ast}}}{\! U}{^{\ast}}
   \xrightarrow{\,\leftidx{^{U^{\ast}}}{\! c}{_{U,U^{\ast}}}\,} U^{\ast}\otimes\lidx{U^{\ast}}{U}
   \xrightarrow{\,U^{\ast}\otimes\leftidx{^{U^{\ast}}}{\!\theta}{_{U}}\,} U^{\ast}\otimes U\biggr)\\
   d'_{U} \!\!\!\! &= \biggl(U\otimes U^{\ast}\xrightarrow{\,\theta_{U}\otimes U^{\ast}\,}
   \lidx{U}{U}\otimes U^{\ast}\xrightarrow{\,c_{\lidx{U}{U},U^{\ast}}\,}
   \leftidx{^{U}}{\! U}{^{\ast}}\otimes\lidx{U}{U}
   \xrightarrow{\,\leftidx{^{U}}{\! d}{_{U}}\,}\tunit\biggr)\text{.}
\end{cases}\end{equation*}

\begin{lemma}\label{l:bd1}
 $U$ is a left dual of $U^{\ast}$ via $(b_{U},d'_{U})$, 
 i.e., $(d'_{U}\otimes U)\circ (U\otimes b'_{U}) = U$ and
 $(U^{\ast}\otimes d'_{U})\circ (b'_{U}\otimes U^{\ast}) = U^{\ast}$,
 that is, $U$ is a stable left dual of $U^{\ast}$ with unit $b'_{U}$ 
 and counit $d'_{U}$.
\end{lemma}
 
The adjunction $(b'_{U},d'_{U})$ is
said \textit{reversed adjunction of $(b_{U},d_{U})$.}
 
\paragraph{\scshape Good left duals}
Let $U$ be a reflexive object. We say that $U^{\ast}$ 
is a \textit{good left dual} if further we have
\begin{equation}\label{e:a-tortility}
  \theta_{U^{\ast}}=\leftidx{^{U^{\ast}}}{\!\theta}{^{\ast}_{U}}\text{.}
\end{equation}

\begin{lemma}\label{l:a-tortility}
 Let $U$ be an object in a balanced \Tcat $\mathcal{T}$ 
 endowed with a stable adjunction 
 $(b_{U},d_{U})\colon\adjunction{U^{\ast}}{U}$. 
 The ribbon condition~\eqref{e:tortility} 
 is satisfied if and only 
 if~\eqref{e:a-tortility} is satisfied. In particular, 
 $T$ is ribbon if and only if any object $U\in\mathcal{T}$ 
 satisfies~\eqref{e:a-tortility}.
\end{lemma}

\begin{lemma}\label{l:starstar}
 Let $U^{\ast}$ be a good left dual of $U\in\mathcal{T}$. 
 If we set $U^{\ast\ast}=U$ via the reversed adjunction $(b'_{U},d'_{U})$,
 then we have $b''_{U}= b_{U}$ and $d''_{U}=d_{U}$.
 If $V^{\ast}$  is a good left dual of $V\in\mathcal{T}$ and 
 we set $V^{\ast\ast}=V$ via the reversed adjunction $(b'_{V},d'_{V})$,
 then, for any $f\in\mathcal{T}(U,V)$, we have $f^{\ast\ast}=f$.
\end{lemma}

\paragraph{\scshape The category $\mathcal{N}(\mathcal{T})$} 
 Let $\mathcal{T}$ be a balanced \Tcat. By definition, 
 $\mathcal{N}(\mathcal{T})$ is the full subcategory of 
 $\mathcal{T}$ of the object $U\in\mathcal{T}$ that admits
 a good left dual. For any class $\Phi(U)$ in $\mathcal{N}(\mathcal{T})$
 we also fix an object $U_{0}\in\Phi(U)$ and a good left dual 
 $U^{\ast}_{0}$ of $U_{0}$, obtaining, in that way, a good left 
 dual $V^{\ast}$ for any $V\in\Phi(U)$.
  
\begin{thm}\label{thm:N}
 $\mathcal{N}(\mathcal{T})$ has a structure of balanced \Tcat.
 Moreover, $\mathcal{N}(\mathcal{T})$ is a ribbon \Tcat
 and any other ribbon subcategory of $\mathcal{T}$ is 
 included in $\mathcal{N}(\mathcal{T})$.
\end{thm}
  
\begin{proof}
 The proof that $\mathcal{N}(\mathcal{T})$ has 
 a structure of balanced \Tcat is given in Lemma~\ref{l:N-1}.
 The proof that $\mathcal{N}(\mathcal{T})$ is autonomous 
 is given in Lemma~\ref{l:N-2}.
 Since, by hypothesis, any object of $\mathcal{N}(\mathcal{T})$ 
 satisfies~\eqref{e:a-tortility}, then, by Lemma~\ref{l:a-tortility}, 
 $\mathcal{N}(\mathcal{T})$ is ribbon.
 The fact that any other ribbon \Tcat included in $\mathcal{T}$
 is also included in $\mathcal{N}(\mathcal{T})$ follows 
 by Lemma~\ref{l:reflexivity} and Lemma~\ref{l:a-tortility}.
\end{proof}
  
To prove that $\mathcal{N}(\mathcal{T})$ is a tensor category, 
we need the following observation.
 
\begin{lemma}\label{l:Protagoras}
 Let $U$ and $V$ be objects in $\mathcal{T}$, let $U^{\ast}$ be 
 a stable left dual of $U$, and let $V^{\ast}$ be a stable left 
 dual of $V$. Consider the dual $(U\otimes V)^{\ast}=V^{\ast}\otimes U^{\ast}$ of $U\otimes V$ 
 via the unit $b_{U\otimes V}=(U\otimes b_V\otimes U^{\ast})\circ b_{U}$ and the counit 
 $d_{U\otimes V}=d_{U}\circ (U^{\ast}\otimes d_V\otimes U)$. We have
 $c_{V^{\ast},U^{\ast}}=c^{\ast}_{V,\lidx{V^{\ast}}{U}}$.
\end{lemma}

\begin{lemma}\label{l:N-1}
 The structure of balanced \Tcat of $\mathcal{T}$
 induces a structure of balanced \Tcat on
 $\mathcal{N}(\mathcal{T})$.
\end{lemma}
 
\begin{proof}
 The only non trivial part is to show that $\mathcal{N}(\mathcal{T})$ 
 is a tensor category. Since $\mathcal{N}(\mathcal{T})$
 is a full subcategory of $\mathcal{T}$, we only need to show that
 the tensor product of  $U,V\in\mathcal{N}(\mathcal{T})$ 
 lies in $\mathcal{N}(\mathcal{T})$, i.e., that $U\otimes V$ admits
 a good left dual. Let $U^{\ast}$ be a good left dual of $U$ and $V^{\ast}$ a
 good left dual $V$. Given $V^{\ast}\otimes U^{\ast}$ as a stable left dual of $U\otimes V$
 with unit $b_{U\otimes V}$ and counit $d_{U\otimes V}$ as in
 Lemma~\ref{l:Protagoras}, then
 $V^{\ast}\otimes U^{\ast}$ is a good left dual of $U\otimes V$ since, by 
 Lemma~\ref{l:Protagoras},
 \begin{equation*}\begin{split}
  \theta_{(U\otimes V)^{\ast}} & = \theta_{V^{\ast}\otimes U^{\ast}}
  = c_{\leftidx{^{V^{\ast}\otimes U^{\ast}}}{\! U}{^{\ast}},\leftidx{^{V^{\ast}}}{\! V}{^{\ast}}}
  \circ\biggl(\Bigl(\leftidx{^{V^{\ast}}}{\!\theta}{_{U^{\ast}}}\Bigr)\otimes\theta_{V^{\ast}}\biggr)
  \circ c_{V^{\ast},U^{\ast}}\\
  & = c_{\leftidx{^{V^{\ast}\otimes U^{\ast}}}{\! U}{^{\ast}},\leftidx{^{V^{\ast}}}{\! V}{^{\ast}}}
  \circ\biggl(\Bigl(\leftidx{^{V^{\ast}\otimes U^{\ast}}}{\!\theta}{^{\ast}_{U}}\Bigr)\otimes%
  \leftidx{^{V^{\ast}}}{\!\theta}{^{\ast}_{V}}\biggr) \circ c_{V^{\ast},U^{\ast}}\\
  & = c^{\ast}_{\lidx{V^{\ast}\otimes U^{\ast}}{U},\lidx{V^{\ast}\otimes U^{\ast}}{V}}
  \circ\biggl(\Bigl(\leftidx{^{V^{\ast}\otimes U^{\ast}}}{\!\theta}{^{\ast}_{U}}\Bigr)\otimes%
  \leftidx{^{V^{\ast}}}{\!\theta}{^{\ast}_{V}}\biggr)\circ c^{\ast}_{V,\lidx{V^{\ast}}{U}}\\
  & = \Biggl(c_{V,\lidx{V^{\ast}}{U}}
  \circ\biggl(\Bigl(\leftidx{^{V^{\ast}}}{\!\theta}{_{V}}\Bigr)\otimes
  \leftidx{^{V^{\ast}\otimes U^{\ast}}}{\!\theta}{_{U}}\biggr)
  \circ c_{\lidx{V^{\ast}\otimes U^{\ast}}{U},\lidx{V^{\ast}\otimes U^{\ast}}{V}}\Biggr)\\
  & = \leftidx{^{V^{\ast}\otimes U^{\ast}}}{\!\Biggl(c_{\lidx{U\otimes V}{V},\lidx{U}{U}}
  \biggl(\Bigl(\lidx{^{U}}{\!\theta}{_{V}}\Bigr)\otimes\theta_{U}\biggr)
  \circ c_{U,V}\Biggr)}{^{\ast}}= \leftidx{^{V^{\ast}\otimes U^{\ast}}}{\!\theta}{^{\ast}_{U\otimes V}}\text{.}
 \end{split}\end{equation*}        
\end{proof}
  
\begin{lemma}\label{l:N-2}
 $\mathcal{N}(\mathcal{T})$ is an autonomous \Tcat.
\end{lemma}
 
\begin{proof} 
 Given $U\in\mathcal{N}(\mathcal{T})$ and a good left dual $U^{\ast}$ of $U$,
 we need to prove that also $U^{\ast}$ is an object in 
 $\mathcal{N}(\mathcal{T})$. Since $U^{\ast}$ is a good left dual of $U$,
 by Lemma~\ref{l:a-tortility} it satisfies the ribbon 
 condition~\eqref{e:tortility}. So, by Lemma~\ref{l:reflexivity}, 
 $U^{\ast}$ is reflexive and so $U$ is a stable left dual of $U^{\ast}$ under 
 the reversed duality. We only need to show
 that, if we set $U^{\ast\ast}=U$ via the reversed duality, 
 then~\eqref{e:a-tortility} is satisfied. Now, 
 by Lemma~\ref{l:starstar}, 
 we have $(\theta_{U}^{\ast})^{\ast}=\theta_{U}$, so we only need
 to check
 \begin{equation}\label{l:s0s}
  \theta_{U}=\leftidx{^{U}}{\!\theta}{^{\ast}_{U^{\ast}}}\text{.}
 \end{equation}
 Applying $\lidx{U}{(\cdot)}$ to~\eqref{e:a-tortility}, we get 
 $\leftidx{^{U}}{\!\theta}{_{U^{\ast}}}=\theta^{\ast}_{U}$.
 By duality, we find~\eqref{l:s0s}.
\end{proof}  
  
\section{The quantum double of a \protect\Tcat}\label{sec:cat-d}
Let $\mathcal{T}$ be \Tcat (that again, for simplicity, we suppose
strict).
Apply the center construction obtaining the braided \Tcat 
$\mathcal{Z}(\mathcal{T})$. Then consider its twist extension 
$\bigl(\mathcal{Z}(\mathcal{T})\bigr)^{Z}$. Finally, consider its maximal
ribbon subcategory $\mathcal{D}(\mathcal{T})=
\mathcal{N}\Bigl(\bigl(\mathcal{Z}(\mathcal{T})\bigr)^{Z}\Bigr)$.
Starting from $\mathcal{T}$, we obtained a ribbon \Tcat.
However, generalizing the construction~\cite{KasTur}, we can
directly construct a ribbon \Tcat $\mathcal{D}(\mathcal{C})$,
the \textit{quantum double of $\mathcal{T}$,} isomorphic to 
$\mathcal{N}\Bigl(\bigl(\mathcal{Z}(\mathcal{T})\bigr)^{Z}\Bigr)$. 
One of the advantages is that a choice of dualities in $\mathcal{T}$,
i.e., a structure of autonomous \Tcat, induces a choice of dualities in
$\mathcal{D}(\mathcal{T})$.
\begin{itemize}
 \item The objects of $\mathcal{D}(\mathcal{T})$ are the triples
       $D=(U,\halfbrd_{\_},t)$, where
       \begin{itemize}
        \item $U$ is an object in $\mathcal{T}$,
        \item $\halfbrd_{\_}\colon U\otimes\_\to\lidx{U}{(\_)}\otimes U$ is a natural 
              isomorphism that satisfies the
              half-braiding axiom~\eqref{e:ax-c},
        \item $t\in\mathcal{T}\Bigl(U,\lidx{U}{U}\Bigr)$ 
              is an isomorphism such that
              $\biggl(\Bigl(\leftidx{^{U}}{\!t}{}\Bigr)\circ t\biggr)^{-1} 
                =\omega_{U}$.
       \end{itemize}
 \item Given two objects $D_{1}=(U_{1},\halfbrd_{\_},t_{1}), 
       D_{2}=(U_{2},\halfbrdtwo_{\_},t_{2})\in\mathcal{D}(\mathcal{T})$,
       an arrow $f\in\mathcal{D}(\mathcal{T})(D_{1},D_{2})$ is an arrow 
       $f\in\mathcal{T}(U_{1},U_{2})$ that is compatible with the
       half-braidings and the twist, i.e., it satisfies~\eqref{e:ax-arr}
       and~\eqref{e:twistator-arrow}.
 \item The tensor product of two object $D_{1}=(U_{1},\halfbrd_{\_},t_{1}),\ 
       D_{2}=(U_{2},\halfbrdtwo_{\_},t_{2})\in\mathcal{D}(\mathcal{T})$, 
       is the triple
       $D_{1}\otimes D_{2}= \bigl(U_{1}\otimes U_{2},(\halfbrd\tenbrd\halfbrdtwo)_{\_},
                   t_{1}\twten t_{2}\bigr)$.
       where $\tenbrd$ is defined as in~\eqref{e:ax-prod} and 
       $\twten$ as in~\eqref{e:twistator}.
       The tensor product of arrows is obtained by requiring 
       that the functor 
       $\mathcal{D}(\mathcal{T})\to\mathcal{T}\colon (U,\halfbrd,t)\to U$ 
       is a tensor functor.
 \item The conjugation is given by
       $\varphi_{\mathcal{D}.\beta}(U,\halfbrd_{\_},t)\eqdef \bigl(\varphi_{\beta}(U), 
          \varphi_{\beta}(\mathcal{Z}.\halfbrd_{\_}),\varphi_{\beta}(t)\bigr)$,
       for any $\beta\in\pi$ and $(U,\halfbrd_{\_},t_{U})\in\mathcal{D}(\mathcal{T})$,
       where $\varphi_{\mathcal{Z}.\beta}(\halfbrd_{\_})$ is defined as in the case
       of the center of $\mathcal{T}$ (see~\eqref{e:ax-aut}).
       For any arrow $f$ in $\mathcal{D}(\mathcal{C})$, we set
       $\varphi_{\mathcal{D}.\beta}(f)=\varphi_{\beta}(f)$.
 \item Let $D=(U,\halfbrd_{\_},\theta_{D})$ be an object in 
       $\mathcal{D}(\mathcal{T})$. We obtain a stable left dual 
       $D^{\ast}$ of $D$ in $\mathcal{D}(\mathcal{T})$ by setting
       $D^{\ast}=\Bigl(U^{\ast},\hat{\halfbrd}_{\_},\leftidx{^{U^{\ast}}}%
         {\!\! t}{^{\ast}}\Bigr)$
       and $b_{D}=b_{U}$, $d_{D}=d_{U}$,
       where $b_{U}$ and $d_{U}$ are the unit and the counit of $U$ 
       in $\mathcal{T}$.
 \item The braiding is obtained by
       $c_{D_{1},D_{2}}=\halfbrd_{U_{2}}$
       for any $D_{1}=(U_{1},\halfbrd_{\_},t_{1}), D_{2}=(U_{2},
       \halfbrdtwo_{\_},t_{2})\in\mathcal{D}(\mathcal{T})$.
 \item The twist is obtained by setting $\theta_{D}=t$
       for any $D=(U,\halfbrd_{\_},t)\in\mathcal{D}(\mathcal{T})$.
\end{itemize}
  
\begin{thm}\label{T:7}
 $\mathcal{D}(\mathcal{T})$ is a ribbon \Tcat isomorphic  
 to $\mathcal{N}\Bigl(\bigl(\mathcal{Z}(\mathcal{T})\bigr)^{Z}\Bigr)$ 
 as balanced \Tcat.
\end{thm}

\section{Yetter-Drinfeld modules and the center of 
         $\protect\Rep(H)\/$}\label{sec:YD}  
  We give the definition of a Yetter-Drinfeld module over a \Tcoalg $H$
  and we prove that the category $\YD\,(H)$ of Yetter-Drinfeld modules
  over $H$ is a braided \Tcat isomorphic to 
  $\mathcal{Z}\bigl(\Rep(H)\bigr)$.

\paragraph{\scshape Definition of a Yetter-Drinfeld module} 
  Let us fix $\alpha$ in $\pi$. 
  An \textit{\aYetter module}
  or, simply, a \textit{\Yettera} is a couple
  $V=\bigl(V,\Delta_{V}=\{ \Delta_{V,\lambda} \}_{\lambda\in\pi} \bigr)$, where  
  $V$ is an \Hamodule and, for any $\lambda\in\pi$,
  $\Delta_{V,\lambda}\colon V\to V\otimes H_{\lambda}$
  is a \klin morphism such that the following conditions are satisfied.
  \begin{itemize}\begin{subequations}
   \item $V$ is \textit{coassociative} in the sense that, 
   for any $\lambda_{1},\lambda_{2}\in\pi$, we have
        \begin{equation}\label{e:YD-coaassoc}
         (V\otimes\Delta_{\lambda_{1},\,\lambda_{2}}) \circ \Delta_{V,\,\lambda_{1}\lambda_{2}} =
         (\Delta_{V,\,\lambda_{1}}\otimes H_{\lambda_{2}}) \circ \Delta_{V,\,\lambda_{2}}\text{.}                   
         \end{equation}
   \item $V$ is \textit{counitary} in the sense that 
         \begin{equation}\label{e:YD-unit}
          (V\otimes \varepsilon)\circ \Delta_{V,1} = \Id\text{.}
         \end{equation}
   \item $V$ is \textit{crossed} in the sense that, for any $\lambda\in\pi$, 
         the diagram
         \begin{equation}\label{e:YD-crossing}
          \qquad\qquad\vcenter{\xymatrix @C=-3pc{
          \ & H_{\alpha}\otimes H_{\lambda}\otimes V\otimes H_{\lambda}
          \ar@/^/[rd]^(.65){H_{\alpha}\otimes\,\sigma\,\otimes H_{\lambda}} & \ \\ H_{\alpha\lambda}\otimes V
          \ar@/^/[ur]^(.35){\Delta_{\alpha,\,\lambda}\,\otimes\,\Delta_{V,\,\lambda}} & \ & H_{\alpha}\otimes V\otimes H_{\lambda}\otimes H_{\lambda}
          \ar[d]|{\!\mu_{V}\,\otimes\,\mu_{\lambda}} \\
          H_{\alpha\lambda\alpha^{-1}}\otimes H_{\alpha}\otimes V
          \ar@{<-}[u]|{\!\!\!\!\!\!\!\!\!\!\!\!\!\!\!\!\Delta_{\alpha\lambda\alpha^{-1},\,\alpha}\otimes V}
          \ar[d]|{\!\!\!\!\!\!\!\!\!\!\! H_{\alpha\lambda\alpha^{-1}}\otimes \mu_{V}}& \ & V\otimes H_{\lambda}\\
          H_{\alpha\lambda\alpha^{-1}}\otimes V\ar@/_/[dr]_(.35){\sigma} & 
          &  V\otimes H_{\lambda}\otimes H_{\lambda}\ar[u]|{\,V\,\otimes\,\mu_{\lambda}} \\ \
          & \ \ \ \ V\otimes H_{\alpha\lambda\alpha^{-1}}\ar@/_/[ur]_(.65){\Delta_{V,\,\lambda}\,\otimes\,\varphi_{\alpha^{-1}}}& \ }}
         \end{equation}
         commutes ($\mu_{\lambda}$ is the product of $H_{\lambda}$ while $\mu_{V}$ 
         is the \Hamodule structural map of $V$\/).
       \end{subequations}
  \end{itemize}   
       
  If, for any $v\in V$, we set
  \begin{equation}\label{e:sigma-YD}
   v_{(V)}\otimes v_{(\lambda)}\eqdef\Delta_{V,\lambda}(v)\text{,}
  \end{equation}
  then we can rewrite axiom~\eqref{e:YD-coaassoc}, as
  $(v_{(V)})_{(V)}\otimes (v_{(V)})_{(\lambda_{1})}\otimes 
   v_{(\lambda_{2})} = v_{(V)}\otimes (v_{(\lambda_{1}\lambda_{2})})'_{(\lambda_{1})}\otimes 
   (v_{(\lambda_{1}\lambda_{2})})''_{(\lambda_{2})}$.
  Similarly, we can rewrite axiom~\eqref{e:YD-unit}, as
  $\varepsilon(v_{(1)})v_{(V)}=v$.
  Finally, we can rewrite the commutativity of~\eqref{e:YD-crossing} as
  $h'_{(\alpha)}v_{(V)}\otimes h''_{(\lambda)}v_{(\lambda)} = 
   (h''_{(\alpha)}v)_{(V)}\otimes (h''_{(\alpha)}v)_{(\lambda)}
   \varphi_{\alpha^{-1}}(h'_{(\alpha\lambda\alpha^{-1})})$,
  for any $\lambda\in\pi$ and $h\in H_{\alpha\lambda}$.
  
  Given two \Yettersa $(V,\Delta_{V})$ and $(W,\Delta_{W})$, a 
  \textit{morphism of \Yettersa} $f\colon (V,\Delta_{V})\to (W,\Delta_{W})$,
  is a \Halinear map $f\colon V\to W$ such that, for any $\lambda\in\pi$,
   \begin{equation}\label{e:YD-morphisms}
   \Delta_{W,\lambda}\circ f = (f\otimes H_{\lambda})\circ \Delta_{V,\lambda}\text{.}
   \end{equation}
  With the notation provided in~\eqref{e:sigma-YD}, 
  axiom~\eqref{e:YD-morphisms} can be rewritten as
  $f(v_{(V)})\otimes v_{(\lambda)}= \bigl(f(v)\bigr)_{(W)}\otimes\bigl(f(v)\bigr)_{(\lambda)}$,
  for any $v\in V$.
  We complete the structure of the category $\YD_{\alpha}(H)$ 
  by defining the composition of morphisms 
  of \Yettersa via the standard composition of the 
  underlying linear maps, i.e., by requiring that the 
  forgetful functor 
  $\YD_{\alpha}(H)\to\Rep_{\alpha}(H)\colon(V,\Delta_{V})\mapsto V$ is faithful.
  
  Let $\YD\,(H)$ be the disjoint union of the categories 
  $\YD_{\alpha}(H)$ for all $\alpha\in\pi$. 
  The category $\YD\,(H)$ admits a structure of braided \Tcat as follows.
  
  \begin{itemize}\begin{subequations}
   \item The tensor product of a 
         \Yettera $(V,\Delta_{V})$  and a \Yetterb $(W,\Delta_{W})$ 
         (with $\alpha,\beta\in\pi$) is the \Yetterab 
         $(V\otimes W,\Delta_{V\otimes W})$, where, for any $v\in V$, $w\in W$, and $\lambda\in\pi$,
         we have
         \begin{equation}\label{e:Epicurus}
         \Delta_{V \otimes W, \lambda}(v\otimes w)
           =v_{(V)}\otimes w_{(W)}\otimes w_{(\lambda)}\varphi_{\beta^{-1}}(v_{(\beta\lambda\beta^{-1})})\text{.}
         \end{equation} 
         
         The tensor unit of $\YD\,(H)$ is the couple 
         $\tunit_{\YD}=(\Bbbk,\Delta_{\Bbbk})$, where, for any $\lambda\in\pi$ and $k\in\Bbbk$,
         we have $\Delta_{\Bbbk,\lambda}(k) = k\otimes 1_{\lambda}$.
         
         Finally, the tensor product of arrows is given 
         by the tensor product of \klin maps, i.e., by requiring that
         the forgetful functor $\YD\,(H)\to\Rep(H)\colon (V,\Delta_{V})\to V$ 
         is a tensor functor.
   \item Given $\beta\in\pi$, the conjugation functor $\lidx{\beta}{(\cdot)}$
         is obtained as follows.
         Let $\alpha$ be in $\pi$ and let $(V,\Delta_{V})$ be a \Yettera. We set 
         $\lidx{\beta\,}{(V,\Delta_{V})}=\Bigl(\lidx{\beta}{V},\Delta_{\lidx{\beta}{V}}\Bigr)$,
         where, for any $\lambda\in\pi$ and $w\in \lidx{\beta}{V}$,
         \begin{equation}\label{e:Epicurus-2}
          \Delta_{\lidx{\beta}{V},\lambda}(w) =\Biggl(\lidx{\beta}{\biggl(\Bigl(
          \lidx{\beta^{-1}}{w}\Bigr)_{(V)}\biggr)}\Biggr)
          \otimes\varphi_{\beta}\biggl(\Bigl(\lidx{\beta^{-1}}{w}\Bigr)_{(\beta^{-1}\lambda\beta)}\biggr)
          \text{.}
         \end{equation}    
         Given a morphism $f\colon (V,\Delta_{V})\to (W,\Delta_{W})$ of \Yetters, 
         for any $v\in V$, we set 
         $\bigl(\lidx{\beta}{f}\bigr)\bigl(\lidx{\beta}{v}\bigr)=
         \lidx{\beta\,}{\bigl(f(v)\bigr)}$,
         i.e., we require that the forgetful functor from 
         $\YD\,(H)\to\Rep(H)$ is a \Tfun.
   \item The braiding $c$ is obtained by setting, 
         for any \Yettera $(V,\Delta_{V})$, any 
         \Yetterb $(W,\Delta_{W})$, and any $v\in V$ and $w\in W$,
         \begin{equation}\label{e:Zatatustra}
          c_{(V,\Delta_{V}),(W,\Delta_{W})}(v\otimes w) = 
          \lidx{\alpha\,}{\bigl(s_{\beta^{-1}}(v_{(\beta^{-1})})w\bigr)}\otimes v_{(V)}\text{.}
         \end{equation}
 \end{subequations}\end{itemize}
 
 To prove that $\YD\,(H)$ is a \Tcat and that it is braided, 
 we prove before that $\YD\,(H)$ is isomorphic to 
 $\mathcal{Z}\bigl(\Rep(H)\bigr)$ as a category.
  
 \begin{thm}\label{thm:C-moll-Messe}
   The category $\YD\,(H)$ is isomorphic to the category 
   $\protect\mathcal{Z}\protect\bigl(\protect\Rep(H)\protect\bigr)$.
   This isomorphism induces on $\YD\,(H)$ the structure of crossed 
   \Tcat described above.
 \end{thm}
  
  Firstly, we construct two functors
  $\Functor_{1}\colon \mathcal{Z}\bigl(\Rep(H)\bigr)\to\YD\,(H)$, and
  $\hat{\Functor}_{1}\colon\YD\,(H)\to\mathcal{Z}\bigl(\Rep(H)\bigr)$.
  Then we prove $\Functor_{1}\circ \hat{\Functor}_{1}=\Id_{\YD\,(H)}$ 
  and $\hat{\Functor}_{1}\circ \Functor_{1}=\Id_{\bigl(\Rep(H)\bigr)}$.
  Via this isomorphism, $\YD\,(H)$ becomes a braided \Tcat. 
  We complete the proof of Theorem~\ref{thm:C-moll-Messe}
  by proving that this structure of \Tcat 
  is the structure described above.
  
  \paragraph{\scshape The functor $\Functor_{1}$}  
   Let $\alpha$ be in $\pi$ and let $(V,\halfbrd_{V})$ 
   be an object in $\mathcal{Z}_{\alpha}\bigl(\Rep(H)\bigr)$. 
   For any $\lambda\in\pi$, we set
   \begin{equation}\label{e:Mozart-1}
    \Delta_{V,\lambda}(v)=\halfbrd^{-1}_{H_{\lambda}}\biggl(\Bigl(\leftidx{^{\alpha}\!}{1}{_{\lambda}}\Bigr)
    \otimes v\biggr)\text{.}
   \end{equation}
   
  \begin{lemma}\label{l:Mozart-1}
   The couple $\bigl(V,\Delta_{V}=\{\Delta_{V,\lambda}\}_{\lambda\in\pi}\bigr)$ is a \Yettera.
   In that way, we obtain a structure of \Yetter for any object 
   in the center of $\Rep(H)$.
   With respect to this natural structure, any morphism 
   in the center of $\Rep(H)$ is also a morphism of \Yetters.
   By setting
   $
    \Functor_{1}(V,\halfbrd_{\_})=(V,\Delta_{V})$ and
    $\Functor_{1}(f)=f$,
   we obtain a functor $\Functor_{1}\colon\mathcal{Z}\bigl(\Rep(H)\bigr)
   \to\YD\,(H)$.
  \end{lemma}
   
   To prove Lemma~\ref{l:Mozart-1} we need some preliminary results.
  
   \begin{rmk}\label{rnk:Mozart}
   Given $\lambda\in\pi$, the algebra $H_{\lambda}$ is a left module over itself 
   via the action provided by the multiplication. Similarly,
   $H_{\alpha\lambda\alpha^{-1}}$ is a left module over itself. 
   By definition~\eqref{e:shift-prod} of the
   action of $H_{\alpha^{-1}\lambda\alpha}$ on the module $\leftidx{^{\alpha^{-1}}\!}{H}{_{\lambda}}$, the \klin map
   \begin{equation*}\displaystyle
     \map{\hat{\varphi}_{\alpha}}{H_{\alpha^{-1}\lambda\alpha}}{\leftidx{^{\alpha^{-1}}\!}{H}{_{\lambda}}}{h}%
     {\lidx{\alpha^{-1}}{\Bigl(\varphi_{\alpha}(h)\Bigr)}=h\Bigl(\lidx{^{\alpha^{-1}}}%
     {\!1}{_{\lambda}}\Bigr)}\text{.}
   \end{equation*}
   is \Haminusbalinear and so it is an isomorphism of \Haminusbamodules. 
   Notice that
   $\lidx{\alpha\,}{\bigl(\hat{\varphi}_{\alpha}(h)\bigr)}=\varphi_{\alpha}(h)$
   and that, for any $\alpha_{1}, \alpha_{2}\in\pi$, we have
   $\hat{\varphi}_{\alpha_{1}\alpha_{2}}=\Bigl(\leftidx{^{\alpha_{2}^{-1}}}{\!\hat{\varphi}}{_{\alpha_{1}}}\Bigr)
    \circ\hat{\varphi}_{\alpha_{2}}$.
   
   Let $X$ be an \Hlmodule (with $\lambda\in\pi$) and let
   $\tilde{x}\colon H_{\lambda}\to X$
   be the unique \Hllinear map sending $1_{\lambda}$ to $x$. We set
   \begin{equation*}
      \tilde{x}^{(\alpha)}\colon H_{\alpha\lambda\alpha^{-1}}
      \xrightarrow{\,\hat{\varphi}_{\alpha^{-1}}\,}\leftidx{^{\alpha}}{H}{_{\lambda}}
      \xrightarrow{\lidx{\alpha}{\,\tilde{x}\,}}\lidx{\alpha}{\, X}\text{.}
   \end{equation*}
   Since, for any $h\in H_{\alpha\lambda\alpha^{-1}}$,
   $\tilde{x}^{(\alpha)}(h)=\lidx{\alpha}{\,\bigl((\tilde{x}\circ\hat{\varphi}_{\alpha^{-1}})(h)\bigr)}
   =hx\text{,}$
   we have that $\tilde{x}^{(\alpha)}$ is the unique \Habalinear map sending 
   $1_{\alpha\lambda\alpha^{-1}}$ to $x\in\lidx{\alpha}{\, X}$, i.e., 
   $\tilde{x}^{(\alpha)}=\widetilde{\,\lidx{\alpha}{x}\,\,}$\!\!.
  \end{rmk}
     
  \begin{lemma}\label{l:Mozart}
   Let $V$ be a \Yettera. For any $v\in V$ and $x\in X$ we have
   \begin{equation}\label{e:l:Mozart}
    \halfbrd^{-1}_{X}(y\otimes v) = v_{(V)}\otimes v_{(\lambda)}\Bigl(\lidx{\alpha^{-1}\!}{y}\Bigr)\text{.}
   \end{equation}
  \end{lemma}
  
  \begin{proof} 
   The proof follows by the commutativity of the diagram\vspace{-18ex}
    
   \noindent\begin{minipage}{\the\textwidth}
    \begin{equation}\label{e:kouasdft}
    \quad\vcenter{\xymatrix @ur  @R=3.6pc @C=3.6pc
     {V\ar[dr]_{v\mapsto 1_{\alpha\lambda\alpha^{-1}}\otimes v\ \ \ \ \ \ \ \ }
      \ar@/^1.5pc/[ddrrr]^{\Delta_{V,\lambda}} & \ & \ & \ \\
     \ & H_{\alpha\lambda\alpha^{-1}}\otimes V \ar[dr]|{\hat{\varphi}_{\alpha^{-1}}\otimes V} 
       \ar@/_1.5pc/[ddr]_{\tilde{x}^{(\alpha)}\otimes V}
       & \ & \ \\
     \ & \ & \bigl(\leftidx{^{\alpha}}{H}{_{\lambda}}\bigr)\otimes V 
     \ar@/^1pc/[r]|(.4){^{\halfbrd^{-1}_{H_{\lambda}}}}
       \ar@/_1pc/[d]|(.4){\bigl(\lidx{\alpha}{\tilde{x}}\bigr)\otimes V\ \ \ }
       & V\otimes H_{\lambda}\ar@/^1pc/[d]^(.7){V\otimes\tilde{x}} \\
     \ & \ & \bigl(\lidx{\alpha}{\, X}\bigr)\otimes V\ar@/_1pc/[r]_(.7){\halfbrd^{-1}_{X}}
       & V\otimes X}}\end{equation}
    \end{minipage}\vspace{-18ex}
          
   \noindent for $x=\lidx{\alpha^{-1}}{y}$. 
   The top triangle commutes by definition of $\Delta_{V,\lambda}$. 
   The bottom triangle commutes by definition of $\tilde{x}^{(\alpha)}$.
   The square commutes because 
   $\halfbrd_{\_}$ is an isomorphism of functors.
  \end{proof}

 \begin{proof}[Proof of Lemma~\ref{l:Mozart-1}] 
   We check that $(V,\Delta_{V,\_})$ satisfies 
   the axioms of \Yettera, then we conclude the proof 
   of Lemma~\ref{l:Mozart-1} with the part concerning morphisms.
  
   \begin{sentence}{Coassociativity}
   Let $X_{1}$ be a \Hlonemodule and let $X_{2}$ be a \Hltwomodule, 
   with $\lambda_{1},\lambda_{2}\in\pi$. 
   By~\eqref{e:ax-c}, we have
   $\halfbrd^{-1}_{X_{1}\otimes X_{2}}= \Bigl(\halfbrd^{-1}_{X_{1}}\otimes X_{2}\Bigr)\circ 
      \Bigl(\bigl(\lidx{\alpha}{\, X_{1}}\bigr)\otimes \halfbrd^{-1}_{X_{2}}\Bigr)$, so,
      for any $v\in V$, $x_{1}\in X_{1}$ and $x_{2}\in X_{2}$, we get
   $
   v_{(V)}\otimes (v_{(\lambda_{1}\lambda_{2})})'_{(\lambda_{1})}x_{1}\otimes (v_{(\lambda_{1}\lambda_{2})})''_{(\lambda_{2})}x_{2}
     = (v_{V})_{(V)}\otimes(v_{(V)})_{(\lambda_{1})}x_{1}\otimes v_{(\lambda_{2})}x_{2}
   $.
   If we evaluate this formula for $X_{1}=H_{\lambda_{1}}$, $X_{2}=H_{\lambda_{2}}$,
   $x_{1}=\leftidx{^{\alpha}\!}{1}{_{\lambda_{1}}}$, and
   $x_{2}=\leftidx{^{\alpha}\!}{1}{_{\lambda_{2}}}$, then we get
   $
      v_{(V)}\otimes (v_{(\lambda_{1}\lambda_{2})})'_{(\lambda_{2})}\otimes (v_{\lambda_{1}\lambda_{2}})''_{(\lambda_{2})}=
      (v_{(V)})_{(V)}\otimes(v_{(V)})_{(\lambda_{1})}\otimes v_{(\lambda_{2})}$.
   \end{sentence}
  
   \begin{sentence}{Counit}
    Since we have $v_{(V)}\otimes v_{(1)}=v\otimes 1$, we get
    $\varepsilon(v_{(1)})v_{(V)}=\varepsilon(1)v=v$.
   \end{sentence}
  
   \begin{sentence}{Crossing property}
   Let $X$ be a \Hlmodule. For any $v\in V$ and $x\in X$ we have
   $
      h \halfbrd^{-1}_{X}\biggl(\Bigl(\lidx{\alpha}{x}\Bigr)\otimes v\biggr)=
      \Delta_{\alpha,\lambda}(h)\,\halfbrd^{-1}_{X}\biggl(\Bigl(\lidx{\alpha}{x}\Bigr)\otimes v\biggr)
      =\Delta_{\alpha,\lambda}(h)(v_{(V)}\otimes v_{(\lambda)})x=
      (h'_{(\alpha)}v_{(V)}\otimes h''_{(\lambda)}v_{(\lambda)})x
   $
   and
   $
      \halfbrd^{-1}_{X}\Biggl(h\biggl(\Bigl(\lidx{\alpha}{x}\Bigr)\otimes
      v\biggr)\Biggr)=
      \halfbrd^{-1}_{X}\biggl(h'_{(\alpha\lambda\alpha)}\Bigl(\lidx{\alpha}{x}\Bigr)\otimes
      h''_{(\alpha)}v\biggr)=
      (h''_{(\alpha)}v)_{(V)}\otimes (h''_{(\alpha)}v)_{(\lambda)}\varphi_{\alpha^{-1}}(h'_{(\alpha\lambda\alpha^{-1})})x$,
   so the crossing property~\eqref{e:YD-crossing} 
   follows by the \Hallinearity of $\halfbrd^{-1}_{X}$.
   This completes the proof that $(V,\Delta_{V})$ is a \Yettera. 
   \end{sentence}
   
   \begin{sentence}{Morphisms} Let $(W,\halfbrdtwo_{\_})$ 
    be another object in $\mathcal{Z}_{\alpha}\bigl(\Rep(H)\bigr)$.
    Define $\Delta_{W}$ as above for $\Delta_{V}$.
    Given any $f\colon V\to W$ in $\mathcal{Z}\Bigl(\Rep(H)\Bigr)$,
    we prove that $f$ gives rise to a morphism of \Yettersa, i.e.,
    that~\eqref{e:YD-morphisms} 
    is satisfied.  
   By the commutativity of
   $\halfbrdtwo^{-1}_{H_{\lambda}}\circ\bigl((\leftidx{^{\alpha}\!}{H}{_{\lambda}})\otimes f\bigr)
   =(f\otimes H_{\lambda})\circ \halfbrd^{-1}_{H_{\lambda}}$,
   we have
   $\bigl((f\otimes H_{\lambda})\circ \Delta_{V,\lambda}\bigr)(v) =
      \bigl((f\otimes H_{\lambda})\circ \halfbrd^{-1}_{H_{\lambda}}\bigr)\biggl(\Bigl(%
      \leftidx{^{\alpha}\!}{1}{_{\lambda}}\Bigr)\otimes v\biggr)=
      \biggl(\halfbrdtwo^{-1}_{H_{\lambda}}\circ\Bigl(\bigl(\leftidx{^{\alpha}\!}%
      {H}{_{\lambda}}\bigr)\otimes f\Bigr)\biggr)
       \biggl(\Bigl(\leftidx{^{\alpha}\!}{1}{_{\lambda}}\Bigr)\otimes v\biggr)
       = \halfbrdtwo^{-1}_{H_{\lambda}}\biggl(\Bigl(\leftidx{^{\alpha}}{1}%
       {_{\lambda}}\Bigr)\otimes f(v)\biggr)=(\Delta_{V,\lambda}\circ f)(v)$.
  \end{sentence}
   
   The proof that $\Functor_{1}$ is a functor is now trivial.
  \end{proof}

  \paragraph{\scshape The functor $\hat{\Functor}_{1}$} 
  Let $(V,\Delta_{V})$ be any \Yettera. Given $\lambda\in\pi$, for 
  any representation $X$ of $H_{\lambda}$ set
  \begin{equation}\label{e:wzn}
      \map{\halfbrd_{X}}{V\otimes X}{\bigl(\lidx{\alpha}{\, X}\bigr)\otimes V}{v\otimes x}
      {\Bigl(\lidx{\alpha}{\,\bigl(s_{\lambda^{-1}}%
      (v_{(\lambda^{-1})})x\bigr)}\Bigr)\otimes v_{(V)}}\text{.}
  \end{equation}
  
  \begin{lemma}\label{l:Mozart-2}
   The couple $(V,\halfbrd_{\_})$ is an object in
   $\mathcal{Z}\bigl(\Rep(H)\bigr)$. 
   In particular,
   \begin{equation*}\label{e:ILSS}
    \halfbrd^{-1}_{X}(y\otimes v)=v_{(V)}\otimes v_{(\lambda)}\Bigl(\lidx{\alpha^{-1}}{y}\Bigr)
   \end{equation*}
   for any $y\in\lidx{\alpha}{X}$ and $v\in V$.
   With respect to this natural structure, any morphism of 
   \Yetters gives rise to an arrow in $\mathcal{Z}\bigl(\Rep(H)\bigr)$.
   By setting
   $
    \hat{\Functor}_{1}(V,\Delta_{V})=(V,\halfbrd_{\_})$ and 
    $\hat{\Functor}_{1}(f)=f$,
   we obtain a functor from $\YD\,(H)$ to 
   $\mathcal{Z}\bigl(\Rep(H)\bigr)$. 
   The functors $\Functor_{1}$ and $\hat{\Functor}_{1}$ are
   mutually inverses.
  \end{lemma}
  
  To prove Lemma~\ref{l:Mozart-2}, we need another preliminary lemma.
  
  \begin{lemma}\label{l:KK}
   For any $v\in V$ we have
   $(v_{(V)})_{(V)}\otimes(v_{(V)})_{(\lambda)}s_{(\lambda^{-1})}(v_{(\lambda^{-1})})=v\otimes 1_{\lambda}$
   and
   $(v_{(V)})_{(V)}\otimes s_{(\lambda^{-1})}\bigl((v_{(V)})_{(\lambda^{-1})}\bigr)v_{(\lambda)}= 
    v\otimes 1_{\lambda}$.
  \end{lemma}
  
  \begin{proof}
   Since $\Delta_{V}$ is counitary, the proof follows by the commutativity
   of
   \begin{equation*}\vcenter{\xymatrix@C=5pc{%
            V\otimes H_{\lambda^{-1}}\ar[r]^{\Delta_{V,\lambda}\otimes H_{\lambda^{-1}}\ \ \ \ \ }
          & V\otimes H_{\lambda}\otimes H_{\lambda^{-1}}\ar[r]^{\ V\otimes H_{\lambda}\otimes s_{\lambda^{-1}}}
          & V\otimes H_{\lambda}\otimes H_{\lambda}\ar[d]^{V\otimes\mu_{\lambda}} \\
            V\ar[u]^{\Delta_{V,\lambda^{-1}}}\ar[d]_{\Delta_{V,\lambda}}\ar[r]|(.43){\Delta_{V,1}}
          & V\otimes H_{1}\ar[u]|{\ \ \ \ \  V\otimes\Delta_{\lambda,\lambda^{-1}}}
            \ar[d]|{\ \ \ \ \  V\otimes\Delta_{\lambda^{-1},\lambda}}
            \ar[r]|{V\otimes(\eta_{\lambda}\circ\varepsilon)} & V\otimes H_{\lambda} \\
            V\otimes H_{\lambda}\ar[r]_{\Delta_{V,\lambda^{-1}}\otimes H_{\lambda}\ \ \ \ \ } 
          & V\otimes H_{\lambda^{-1}}\otimes H_{\lambda}\ar[r]_{\ V\otimes s_{\lambda^{-1}}\otimes H_{\lambda}}
          & V\otimes H_{\lambda}\otimes H_{\lambda}\ar[u]_{V\otimes\mu_{\lambda}}}}\text{,}\end{equation*}
      where the two squares on the right-hand side are commutative 
      because of the coassociativity of $\Delta_{V}$,
      while the two squares on the left-hand  side are 
      commutative because $s$ is the antipode of a \Tcoalg.
  \end{proof}
  
  \begin{proof}[Proof of Lemma~\ref{l:Mozart-2}]
   Let us check that $(V,\halfbrd_{\_})$ is an object 
   in $\mathcal{Z}_{\alpha}\bigl(\Rep(H)\bigr)$.
   
   \begin{sentence}{Invertibility}
    Let $X$ be a representation of $H_{\lambda}$, with $\lambda\in\pi$. We set
    \begin{equation*}
      \map{\hat{\halfbrd}_{X}}{\bigl(\lidx{\alpha}%
      {\, X}\bigr)\otimes V}{V\otimes X}{y\otimes v}
      {v_{(V)}\otimes v_{(\lambda)}\Bigl(\lidx{^{\alpha^{-1}}}{\! y}\Bigr)}\text{.}
    \end{equation*}
    Let us prove that $\hat{\halfbrd}_{X}$ is the inverse 
    of $\halfbrd_{X}$. For any $v\in V$ and $x\in X$ we have
     $
       v\otimes x\overmapsto[2pc]{\halfbrd_{X}}%
       \Bigl(\lidx{\alpha}{\,\bigl(s_{\lambda^{-1}}(v_{(\lambda^{-1})})x\bigr)}\Bigr)\otimes v_{(V)}
       \overmapsto[2pc]{\hat{\halfbrd}_{X}}
       (v_{(V)})_{(V)}\otimes(v_{(V)})_{(\lambda)}s_{(\lambda^{-1})}(v_{(\lambda^{-1})})x=v\otimes x
     $
    (where the last passage follows by Lemma~\ref{l:KK}). 
    Similarly, for any $v\in V$ and $y\in\lidx{\alpha}{\, X}$ we have
     $
      y\otimes v\overmapsto[2pc]{\hat{\halfbrd}_{X}} v_{(V)}\otimes v_{(\lambda)}%
      \Bigl(\lidx{^{\alpha^{-1}}}{\! y}\Bigr)\overmapsto[2pc]{\halfbrd_{X}}
      s_{(\lambda^{-1})}\bigl((v_{(V)})_{(\lambda^{-1})}\bigr)v_{(\lambda)}y\otimes v_{(V)}
      =y\otimes\varepsilon(v_{(1)})v_{(V)}=y\otimes v$.
   \end{sentence}
  
  \begin{sentence}{Linearity}
   Let $X$ be a representation of $H_{\lambda}$, with $\lambda\in\pi$.
    It is a bit easier to prove that $\hat{\halfbrd}_{X}$ 
    (instead of $\halfbrd_{X}\/$) is \Hallinear.
     For any $v\in V$, $y\in \lidx{\alpha^{-1}}{X}$, and $h\in H_{\alpha\lambda}$, we have
    $
      h \hat{\halfbrd}_{V,X}(y\otimes v) =
      h\biggl(v_{(V)}\otimes v_{(\lambda)}\Bigl(\lidx{\alpha^{-1}}{y}\Bigr)\biggr)=
      h'_{(\alpha)}v_{(V)}\otimes h''_{(\lambda)}v_{(\lambda)}\Bigl(\lidx{\alpha^{-1}}{y}\Bigr)
    $
   and
   $
      \hat{\halfbrd}_{X}\bigl(h(y\otimes v)\bigr)=
      \hat{\halfbrd}_{X}(h'_{(\alpha\lambda\alpha^{-1})}y\otimes h''_{(\lambda)}v)
      =(h''_{(\alpha)}v)_{(V)}\otimes (h''_{(\alpha)}v)_{(\lambda)}\varphi_{\alpha^{-1}}(h'_{(\alpha\lambda\alpha^{-1})})y
  $.
   By the crossing property~\eqref{e:YD-crossing} of $(V,\Delta_{V})$, 
   these two expressions are equal.
  \end{sentence}  
 
  \begin{sentence}{Naturality}
  Let us check that $\halfbrd_{\_}$ is a natural transformation 
  from the functor $V\otimes\_$ to the functor $\lidx{V}{(\_)}\otimes V$.
  Given two representations $X_{1}$ and $X_{2}$ of $H_{\lambda}$ and a 
  \Hllinear map $f\colon X_{1}\to X_{2}$, for any $v\in V$ and $x\in X_{1}$ we have
  $\Biggl(\biggl(\Bigl(\lidx{\alpha}{f}\Bigr)\otimes V\biggr)
      \circ\halfbrd_{X_{1}}\Biggr)(v\otimes x) =
      \biggl(\Bigl(\lidx{\alpha}{f}\Bigr)\otimes V\biggr)
      \biggl(\Bigl(\lidx{\alpha}{\,\bigl(s_{\lambda^{-1}}
      (v_{(\lambda^{-1})})x\bigr)}\Bigr)\otimes v_{(V)}\biggr)
      =\Bigr(\lidx{\alpha}{\,\bigl(s_{\lambda^{-1}}(v_{(\lambda^{-1})})f(x)\bigr)}\Bigr)\
       \otimes v_{(V)}
      =\Bigl(\halfbrd_{X_{2}}\circ(V\otimes f)\Bigr)(v\otimes x)\text{.}$
  \end{sentence}
  
 \begin{sentence}{Half-braiding axiom}
  We still have to check that $(V,\halfbrd_{\_})$ satisfies the 
  half-braiding axiom~\eqref{e:ax-c}. 
  Let $X_{1}$ be a \Hlonemodule and let $X_{2}$ be a \Hltwomodule, 
  with $\lambda_{1},\lambda_{2}\in\pi$. We want
  \begin{equation*}
      \halfbrd_{X_{1}\otimes X_{2}}(v\otimes x_{1}\otimes x_{2})=
      \Biggl(\biggl(\Bigl(\lidx{\alpha}{\, X_{1}}\Bigr)\otimes\halfbrd_{X_{2}}\biggr)\circ 
      (\halfbrd_{X_{1}}\otimes X_{2})\Biggr)(v\otimes x_{1}\otimes x_{2})\text{.}
  \end{equation*}
  for any $x_{1}\in X_{1}$, $x_{2}\in X_{2}$, and $v\in V$.
  We have
  \begin{multline*}
     c_{V,X\otimes X_{1}}(v\otimes x_{1}\otimes x_{2}) =\biggl(\lidx{\alpha}{\,\Bigl(s_{\lambda_{1}^{-1}}
     \bigl((v_{\bigl((\lambda_{1}\lambda_{2})^{-1}\bigr)})''_{(\lambda_{1})}\bigr)x_{1}\Bigr)}\biggr)\otimes\\
     \otimes \biggl(\lidx{\alpha}{\,\Bigl(s_{\lambda_{2}^{-1}}
     \bigl((v_{\bigl((\lambda_{1}\lambda_{2})^{-1}\bigr)})''_{(\lambda_{2})}\bigr)x_{2}\Bigr)}\biggr)
     \otimes v_{(V)}
  \end{multline*}
  and
  \begin{multline*}\Biggl(\biggl(\Bigl(\lidx{\alpha}%
  {\, X_{1}}\Bigr)\otimes\halfbrd_{X_{2}}\biggr)\circ 
   (\halfbrd_{X_{1}}\otimes X_{2})\Biggr)(v\otimes x_{1}\otimes x_{2}) 
   = \Bigl(\lidx{\alpha\,}{\bigl(s_{\lambda_{1}^{-1}}(v_{(\lambda_{1}^{-1})})x_{1}\bigr)}\Bigr)\otimes\\
   \otimes \biggl(\lidx{\alpha\,}{\Bigl(s_{\lambda_{2}^{-1}}%
   \bigl((v_{(V)})_{(\lambda_{2}^{-1})}\bigr)x_{2}\Bigr)}\biggr)\otimes (v_{(V)})_{(V)}
  \text{.}\end{multline*}
  By the coassociativity of $\Delta_{V}$ these two expressions are equal. 
  \end{sentence}
  
  This concludes the proof that $(V,\halfbrd_{\_})$ is an object in 
  $\mathcal{Z}_{\alpha}\bigl(\Rep(H)\bigr)$.
  
  \begin{sentence}{Morphisms}
   Let $(W,\halfbrdtwo_{\_})$ be another object in 
   $\mathcal{Z}_{\alpha}\bigl(\Rep(H)\bigr)$.
   Define $\Delta_{W}$ as for $\Delta_{V}$ above. 
   Given $f\colon V\to W$ in $\mathcal{Z}\Bigl(\Rep(H)\Bigr)$,
   we prove that $f$ gives rise to a morphism of \Yettersa. 
   Let $X$ be a \Hlmodule, with $\lambda\in\pi$. Given $v\in V$ and $x\in X$ we have
   $
    \bigl(\halfbrd_{W}\circ(f\otimes X)\bigr)(v\otimes x) =
    \biggl(\lidx{\alpha\,}{\Bigl(s_{\lambda^{-1}}\bigl((f(v))_{(\lambda^{-1})}
    \bigr)x\Bigr)}\biggr)\otimes \bigl(f(v)\bigr)_{(W)}=
    \Bigl(\lidx{\alpha\,}{\bigl(s_{\lambda^{-1}}(v_{(\lambda^{-1})})x\bigr)}\Bigr)
    \otimes f\bigl(v_{(V)}\bigr)
    =\biggl(\Bigl(\lidx{\alpha}{\, X}\Bigr)\otimes f\biggr)
     \biggl(\Bigl(\lidx{\alpha\,}{\bigl(s_{\lambda^{-1}}(v_{(\lambda^{-1})})x\bigr)}
     \Bigr)\otimes v_{(V)}\biggr)=
    \Biggl(\biggl(\Bigl(\lidx{\alpha}{\, X}\Bigr)\otimes f\biggr)\circ
    \halfbrd_{V}\biggr)(v\otimes x)
   $,
   where in the second passage we used~\eqref{e:YD-morphisms}.
   \end{sentence}
   
   The proof that $\hat{\Functor}_{1}$ is a functor is now trivial. 
   We still have to check that $\Functor_{1}$ and $\hat{\Functor}_{1}$ 
   are mutually inverse.
   
   \begin{sentence}{Isomorphism}
    Let us prove that $\hat{\Functor}_{1}\circ\Functor_{1}=
    \Id_{\mathcal{Z}\bigl(\Rep(H)\bigr)}$.  Let $(V,\halfbrd_{\_})$ be an
    object in $\mathcal{Z}_{\alpha}\bigl(\Rep(H)\bigr)$, with $\alpha\in\pi$. 
    We have $\Functor_{1}(V,\halfbrd_{\_})=(V,\Delta_{V})$, where, for any $\lambda\in\pi$, 
    $v_{(V)}\otimes v_{(\lambda)}=\halfbrd_{H_{\lambda}}^{-1}\biggl(\Bigl(\leftidx{^{\alpha}}{1}{_{\lambda}}\Bigr)
    \otimes v\biggr)$.
    We set $(V,\dot{\halfbrd}_{\_})=
    (\hat{\Functor}_{1}\circ \Functor_{1})(V,\halfbrd_{\_})=\hat{\Functor}_{1}(V,\Delta_{V})$.
    Given any \Hlmodule $X$, with $\lambda\in\pi$, for any $v\in V$ and $x\in X$ we have
    $
     \dot{\halfbrd}_{X}^{-1}\biggl(\Bigl(\lidx{\alpha}{\, x}\bigr)\otimes v\biggr)
     = v_{(V)}\otimes v_{(\lambda)}x=
     \halfbrd^{-1}_{H_{\lambda}}\biggl(\Bigl(\leftidx{^{\alpha}}{1}{_{\lambda}}\Bigr)\otimes v\biggr)x
     =\halfbrd_{X}^{-1}\biggl(\Bigl(\lidx{\alpha}{\, x}\bigr)\otimes v\biggr)
    $,
    where the last passage follows by the commutativity of the 
    square in diagram~\eqref{e:kouasdft}.
   
    Let us prove that $\Functor_{1}\circ\hat{\Functor}_{1}=\Id_{\YD\,(H)}$. 
    Given $\alpha\in\pi$ and a \Yettera $X$, we have 
    $\Functor_{1}(V,\Delta_{\_})=(V,\halfbrd_{V})$, 
    where, for any representation $X$ of $H_{\lambda}$ (with $\lambda\in\pi\/$) 
    and for any $x\in X$ and $v\in V$, we have
    $\halfbrd_{X}(v\otimes x)=\Bigl(\lidx{\alpha}{\,\bigl(s_{\lambda^{-1}}%
    (v_{(\lambda^{-1})})x\bigr)}\Bigr)\otimes v_{(V)}$.
    If we set $(V,\dot{\Delta}_{V})=(\Functor_{1}\circ\hat{\Functor}_{1})(V,\Delta_{\_})
    =\Functor_{1}(V,\halfbrd)$, then we obtain
    $
     \dot{\Delta}_{\lambda}(v)=\halfbrd^{-1}_{H_{\lambda}}\biggl(%
     \Bigl(\leftidx{^{\alpha}}{1}{_{\lambda}}\Bigr)
     \otimes v\biggr)=v_{(V)}\otimes v_{(\lambda)}1_{\lambda}=v_{(V)}\otimes v_{(\lambda)}=\Delta_{\lambda}(v)
    $, where the second passage follows by~\eqref{e:ILSS}.
   \end{sentence}
   
   This concludes the proof of Lemma~\ref{l:Mozart-2}.
  \end{proof}
   
  \begin{proof}[Proof of Theorem~\ref{thm:C-moll-Messe}]
   By Lemma~\ref{l:Mozart-2}, the categories 
   $\mathcal{Z}\Bigl(\Rep(H)\Bigr)$ and $\YD\,(H)$ 
   are isomorphic via the functor $\Functor_{1}$ and the 
   functor $\hat{\Functor}_{1}$. 
   This isomorphism induces on $\YD\,(H)$ a structure of a strict \Tcat.
   
  \begin{sentence}{Components}
   Let $\alpha$ be in $\pi$. Since 
   $\YD_{\alpha}(H)=(\Functor_{1}\circ\hat{\Functor}_{1})\bigl(\YD_{\alpha}(H)\bigr)$, 
   the \alphath component of $\YD\,(H)$ is $\YD_{\alpha}(H)$. 
  \end{sentence}
 
  \begin{sentence}{Tensor category structure}
 Let $(V,\Delta_{V})$ be a \Yettera and let $(W,\Delta_{W})$ be a \Yetterb.
 Suppose $(V,\halfbrd_{\_})=\hat{\Functor}_{1}(V,\Delta_{V})$ and 
 $(W,\halfbrdtwo_{\_})=\hat{\Functor}_{1}(W,\Delta_{W})$.
 We set
   \begin{equation*}
    (V,\Delta_{V})\otimes(W,\Delta_{W})=
    \Functor_{1}\bigl(\hat{\Functor}_{1}(V,\Delta_{V})\otimes 
    \hat{\Functor}_{1}(W,\Delta_{W})\bigr)=
    \Functor_{1}\bigl(V\otimes V', 
    (\halfbrd\tenbrd \halfbrdtwo)_{\_}\bigr)\text{.}
   \end{equation*}
   Since, for any $v\in V$,
   $(\halfbrd_{\leftidx{^{\beta}}{\! H}{_{\lambda}}})^{-1}
     \biggl(\Bigl(\leftidx{^{\alpha\otimes\beta}}{1}{_{\lambda}}\Bigr)\otimes v\biggr)
    = v_{(V)}\otimes\hat{\varphi}_{\beta^{-1}}(v_{(\beta\lambda\beta^{-1})})$,
   we obtain
   $\Delta_{V\otimes W,\lambda}(v\otimes w)
     =\bigl((\halfbrd\twten\halfbrdtwo)_{H_{\lambda}}\bigr)^{-1}
    \biggl(\Bigl(\leftidx{^{\alpha\otimes\beta}}{1}{_{\lambda}}\Bigr)\otimes v\otimes w\biggr)
    =\Bigl(\bigl(V\otimes (\halfbrdtwo_{H_{\lambda}})^{-1}\bigr)\circ
    \bigl((\halfbrd_{\leftidx{^{\beta}}{\! H}{_{\lambda}}})^{-1}\otimes W\bigr)\Bigr)
    \biggl(\Bigl(\leftidx{^{\alpha\otimes\beta}}{1}{_{\lambda}}\Bigr)\otimes v\otimes w\biggr)
     =v_{(V)}\otimes w_{(W)}\otimes w_{(\lambda)}\varphi_{\beta^{-1}}(v_{(\beta\lambda\beta^{-1})})$, 
     for any $v\in V$ and $w\in W$.
   The part concerning the tensor unit of $\YD\,(H)$ is trivial.
  \end{sentence}
   
 \begin{sentence}{Conjugation}
  The \Tcat structure of $\YD\,(H)$ is completed by setting, for any $\beta\in\pi$,
  $\lidx{\beta}{(\cdot)}=\Biggl(\YD\,(H)\xrightarrow{\,\hat{\Functor}_{1}\,}
   \mathcal{Z}\bigl(\Rep(H)\bigr)
   \xrightarrow{\,\lidx{\beta}{(\cdot)}\,}\mathcal{Z}\bigl(\Rep(H)\bigr)
   \xrightarrow{\,\Functor_{1}\,}\YD\,(H)\Biggr)$.
  In particular, given $\alpha\in\pi$ and a \Yettera $(V,\Delta_{V})$, if 
  $(V, \halfbrd_{\_})=\hat{\Functor}_{1}(V,\Delta_{V})$, then, for any $\lambda\in\pi$ 
  and $v\in V$, we get
  $\Delta_{\lidx{\beta}{V},\lambda}\Bigl(\lidx{\beta\,}{v}\Bigr)
   =\Bigl(\lidx{\beta\,}{\halfbrd}\Bigr)_{H_{\lambda}}
   \biggl(\Bigl(\leftidx{^{\beta\alpha\beta^{-1}}}{\!\! 1}{_{\lambda}}\Bigr)
    \otimes\Bigl(\lidx{\beta\,}{v}\Bigr)\biggr)
   =\lidx{\beta\,}{\bigl(v_{(V)}\otimes\hat{\varphi}_{\beta}(v_{(\beta^{-1}\lambda\beta)})\bigr)}
   =\Bigl(\lidx{\beta\,}{(v_{(V)})}\Bigr)\otimes v_{(\beta^{-1}\lambda\beta)}$.
  By setting $w=\lidx{\beta}{\, v}$, we get~\eqref{e:Epicurus-2}.
 \end{sentence}
 
 \begin{sentence}{Braiding}
  The braiding in $\YD\,(H)$ is obtained by setting
  $
   c_{(V,\Delta_{V}),(W,\Delta_{W})}=\Functor_{1}(c_{\hat{\Functor}_{1}(V,\Delta_{V}),
   \hat{\Functor}_{1}(W,\Delta_{W})})= \halfbrd_{W}
  $, for any $(V,\Delta_{V}),(W,\Delta_{W})\in\YD\,(H)$,
  where $(V,\halfbrd_{\_})=\hat{\Functor}_{1}(V,\Delta_{V})$.
  By definition~\eqref{e:wzn} of $\halfbrd_{\_}$, we 
  get~\eqref{e:Zatatustra}.
 \end{sentence}
 
   This concludes the proof of the theorem.
 \end{proof}
 
\section{The \protect\Tcoalg $\protect\overline{D}(H)\/$
         and its representations}\label{sec:CR}
 \label{sec:sismology}
 We define a quasitriangular \Tcoalg  $\protect\overline{D}(H)$.
 This \Tcoalg is the mirror of the \Tcoalg $D(H)$ defined 
 in~\cite{Zunino-CHC}. To do this, we need to introduce another
 \Tcoalg $\mHdualcop$ (the mirror of $\Hdualcop$ 
 defined in~\cite{Zunino-CHC}).
 Then, we discuss the structure of a module over 
 $\overline{D}(H)$. More in detail, we prove that a \kvector 
 space $V$ is a \Dbarmodule if and only if it is both a \Hmodule 
 and a \mHmodule and the actions of $H$ and $\mHdualcop$ satisfy 
 a compatibility condition. 
 Finally, we prove that $\mathcal{Z}\bigl(\Rep(H)\bigr)$ and 
 $\Rep\bigl(\overline{D}(H)\bigr)$ are isomorphic as braided 
 \Tcats.
 
 \paragraph{\scshape Definition of $\mHdualcop$}
 The \Tcoalg $\mHdualcop$
 is defined as follows.
 \begin{itemize}
  \item  For any $\alpha\in\pi$, the component $\mHdualcop_{\alpha}$ 
         is equal to $\bigoplus_{\beta\in\pi}H^{\ast}_{\beta}$ as a vector space,
         with the product of $f\in H^{\ast}_{\gamma}$ and $g\in H^{\ast}_{\delta}$ 
         (with $\gamma,\delta\in\pi\/$) given by the linear map $fg\in H^{\ast}_{\gamma\delta}$ 
         defined by
         $\langle fg, x\rangle = \langle f, x'_{(\gamma)}\rangle\,\langle g, x''_{(\delta)}\rangle$
         for any $x\in H_{\gamma\beta}$.
         The unit of $H_{\alpha}$ is the morphism $\varepsilon\in H^{\ast}_{1}\subset 
         \mHdualcop_{1}$.
  \item  The comultiplication $\Delta^{\mdcsym}$ is defined by setting, 
         for any $\alpha,\beta,\gamma\in\pi$ and $f\in H^{\ast}_{\gamma}$,
         $\Delta^{\mdcsym}_{\alpha,\beta}(f)=\Delta_{\beta}(f)\in H^{\ast}_{\beta\gamma\beta^{-1}}\otimes H^{\ast}_{\beta}$, where
         $\bigl\langle\Delta_{\beta}(f), x\otimes y\bigl\rangle =\bigl\langle f, y\varphi_{\beta^{-1}}(x)\bigr\rangle$,
         for any $x\in H_{\beta\gamma\beta^{-1}}$.
         We introduce the notation
         \begin{equation*}
          f_{\primus.\beta}\otimes f_{\semel.\beta}\eqdef\Delta_{\beta}(f)\text{.}
         \end{equation*}
        The counit $\varepsilon^{\mdcsym}\colon \mHdualcop_{1}\to\Bbbk$ is given by
        $\langle\varepsilon^{\mdcsym}, f\rangle = \langle f, 1_{\gamma}\rangle$,
        for any $f\in H^{\ast}_{\gamma}$, with $\gamma\in\pi$.
 \item  For any $\alpha\in\pi$, the component $s^{\mdcsym}_{\alpha}$ of 
        the antipode $s^{\mdcsym}$ of $\mHdualcop$ sends $f\in H_{\gamma}$ to
        $s_{\alpha}^{\mdcsym}(f)=\Bigl\langle f, 
          \varphi_{\alpha^{-1}}\bigl(s^{-1}_{\alpha\gamma\alpha^{-1}}(\_)\bigr)\Bigr\rangle
          \in H^{\ast}_{\alpha\gamma^{-1}\alpha^{-1}}$.
 \item Finally, for any $\beta\in\pi$, we set $\varphi^{\mdcsym}_{\beta}=\varphi^{\ast}_{\beta^{-1}}$.
 \end{itemize}
 
 \paragraph{\scshape 
   The \protect\Tcoalg 
   $\protect\overline{D}(H)\/$}\label{p:indefinitamente}
 The  \Tcoalg $\overline{D}(H)$ is defined as follows
 (see the detailed description of the mirror $D(H)$ of 
 $\overline{D}(H)$ in~\cite{Zunino-CHC}).
 \begin{itemize}
     \item For any $\alpha\in\pi$, the \alphath component
           of $\overline{D}(H)$, denoted $\overline{D}_{\alpha}(H)$,
           is equal to $H_{\alpha}\otimes\bigoplus_{\beta\in\pi}H_{\beta}^{\ast}$ , as a vector space. 
           Given $h\in H_{\alpha}$ and $F\in\bigoplus_{\beta\in\pi}H_{\beta}^{\ast}$, the element in
           $\overline{D}_{\alpha}(H)$ corresponding
           to $h\otimes F$ is denoted $h\circledast F$. 
           The product in $\overline{D}_{\alpha}(H)$ is given by
           \begin{equation*}
               (h\circledast f)\,(k\circledast g) = h''_{(\alpha)}k\circledast f\bigl\langle g, s^{-1}_{\delta^{-1}}
               (h'''_{(\delta^{-1})})\_\varphi_{\alpha^{-1}}(h'_{(\alpha\delta\alpha^{-1})})\bigr\rangle
           \end{equation*}
           for any $h, k\in H_{\alpha}$, $f\in H_{\gamma}^{\ast}$, and $g\in H^{\ast}_{\delta}$, with $\gamma, \delta\in\pi$.
           $H_{\alpha}$ has unit $1_{\alpha}\circledast\varepsilon$.
           
           The algebra structure of $D_{\alpha}(H)$ is uniquely defined 
           by the condition that the inclusions
           $H_{\alpha}, \mHdualcop_{\alpha}\hookrightarrow D_{\alpha}(H)$ are algebra morphisms 
           and that the relations
           \begin{subequations}\label{e:Stranamore}
            \begin{equation}\label{e:Strabamore-a}
             (1_{\alpha}\circledast f)\,(h\circledast\varepsilon) = h\circledast f
           \end{equation}
           and
           \begin{equation}\label{e:Stranamore-b}
            (h\circledast\varepsilon)\,(1_{\alpha}\circledast f)=h''_{(\alpha)}\circledast
            \bigl\langle f, s^{-1}_{\gamma^{-1}}(h'''_{(\gamma^{-1})})\_
            \varphi(h'_{(\alpha\gamma\alpha^{-1})})\bigr\rangle\text{,}
            \end{equation} 
          \end{subequations}
          (for any $h\in H_{\alpha}$ and $f\in H_{\gamma}^{\ast}$, with $\gamma\in\pi\/$) 
          are satisfied.
    \item The comultiplication is given by
           \begin{equation*}
               (h\otimes F)'_{\overline{(\alpha)}}\otimes (h\otimes F)''_{\overline{(\beta)}} = 
               (h'_{(\alpha)}\circledast F_{\primus,\beta})\otimes (h''_{(\beta)}\circledast F_{\semel,\beta})\text{,}
           \end{equation*}
           for any $\alpha,\beta\in\pi$, $h\in H_{\alpha\beta}$ and $F\in\mHdualcop$.
           The counit is given by
           $ \langle\varepsilon,h\circledast f\rangle=\langle\varepsilon, h\rangle\,\langle f,1_{\gamma}\rangle $,
           for any $h\in H_{1}$ and $f\in H^{\ast}_{\gamma}$, with $\gamma\in\pi$.
     \item The antipode is given by
           $
               \overline{s}_{\alpha}(h\circledast F)= \bigl(s_{\alpha}(h)\circledast\varepsilon\bigr)\, 
               \bigl(1_{\alpha}\otimes s^{\mdcsym}_{\alpha}(F)\bigr)
           $,
           for any $\alpha\in\pi$, $h\in H_{\alpha}$, and $F\in\mHdualcop$.
     \item The conjugation is given by 
           $
               \overline{\varphi}_{\beta}(h\circledast f)=\varphi_{\beta}(h)\otimes\varphi^{\ast}_{\beta^{-1}}(f)
           $
           for any $\alpha,\beta\in\pi$, $h\in H_{\alpha}$, and $f\in H^{\ast}_{\gamma}$, with $\alpha,\gamma\in\pi$. 
     \item The universal \Rmatrix $\overline{R}$ of 
           $\overline{D}(H)$ is given by
           \begin{equation*}
             \overline{R}_{\alpha,\beta} = 
             \overline{\xi}_{(\alpha).i}\otimes
             \overline{\zeta}_{(\beta).i}=1_{\alpha}\circledast e^{\beta^{-1}.i}\otimes s_{\beta^{-1}}
             (e_{\beta^{-1}.i})\circledast\varepsilon
           \end{equation*}
           for any $\alpha,\beta\in\pi$, where $(e_{\beta.i})$ is a basis of $H_{\beta}$ 
           and $(e^{\beta.i})$ the dual basis. The inverse of $R_{\alpha,\beta}$ is
           $\hat{R}_{\alpha,\beta}=\hat{\xi}_{(\alpha).i}\otimes\hat{\zeta}_{(\beta).i}
           = 1_{\alpha}\circledast e^{\beta.i}\otimes e_{\beta.i}\circledast\varepsilon$.
 \end{itemize}
  
 \paragraph{\scshape The category $\protect\Rep(H,\protect\mHdualcop,\circledast)\/$}
  Given $\alpha\in\pi$, a \textit{\HHamodule} is a \kvector space $V$ endowed 
  with both a structure of left module over $H_{\alpha}$ and a structure of 
  left module over $\mHdualcop_{\alpha}=\mHdualcop_{1}$ 
  (via an action denoted $\maction\/$) satisfying the compatibility 
  condition 
  \begin{equation}\label{e:PeterSellers}
     h\bigl(f\maction v) = 
     \bigl\langle f, s^{-1}_{\gamma^{-1}}(h'''_{(\gamma^{-1})})\_\varphi_{\alpha^{-1}}(h'_{\alpha\gamma\alpha^{-1}})\bigr\rangle
     \maction (h''_{(\alpha)}v)\text{.}
  \end{equation}
  for any $v\in V$, $h\in H_{\alpha}$, and $f\in H^{\ast}_{\gamma}$, with $\gamma\in\pi$\text{.}
  A morphism of \HHamodules is a morphisms that is both a morphism 
  of \Hamodules and a morphism of \mHmodules.
  In that way, with the obvious composition, we obtain the category 
  $\Rep_{\alpha}(H,\protect\mHdualcop,\circledast)$  of \HHamodules.
  The disjoint union $\protect\Rep(H,\protect\mHdualcop,\circledast)$
  of the categories $\Rep_{\alpha}(H,\protect\mHdualcop,\circledast)$ for all $\alpha\in\pi$ 
  is a braided \Tcat as follows.
 \begin{itemize}
 \item $\Rep_{\alpha}(H,\protect\mHdualcop,\circledast)$ is the \alphath component 
       of $\Rep(H,\protect\mHdualcop,\circledast)$.
 \item Given $\alpha,\beta\in\pi$, let $U$ be an object in 
 $\Rep_{\alpha}(H,\protect\mHdualcop,\circledast)$ and
       let $V$ be an object in $\Rep_{\beta}(H,\protect\mHdualcop,\circledast)$.
       The tensor product $U\otimes V$ of \HHmodules is given by the 
       tensor product of $U$ and $V$ as both \Hamodules and \MHmodules, 
       i.e., given $u\in U$ and $v\in V$, the action of $h\in H_{\alpha\beta}$ and, 
       respectively, $f\in H_{\gamma}^{\ast}$ (with $\gamma\in\pi\/$) on $u\otimes v$
       given by
       $
        h(u\otimes v)= h'_{(\alpha)}u\otimes h''_{(\beta)}v$ and $
        f\maction(u\otimes v)= f_{\primus.\beta}\maction u\otimes f_{\semel.\beta}\maction v\text{.}
       $
 \item The conjugation is obtained in the obvious way by the 
       conjugation of $\Rep(H)$ and the conjugation of
       $\Rep(\mHdualcop)$.
 \item The braiding is obtained by setting, 
       \begin{equation}\label{e:hyerogliph}
        \map{c_{U,V}}{U\otimes V}{\Bigl(\lidx{\alpha\,}{V}\Bigr)\otimes V}{u\otimes v}%
        {\lidx{\alpha}{\,\bigl(s_{\beta^{-1}}(e_{\beta^{-1}.i})v\bigr)}\otimes e^{\beta^{-1}.i}
        \maction u}
       \end{equation}
       for any $U\in\Rep_{\alpha}(H,\protect\mHdualcop,\circledast)$ and 
       $V\in\Rep_{\beta}(H,\protect\mHdualcop,\circledast)$, with $\alpha,\beta\in\pi$.
 \end{itemize}
 
 \begin{thm}\label{thm:ofnghis}
  $\Rep_{\alpha}(H,\protect\mHdualcop,\circledast)$  is
  isomorphic to $\Rep\Bigl(\overline{D}(H)\Bigr)$
  as a braided \Tcat.
 \end{thm}
 
 \begin{proof}
  The simplest way to prove the theorem is to construct 
  an isomorphism of categories
  $
   \Functor_{3}\colon\Rep\Bigl(\overline{D}(H)\Bigr)\to
   \Rep_{\alpha}(H,\protect\mHdualcop,\circledast)
  $
  such that $\Functor_{3}$ induces on $\Rep\Bigl(\overline{D}(H)\Bigr)$ 
  the structure of braided \Tcat described above.
  
  Let $V$ be a \Dbaramodule, with $\alpha\in\pi$. 
  Since both $H_{\alpha}$ and $\mHdualcop_{1}=\mHdualcop_{\alpha}$ can be
  identified with subalgebras of $\overline{D}_{\alpha}(H)$ via the canonical 
  embeddings, $V$ has both a natural structure of left \Hamodule and a 
  natural structure of left \mHmodule. Explicitly, for any
  $v\in V$, $h\in H_{\alpha}$, and $f\in H^{\ast}_{\gamma}\subseteq\mHdualcop_{1}$, with $\gamma\in\pi$, we set
  \begin{equation*}
   hv=(h\circledast\varepsilon)v\qquad\qquad\text{and}\qquad\qquad f\maction v= 
   (1_{\alpha}\circledast f)v\text{.}
  \end{equation*}
  
  Let us prove that the compatibility condition~\eqref{e:PeterSellers}
  is satisfied. By the associativity of the 
  action of $\overline{D}(H)$ on $V$ and by~\eqref{e:Stranamore} we get
  $
   h(f\maction v) = (h\circledast\varepsilon)\bigl((1_{\alpha}\circledast f)v\bigr)=
   \bigl((h\circledast\varepsilon)(1_{\alpha}\circledast f)\bigr)v=
   \Bigl(h''_{(\alpha)}\circledast\bigl\langle f, s_{\gamma^{-1}}(h'''_{(\gamma^{-1})}))\_
   \varphi_{\alpha^{-1}}(h'_{(\alpha\gamma\alpha^{-1}})\bigr\rangle\Bigr)v
   =\Bigl(1_{(\alpha)}\circledast\bigl\langle f, s_{\gamma^{-1}}(h'''_{(\gamma^{-1})}))\_
   \varphi_{\alpha^{-1}}(h'_{(\alpha\gamma\alpha^{-1}})\bigr\rangle\Bigr)\,(h''_{(\alpha)}\circledast\varepsilon)v
   =     \bigl\langle f, s^{-1}_{\gamma^{-1}}(h'''_{(\gamma^{-1})})\_\varphi_{\alpha^{-1}}(h'_{\alpha\gamma\alpha^{-1}})\bigr\rangle
     \maction (h''_{(\alpha)}v)$.
  We set $\Functor_{3}(V)$ equal to $V$ endowed with the structure 
  of \HHamodule described above. 
  
  Given another \Dbaramodule $V$ and a \klin morphism $f\colon V\to W$, 
  it is easy to prove that $f$ is a morphism of \Dbaramodules 
  if and only if it is both a morphism of \Hamodules and a
  morphisms of \mHmodules. 
  By setting $\Functor_{3}(f)=f$, we obviously obtained a functor.
 
  Let us prove that $\Functor_{3}$ is invertible. 
  Given a \HHamodule $W$, we define an action of $\overline{D}_{\alpha}(H)$
  on $W$ via the tensor lift of the linear map
  $H_{\alpha}\times\mHdualcop_{\alpha}\times W\to W\colon (h,F,v)\to F\maction (hv)$, 
  we have to prove that we obtained a \Dbaramodule. 
  For any $h, k\in H_{\alpha}$, $f\in H^{\ast}_{\gamma}$, and $g\in H^{\ast}_{\delta}$, with $\gamma,\delta\in\pi$, we have
  $(1_{\alpha}\circledast\varepsilon)v=\varepsilon\maction (1_{\alpha}v)=\varepsilon\maction v =v$
  and
  $(h\circledast f)\bigl((k\circledast g)v\bigr) =(h\circledast f)\bigl(g\maction (kv)\bigr)
   = f\maction \Bigl(h\bigl(g\maction(kv)\bigr)\Bigr)
   = f\maction \bigl\langle g, s^{-1}_{\delta^{-1}}(h'''_{(\delta^{-1})})\_
   \varphi_{\alpha^{-1}}(h'_{(\alpha\delta\alpha^{-1})})\bigr\rangle\maction (h''_{(\alpha)}k\, v)
    = f\bigl\langle g, s^{-1}_{\delta^{-1}}(h'''_{(\delta^{-1})})\_
   \varphi_{\alpha^{-1}}(h'_{(\alpha\delta\alpha^{-1})})\bigr\rangle\maction (h''_{(\alpha)}k\, v) =
   \Bigl(h''_{(\alpha)}k\circledast f\bigl\langle g, s^{-1}_{\delta^{-1}}(h'''_{(\delta^{-1})})\_
   \varphi_{\alpha^{-1}}(h'_{(\alpha\delta\alpha^{-1})})\bigr\rangle\Bigr)v
    = \bigl((h\circledast f)\,(k\circledast g)\bigr)v$.
  
  To prove that $\Functor_{3}$ is invertible and to complete 
  the proof of the theorem is now trivial.
 \end{proof}  
 
 \section{$\protect\mathcal{Z}\protect\Bigl(%
  \protect\Rep(H)\protect\Bigr)$ and 
  $\protect\Rep\protect\Bigl(\protect\overline{D}(H)\protect\Bigr)$ 
  are isomorphic}\label{s:JLS}
  In this section we prove that 
  $\mathcal{Z}\bigl(\Rep(H)\bigr)$ and 
  $\protect\Rep\protect\Bigl(\protect\overline{D}(H)\protect\Bigr)$ 
  are isomorphic as braided \Tcats. 
  We start by defining a braided \Tfun 
  $\Functor_{2}\colon\YD\,(H)\to\Rep(H,\mHdualcop,\circledast)$.
  After that, we set 
  $\GFunctor=\Functor_{3}\circ \Functor_{2}\circ \Functor_{1}\colon\mathcal{Z}
  \bigl(\Rep(H)\bigr)\to\Rep\bigl(\overline{D}(H)\bigr)$ 
  and we prove that $\GFunctor$ is invertible.
 
 \begin{thm}\label{thm:verpastwt}
  $\mathcal{Z}\bigl(\Rep(H)\bigr)$ and 
  $\Rep\bigl(\overline{D}(H)\bigr)$ are isomorphic braided \Tcats.
 \end{thm}
 
 \paragraph{\scshape The functor $\Functor_{2}$}
 To prove 
 Theorem~\ref{thm:verpastwt}, we start by constructing the functor
 $\Functor_{2}$. For this, we need two preliminary lemmas.
 
 \begin{lemma}\label{e:HammerFilm-1}
  Let $(V,\Delta_{V})$ be a \Yettera \textup{(\/}with $\alpha\in\pi\/$\textup{).} 
  Given $f\in H^{\ast}_{\gamma}$, with $\gamma\in\pi$, for any $v\in V$ we set
   \begin{equation}\label{e:daemon}
     f\maction v\eqdef \langle f, v_{(\gamma)}\bigr\rangle v_{(V)}\text{.}
   \end{equation}
  With this action, $V$ becomes a \mHmodule and a \HHamodule.
 \end{lemma}
 
 \begin{proof} Let us prove that the action $\maction$ 
 is associative and unitary.
 
 \begin{sentence}{Associativity}
  Given $f\in H^{\ast}_{\gamma}$ and $g\in H^{\ast}_{\delta}$, with $\gamma,\delta\in\pi$, for any $v\in V$, we have 
 $
     f\maction (g\maction v)= f\maction \langle g, v_{(\delta)}\rangle v_{(V)}
     = \langle g, v_{(\delta)}\rangle\,\langle f, (v_{(V)})_{(\gamma)}\rangle(v_{(V)})_{(V)}
 $
 and
 $
    (fg)\maction v= \langle fg, v_{(\gamma\delta)}\rangle v_{(V)} = \bigl\langle f, (v_{(\gamma)})'_{(\gamma)}\rangle\,
         \langle g, (v_{(\gamma)})''_{(\delta)}\rangle v_{(V)}
 $.
 By the coassociativity of a \Yetter, these two expressions coincide.
 \end{sentence}
 
 \begin{sentence}{Unit}
 By~\eqref{e:YD-unit}, for any $v\in V$ we have
 $
     \varepsilon\maction v= \langle\varepsilon, v_{(1)}\rangle v_{(V)} = v$,
 i.e., $\maction$ is unitary.
 \end{sentence}
 
  \begin{sentence}{Compatibility condition~\eqref{e:PeterSellers}}
 Given $h\in H_{\alpha}$ and $f\in H^{\ast}_{\gamma}$, with $\gamma\in\pi$, for any $v\in V$,
 by using~\eqref{e:YD-crossing}) we get
  \begin{equation*}\begin{split}
  \ \ h(f\maction v)& = \langle f, v_{(\gamma)}\rangle\, h\, v_{(V)}
    = \langle f, v_{(\gamma)}\rangle\,\langle\varepsilon, h''\rangle h'_{(\alpha)}\, v_{(V)}
   = \bigl\langle f, \langle\varepsilon, h''\rangle v_{(\gamma)}\bigr\rangle h'_{(\alpha)}\, v_{(V)}\\ &
    = \bigl\langle f, s^{-1}_{\gamma^{-1}}(h'''_{(\gamma^{-1})})h''_{(\gamma)}v_{(\gamma)}\bigr\rangle\,
       h'_{(\alpha)} v_{(V)}\\ 
   \intertext{(by the crossing property~\eqref{e:YD-crossing})}
   & = \bigl\langle f, s^{-1}_{\gamma^{-1}}(h'''_{(\gamma^{-1})})\, (h''_{(\alpha)}v)_{(\gamma)}\,
       \varphi_{\alpha^{-1}}(h'_{(\alpha\gamma\alpha^{-1})})\bigr\rangle (h''_{(\alpha)}v)_{(V)}\\ &
   = \bigl\langle f, s^{-1}_{\gamma^{-1}}(h'''_{(\beta^{-1})})\_\varphi_{\alpha^{-1}}(h'_{(\alpha\gamma\alpha^{-1})})\bigr\rangle
       \maction (h''_{(\alpha)}v)\text{.}
 \end{split}\end{equation*}
 \end{sentence}
 \end{proof}
 
 \begin{lemma}\label{e:Otin-to-stauron}
  Take two \Yetters $(V,\Delta_{V})$ and $(W,\Delta_{W})$ and define the action 
  of $\mHdualcop_{1}$ on both $V$ and $W$ via~\eqref{e:HammerFilm-1}.
  A morphism of \Yetters $f\colon V\to W$ is also a morphism of \HHmodules.
 \end{lemma}
 
 \begin{proof}
  We only need to show that $f$ preserves the action of $\mHdualcop_{1}$.
  Let $v\in V$ and $g\in H_{\gamma}^{\ast}$, with $\gamma\in\pi$. Since $f$ is a morphism 
  of \Yetters, we have
  $
   g\maction f(v)  = 
   \Bigl\langle g, \bigl(f(v)\bigr)_{(\gamma)}\Bigr\rangle\bigl(f(v)\bigr)_{(W)}
    = \langle g, v_{(\gamma)}\rangle f(v_{(V)}\bigr) = f\bigl(\langle g,v_{(\gamma)}\rangle\, v_{(V)}\bigr)
     = f(g\maction v)
  $.
 \end{proof}
 
 \begin{lemma}\label{l:F2}
  For any \Yetter module $(V,\Delta_{V})$, set 
  $\Functor_{2}(V,\Delta_{V})=(V,\maction)$, with the action $\maction$
  of $\mHdualcop_{1}$ on $V$ defined as in~\eqref{e:HammerFilm-1}. 
  For any morphism $f$ of \Yetters, set $\Functor_{2}(f)=f$.
  In that way, we obtain a braided \Tfun 
  $\Functor_{2}\colon\YD\,(H)\to\Rep(H,\mHdualcop,\circledast)$.
 \end{lemma}
 
 \begin{proof}
  By Lemma~\ref{e:HammerFilm-1} and Lemma~\ref{e:Otin-to-stauron}, 
  $\Functor_{2}$ is well defined. The proof that it is a functor
  (i.e., that preserves identities and composition), is trivial. 
  We have to check that it is a tensor functor,
  that it commutes with the conjugation and that it is braided.
  
  \begin{sentence}{Tensor product}
   Given $\alpha,\beta\in\pi$, let $(V,\Delta_{V})$ be a \Yettera and let 
   $(W,\delta_{W})$ be a \Yetterb 
   module. By the definition~\eqref{e:Epicurus}
   of the tensor product in $\YD\,(H)$, the action $\maction$ 
   of $\mHdualcop_{1}$ of $V\otimes W$ is given by
   $
    f\maction (v\otimes w) = \bigl\langle f, (v\otimes w)_{(\gamma)}\bigr\rangle=
    \bigl\langle f, w_{(\gamma)}\varphi_{\beta^{-1}}(v_{(\beta\gamma\beta^{-1}})\bigr\rangle\, v_{(V)}\otimes w_{(W)}
    =(f_{\primus.\beta}\maction v)\otimes (f_{\semel.\beta}\maction w)
   $,
   i.e., $\Functor_{2}$ preserves the tensor product.
   The fact that $\Functor_{2}$ preserves the tensor unit is trivial.
  \end{sentence}
  
  \begin{sentence}{Crossing}
   Let $\alpha$ and $\beta$ be in $\pi$ and let $(V,\Delta_{V})$ be a \Yettera.
   The action of $\mHdualcop_{1}$ on 
   $\lidx{\beta}{\bigl(\Functor_{2}(V,\Delta_{V})\bigr)}$ is given by
   $
    f\maction w =\biggl(\varphi^{\mdcsym}_{\beta^{-1}}(f)\maction
    \Bigl(\lidx{\beta^{-1}}{w}\Bigr)\biggr)=
    \biggl(\varphi_{\beta}^{\ast}(f)\maction\Bigl(\lidx{\beta^{-1}}{w}\Bigr)\biggr)=
    \Biggl\langle f,\varphi_{\beta}\biggl(\Bigl(\lidx{\beta^{-1}}{w}\Bigr)_{(\beta^{-1}\gamma\beta)}\biggr)\biggr\rangle\,
    \Biggl(\lidx{\beta}{\biggl(\Bigl(\lidx{\beta^{-1}}{w}\Bigr)_{(V)}
    \biggr)}\Biggr)
   $,
   for any $f\in H_{\gamma}^{\ast}$, with $\gamma\in\pi$, and $w\in \lidx{\beta}{V}$.
  By~\eqref{e:Epicurus-2}, both 
  $\lidx{\beta}{\bigl(\Functor_{2}(V,\Delta_{V})\bigr)}$ and
  $\Functor_{2}\Bigl(\lidx{\beta}{\, (V,\Delta_{V}})\Bigr)$ 
  has the same structure of \mHmodule and so of
  \HHmodule. Since both $\lidx{\beta}{\,(\cdot)}$ and $\Functor_{2}$ 
  are the identity on the morphisms, we conclude that 
  $\Functor_{2}$ commute with the conjugation and that it is a \Tfun.
 \end{sentence}

  \begin{sentence}{Braiding}
   Let $(V_{1},\Delta_{V_{1}})$ be a \Yetteraone and let $(V_{2},\Delta_{V_{2}})$ 
   be a \Yetteratwo module.
   By~\eqref{e:hyerogliph}, for any $v_{1}\in V_{1}$ and $v_{2}\in V_{2}$
   we have
   \begin{equation*}\begin{split}
    c_{\Functor_{2}(V_{1},\Delta_{V_{1}}),\Functor_{2}(V_{2},\Delta_{V_{2}})}(v_{1}\otimes v_{2})
    & = \lidx{\alpha_{1}\,}{\bigl(s_{\alpha_{2}^{-1}}(e_{\alpha_{2}^{-1}.i})v_{2}\bigr)}\otimes 
    e^{\alpha_{2}^{-1}.i}\maction v_{1}\\
   & = \lidx{\alpha_{1}\,}{\bigl(s_{\alpha_{2}^{-1}}(e_{\alpha_{2}^{-1}.i})v_{2}\bigr)}\otimes
      \langle e^{\alpha_{2}^{-1}.i}, (v_{1})_{(\alpha^{-1}_{2})}\rangle (v_{1})_{(V_{1})}\\
   &  =\lidx{\alpha_{1}}{\Bigl(s_{\alpha_{2}^{-1}}\bigl((v_{1})_{(\alpha_{2}^{-1})}\bigr)v_{2}\Bigr)}
    \otimes(v_{1})_{(V_{1})}\text{.}
   \end{split}\end{equation*} By~\eqref{e:hyerogliph} we have 
   \begin{equation*}c_{\Functor_{2}(V_{1},\Delta_{V_{1}}),\Functor_{2}(V_{2},\Delta_{V_{2}})}
   = c_{(V_1,\Delta_{V_1}),(V_2,\Delta_{V_2})}=
   \Functor_{2}(c_{(V_{1},\Delta_{V_{1}}),(V_{2},\Delta_{V_{2}})})\text{.}\end{equation*}
  \end{sentence}
  \end{proof}
 
  \paragraph{\scshape Proof of Theorem~\ref{thm:verpastwt}} 
  To prove Theorem~\ref{thm:verpastwt}, we need a preliminary lemma.
  
  \begin{lemma}\label{l:verpastwt}
   For any $f\in H^{\ast}_{\gamma}$, with $\gamma\in\pi$,
   we have
   $f=\langle f, e_{(\gamma).i}1_{\gamma}\rangle e^{(\gamma).i}$.
  \end{lemma}
  
  \begin{proof}
   By evaluating $\langle f, e_{(\gamma).i}1_{\gamma}\rangle e^{(\gamma).i}$ against a generic 
   $h\in H_{\gamma}$ we obtain
   \begin{equation*}
    \bigl\langle \langle f, e_{(\gamma).i}1_{\gamma}\rangle e^{(\gamma).i}, h\bigr\rangle= 
    \langle f, e_{(\gamma).i}1_{\gamma}\rangle\,\langle e^{(\gamma).i}, h\rangle= \langle f, h\rangle\text{.}
    \end{equation*}
  \end{proof}
  
  \begin{proof}[Proof \textup{(}of Theorem~\ref{thm:verpastwt}\/\textup{)}]
    Let us set $\GFunctor= \Functor_{3}\circ \Functor_{2}\circ 
    \Functor_{1}\colon\mathcal{Z}
    \bigl(\Rep(H)\bigr)\to\Rep\bigl(\overline{D}(H)\bigr)$.
  Since both $\Functor_{1}$ and $\Functor_{2}$ as well as 
  $\Functor_{3}$ are braided \Tfuns, $\GFunctor$ is a braided \Tfun.
   To complete the proof of Theorem~\ref{thm:verpastwt}, we only
   need to show that $\GFunctor$ is invertible.
   Given a \Dbarmodule $V$, we set
   $\hat{\GFunctor}(V)=(V,c_{V,\_})$.
   Of course, $\hat{\GFunctor}(V)$ is an half-braiding and by setting
   $\hat{\GFunctor}(f)=f$,
   for any morphism $f$ of \Yetters, we obtain a functor 
   $\hat{\GFunctor}\colon\Rep\bigl(\overline{D}(H)\bigr)\to
   \mathcal{Z}\bigl(\Rep(H)\bigr)$. Let us prove that
   $\hat{\GFunctor}$ and $\GFunctor$ are mutually inverses.
   
  \begin{sentence}{$\hat{\GFunctor}\circ \GFunctor=\Id$}
   Let $(V,\halfbrd_{\_})$ be an object in 
   $\mathcal{Z}\bigl(\Rep(H)\bigr)$. Since
   $c_{\GFunctor(V,\halfbrd_{\_}),\_}=\GFunctor(c_{(V,\halfbrd_{\_}),\_})
    =\GFunctor(\halfbrd_{\_})=\halfbrd_{\_}$,
   we get
   $(\hat{\GFunctor}\circ \GFunctor)(V,\halfbrd_{\_})=(V,c_{\GFunctor(V),\_})
    =(V,\halfbrd_{\_})$.
  \end{sentence}
  
  \begin{sentence}{$\GFunctor\circ\hat{\GFunctor}=\Id$}
   Let $V$ be a \Dbaramodule, with $\alpha\in\pi$.
   Clearly,  $(\GFunctor\circ\hat{\GFunctor})(V)$ 
   and $V$ have the same structure of \kvector spaces
   and the same structure of \Hamodule 
   (via the embedding $H_{\alpha}\hookrightarrow\overline{D}_{\alpha}(H)\/$).
   To prove $\GFunctor\circ\hat{\GFunctor}(V)=V$, 
   we only need to check that the action $\maction$
   of $\mHdualcop_{1}$ on $V$ and the action 
   $\dot{\maction}$ of $\mHdualcop_{1}$ on 
   $(\GFunctor\circ\hat{\GFunctor})(V)$ (both obtained via the embedding
   $\mHdualcop_{1}\hookrightarrow\overline{D}_{\alpha}(H)\/$) are the same.
   
   Let $f$ be in $H^{\ast}_{\gamma}$, with $\gamma\in\pi$. By observing that, for any $v\in V$,
   $c^{-1}_{V, H_{\gamma}}\Bigl(\leftidx{^{\alpha}}{1}{_{\gamma}}\Bigr)=
    \hat{\xi}_{(\alpha).i}v \otimes \hat{\zeta}_{(\gamma).i}1_{\gamma}=
    e^((\gamma).i\maction v\otimes e_{(\gamma).i}1_{\gamma}$,
   we get
   $f\dot{\maction}v=\langle f,v_{(\gamma)}\rangle v_{V}=\langle e_{(\gamma).i} 1_{\gamma}\rangle e^{(\gamma).i}=f$,
   where the last passage follows by Lemma~\ref{l:verpastwt}.
  \end{sentence}
  \end{proof}
  
  \begin{cor}
   $\mathcal{Z}\bigl(\Rep(H)\bigr)$, $\YD\,(H)$,
   $\Rep(H,\mHdualcop,\circledast)$, and $\Rep\bigl(\overline{D}(H)\bigr)$ 
   are isomorphic braided \Tcats.
  \end{cor}
  
  \begin{proof}
   We have seen that both the functor $\Functor_{1}$ and the functor
   $\Functor_{3}$ are isomorphisms of braided \Tcats.
   By Lemma~\ref{l:F2}, $\Functor_{2}$ is an braided \Tfun and,
   by Theorem~\ref{thm:verpastwt}, $\Functor_{2}$ is invertible
   with inverse 
   $\hat{\Functor}_{2}= \Functor_{1}\circ \hat{\GFunctor}\circ \Functor_{3}$.
  \end{proof}
  
  Let $\YDf(H)$ be the category of \textit{finite-dimensional \Yetters,} 
  i.e., the category of \Yetters $(V,\halfbrd_{\_})$ such that 
  $\dim_{\Bbbk}V\in\N$, and let $\Repf(H,\mHdualcop,\circledast)$ be
  the category of finite-dimensional \HHmodules.
  
  \begin{cor}\label{YDf}
   $\mathcal{Z}\bigl(\Repf(H)\bigr)$, $\YDf(H)$,
   $\Repf(H,\mHdualcop,\circledast)$, and $\Repf\bigl(\overline{D}(H)\bigr)$
   are isomorphic braided \Tcats.
  \end{cor}
 
  \begin{proof}
   The functor $\Functor_{1}$ sends the full subcategory 
   $\mathcal{Z}\bigl(\Repf(H)\bigr)$ of 
   $\mathcal{Z}\bigl(\Rep(H)\bigr)$
   to the full subcategory $\YDf(H)$ of $YD(H)$.
   Similarly, the functor $\Functor_{2}$ sends $\YDf(H)$ 
   to the full subcategory $\Repf(H,\mHdualcop,\circledast)$ 
   of $\Repf(H,\mHdualcop,\circledast)$ and the functor $\Functor_{3}$ sends
   $\Repf(H,\mHdualcop,\circledast)$ to the full subcategory 
   $\Repf\bigl(\overline{D}(H)\bigr)$ of 
   $\Rep\bigl(\overline{D}(H)\bigr)$.
  \end{proof}
  
  \begin{rmk}[modular \protect\Tcats]
  The categorical analog of the notion of modular Hopf algebra 
  is the notion of modular category~\cite{RT,Tur-QG}. 
  A \Tcat $\mathcal{T}$ is \textit{modular} when the component
  $\mathcal{T}_{1}$ is modular as a tensor category~\cite{Tur-CPC}.
  
  Let $\mathcal{R}$ be a semisimple tensor category. 
  It was proved by M{\"u}ger~\cite{Muger} that, under certain 
  conditions on $\mathcal{R}$, the center of $\mathcal{R}$,
  is modular. 
  We expect that it will be possible to generalize the result to 
  the crossed case when $\pi$ is finite.
  On the contrary, when $\pi$ is not finite, since the quantum 
  double of a semisimple \Tcoalg is not modular, the theory fails
  to be applicable to the crossed case. However, is some case,
  for instance when the isomorphism classes of the $H_{\alpha}$ 
  (for all $\alpha\in\pi\/$) are finite, $\mathcal{Z}(\mathcal{R})$ 
  should be modular, or at least, premodular in the sense 
  of Brugui{\`e}res
  and, in that case, they should give rise to a modular category.
 \end{rmk}
 
 \section{Ribbon structures}\label{sec:rstr} 
 We conclude by discussing 
 the relation between algebraic and categorical ribbon extensions.
 Let $H$ be \Tcoalg (not necessarily of a finite-type).
 Firstly, we recall the definition of the ribbon \Tcoalg 
 $\RT(H)$ (see~\cite{Zunino-CHC}). Then, we prove that the 
 categories $\Repf\bigl(\RT(H)\bigr)$ and
 $\mathcal{N}\Bigl(\bigl(\Repf(H)\bigr)^{N}\Bigr)$ are isomorphic 
 as balanced \Tcats. 
 To prove this statement we start by introducing an auxiliary 
 ribbon \Tcat $\Rib\,(H)$. 
 Then we prove that $\Repf\bigl(\RT(H)\bigr)$
 and $\Rib\,(H)$ are isomorphic as ribbon \Tcats while $\Rib\,(H)$ 
 that $\mathcal{N}\Bigl(\bigl(\Repf(H)\bigr)^{N}\Bigr)$ 
 are isomorphic as balanced \Tcats. 
 Finally, we prove that, if $H'$ is a \Tcoalg of finite-type, 
 then $\Repf\Bigl(\RT\bigl(D(H')\bigr)\Bigr)$ and 
 $\mathcal{D}\bigl(\Repf(H')\bigr)$ are isomorphic
 as ribbon \Tcats.
 
 \paragraph{\scshape The \protect\Tcoalg $\RT(H)$}
 The ribbon \Tcoalg $\RT(H)$ is defined as follow.
  \begin{itemize}
  \item For any $\alpha\in\pi$, the \alphath component of $RT(H)$, 
        denoted $RT_{\alpha}(H)$, is the vector space whose elements 
        are formal expressions $h+kv_{\alpha}$, with $h, k\in H_{\alpha}$. 
        The sum is given by
        $(h+kv_{\alpha}) + (h'+k'v_{\alpha})\eqdef (h + h') + (k + k')v_{\alpha}$,
        for any $h, h', k, k'\in H_{\alpha}$. 
        The multiplication is obtained by requiring 
        $v^{2}_{\alpha}=u_{\alpha}s_{\alpha^{-1}}(u_{\alpha^{-1}})$, i.e., by setting, 
        for any $h, h', k, k'\in H_{\alpha}$,
        $(h+kv_{\alpha})\,(h'+k'v_{\alpha}) \eqdef hh' + hk'v_{\alpha}+k\varphi_{\alpha}(h')v_{\alpha}+k\varphi_{\alpha}(k')
         u_{\alpha}s_{\alpha^{-1}}(u_{\alpha^{-1}})
          = \bigl(hh'+k\varphi_{\alpha}(k')u_{\alpha}s_{\alpha^{-1}}(u_{\alpha^{-1}})\bigr)+
           \bigl(hk'+k\varphi_{\alpha}(k')\bigr)v_{\alpha}$.

        We identify $H_{\alpha}$ with the subset $\{h+0v_{\alpha}\vert h\in H_{\alpha}\}$ 
        of $RT_{\alpha}(H)$. 
        The algebra $RT_{\alpha}(H)$ is unitary with unit $1_{\alpha}= 1_{\alpha}+0v_{\alpha}$. 
        Moreover, for any $\alpha,\beta\in\pi$, 
        we have $R_{\alpha,\beta}\in H_{\alpha}\otimes H_{\beta}\subset RT_{\alpha}(H)\otimes RT_{\beta}(H)$.
  \item The comultiplication is given by
        $\Delta_{\alpha,\beta}(h+kv_{\alpha\beta}) =\Big(h'_{(\alpha)}+k'_{(\alpha)}\tilde{\xi}_{(\alpha).i}
          \tilde{\zeta}_{(\alpha).j}v_{\alpha}\Bigr)\otimes\biggl(h''_{(\beta)}+k''_{(\beta)}\tilde{\zeta}_{(\beta).i}\varphi_{\alpha^{-1}}
         \Bigl(\tilde{\xi}_{(\alpha\beta\alpha^{-1}).j}\Bigr)v_{\beta}\biggr)$,
        for any $h, k\in H_{\alpha}$, and $\alpha,\beta\in\pi$.
        Further, the counit is given by
        $\langle\varepsilon, h+kv_{\alpha}\rangle =\langle\varepsilon, h\rangle+\langle\varepsilon,k\rangle$, for any $h, k\in H_{1}$.
  \item The antipode is given by
        $s_{\alpha}(h+kv_{\alpha}) = s_{\alpha}(h)+(s_{\alpha}\circ\varphi_{\alpha^{-1}})(k)v_{\alpha^{-1}}$,
        for any $h, k\in H_{\alpha}$ and $\alpha\in\pi$.
  \item Finally, the conjugation is given by
        $\varphi_{\beta}(h+kv_{\alpha})=\varphi_{\beta}(h)+\varphi_{\beta}(k)v_{\beta\alpha\beta^{-1}}$,
        for any $h, k\in H_{\alpha}$ and $\alpha,\beta\in\pi$.
 \end{itemize}
 
 \paragraph{\scshape The category $\Rib\,(H)$.} 
 The ribbon \Tcat $\Rib\,(H)$, is defined as follows.
 \begin{itemize}\begin{subequations}\label{e:povera-civetta}
  \item For any $\alpha\in\pi$, the objects of the component $\Rib_{\alpha}(H)$ 
        of $\Rib\,(H)$ are the couples $(M,\ribbon)$, where $M$ 
        is a finite-dimensional representation
        of $H_{\alpha}$ and $\ribbon\colon M\to\lidx{M}{M}$ is a \Halinear 
        isomorphism such that, if we set
        \begin{equation*}
         \ribbon^{2}=\Biggl(M\xrightarrow{\,\,\ribbon\,\,}
         \lidx{M}{M}\xrightarrow{\,\lidx{M}{\ribbon}\,}
         \lidx{M\otimes M}{M}\Biggr)
        \end{equation*}
        and
        $
         \ribbon^{-2}=(\ribbon^{2})^{-1}
        $, 
        then we have
        \begin{equation}\label{e:RibCat}
         \ribbon^{-2}\Bigl(\lidx{\alpha^{2}}{m}\Bigr)= u_{\alpha}s_{\alpha^{-1}}(u_{\alpha^{-1}})m
        \end{equation}
        for any $m\in M$ (where $u_{\alpha}$ is the \alphath Drinfeld element).
  \item Given two objects $(M_{1},\ribbon_{1})$,
        $(M_{2},\ribbon_{2})\in\Rib_{\alpha}(H)$, a morphism 
        $f\colon (M_{1},\ribbon_{1})\to (M_{2}, \ribbon_{2})$ is a \Halinear 
        map $f\colon M_{1}\to M_{2}$ such that
        $t_2\circ f= \bigl(\lidx{\alpha}{f}\bigr)\circ t_{1}$.
  \item The composition of morphisms in $\Rib_{\alpha}(H)$ is obtained
        in the obvious way via the compositions of \Hamodules.
  \item The tensor product of two objects 
        $(M,\ribbon),(M',\ribbon')\in\Rib\,(H)$, is given by
        \begin{equation}\label{e:povera-civetta-1}
         (M,\ribbon)\otimes (M',\ribbon') = 
         (M\otimes M',\ribbon\twten\ribbon')\text{,}
        \end{equation}
        where we recall that
        \begin{equation}\label{e:povera-civetta-2}
         \ribbon\twten\ribbon'
          = \Biggl(\biggl(\leftidx{^{\bigl(\leftidx{^{M}}{M}{'}\bigr)}}%
          {\!\ribbon}{}\biggr)\otimes\leftidx{^{M}}{\!\ribbon'}{}\Biggr)\circ 
         c_{\leftidx{^{M}}{\! M}{'},M}\circ c_{M,M'}\text{,}
        \end{equation}
        where $c$ is the standard braiding in $\Repf(H)$.
  \item The tensor unit of $\Rib\,(H)$ is the couple $(\Bbbk,\Id_{\Bbbk})$,
        where $\Bbbk$ is a \Honemodule via the counit of $H$.
  \item The action of the crossing on an object $(M,\ribbon)\in\Rib\,(H)$ 
        is obtained by setting
        $
         \lidx{\beta}{(M,\ribbon)}=\Bigl(\lidx{\beta}{M},\lidx{\beta}{f}\Bigr)
        $ for any $\beta\in\pi$,
        while the action of the crossing on morphisms in obtained 
        by requiring that the forgetful functor
        $\Rib\,(H)\to\Repf(H)\colon(M,\ribbon)\to M$ is a \Tfun.
  \item The braiding is given by
        $
         c_{(M,\ribbon),(M',\ribbon')}=c_{M,M'}
        $
        for any $(M,\ribbon),(M',\ribbon')\in\Rib\,(H)$.
  \item The twist is given by
        $
         \theta_{(M,\ribbon)}=\ribbon
        $
        for any $(M,\ribbon)\in\Rib\,(H)$.
  \item The duality in $\Rib\,(H)$ is obtained as follows.
        Let $(M,\ribbon)$ be an object in $\Rib\,(H)$.
        The dual object of $(M,\ribbon)$ is given by the couple
        $\Bigl(M^{\ast},\leftidx{^{M^{\ast}}}{\!\ribbon}{^{\ast}}\Bigr)$
        where $M^{\ast}$ is the dual \Hmodule of $M$ (via
        unit the $b_{M}$ and the counit $d_{M}$ defined 
        in~\eqref{e:tardi}). 
        Finally we set
        $b_{(M,\ribbon)}=b_{M}$ and $d_{(M,\ribbon)}=d_{M}$.
 \end{subequations}\end{itemize}
 
 \begin{thm}\label{thm:-pre-final}
  $\Rib\,(H)$ is a ribbon \Tcat isomorphic to 
  $\Repf\bigl(\RT(H)\bigr)$. Moreover, $\Rib\,(H)$ is isomorphic to 
  $\mathcal{N}\Bigl(\bigl(\Repf(H)\bigr)^{N}\Bigr)$ as a balanced \Tcat.
 \end{thm}
 
 To prove Theorem~\ref{thm:-pre-final} 
 we need three preliminary lemmas.
 
 \begin{lemma}\label{l:final-pre-final}
 For any $\alpha\in\pi$ we have
 \begin{equation}\label{e:final-pre-final}
  s_{\alpha^{-1}}(u_{\alpha^{-1}})=\xi_{(\alpha).i}s_{\alpha^{-1}}(\zeta_{(\alpha^{-1}).i})\text{.}
 \end{equation}
 \end{lemma}
 
\begin{proof}
By observing that
$s_{\alpha}(\xi_{(\alpha).i})\otimes\zeta_{(\beta).i}=\varphi_{\alpha}(\xi_{(\alpha^{-1}.i)})\otimes
\zeta_{(\beta^{-1}).i}$ 
(see~\cite{Virelizi}), we have
$\xi_{(\alpha^{-1}).i}\otimes\zeta_{(\alpha).i}
=(s^{-1}_{\alpha^{-1}}\circ\varphi_{\alpha^{-1}})(\xi_{(\alpha).i})\otimes s^{-1}_{\alpha}(\zeta_{(\alpha^{-1}).i})$,
so we get
$u_{\alpha^{-1}} = (s_{\alpha}\circ\varphi_{\alpha^{-1}})(\zeta_{(\alpha).i})\xi_{(\alpha^{-1}).i}
                 = \varphi_{\alpha^{-1}}(\zeta_{(\alpha^{-1}).i})(s^{-1}_{\alpha^{-1}}\circ\varphi_{\alpha^{-1}})(\xi_{(\alpha).i})
                 = \zeta_{(\alpha^{-1}).i}s^{-1}_{\alpha^{-1}}(\xi_{(\alpha).i})$.
By composing both sides by $s_{\alpha^{-1}}$ we get~\eqref{e:final-pre-final}.
\end{proof}

 \begin{lemma}\label{l:final-pre-final-1}
 For any $\alpha\in\pi$ and $h\in H_{\alpha}$, we have
 \begin{equation}\label{e:final-pre-final-1}
  s_{\alpha^{-1}}(u_{\alpha^{-1}})h=(s^{-1}_{\alpha}\circ s^{-1}_{\alpha^{-1}}\circ\varphi_{\alpha})(h)s_{\alpha^{-1}}(u_{\alpha^{-1}})\text{.}
 \end{equation}
\end{lemma}

\begin{proof}
 Let $k$ be in $H_{\alpha^{-1}}$. By~\eqref{Ovid-7}, we have 
 $(s_{\alpha}\circ s_{\alpha^{-1}}\circ\varphi_{\alpha^{-1}})(k)=u_{\alpha^{-1}}ku^{-1}_{\alpha^{-1}}$. By
 composing both sides by $s^{-1}_{\alpha}$ and observing that, by~\eqref{Ovid-5},
  $s_{\alpha^{-1}}(u_{\alpha^{-1}})=s^{-1}_{\alpha}(u_{\alpha^{-1}})$, we get
  $s_{\alpha^{-1}}(u_{\alpha^{-1}})(s_{\alpha^{-1}}\circ\varphi_{\alpha^{-1}})(k)=s^{-1}_{\alpha}(k)s_{\alpha^{-1}}(u_{\alpha^{-1}})$.
  For $k=(\varphi_{\alpha}\circ s^{-1}_{\alpha^{-1}})(h)$, we get~\eqref{e:final-pre-final-1}.
\end{proof}
 
 \begin{lemma}\label{l:final-lemma}
  Let $M$ be a finite-dimensional representation of $H$ 
  and let $\omega_{M}$ defined as in~\eqref{e:omega}.
  For any $m\in M$ we have
  $\Omega_{M}\Bigl(\lidx{\alpha^{2}}{m}\Bigr)=u_{\alpha}s_{\alpha^{-1}}(u_{\alpha})m$.
 \end{lemma}
 
 The proof is a long but not difficult computation and is
 omitted.

 \begin{proof}[Proof of Theorem~\ref{thm:-pre-final}] 
  $\Rib\,(H)$ is obviously a well defined category. 
  We start by proving that $\Rib\,(H)$ is isomorphic 
  to $\Rep\bigl(\RT(H)\bigr)$ as a category. 
  Let $M$ be a finite-dimensional representation of $\RT(H)$.
  Set $\map{\theta_{M}}{M}{M}{x}{\lidx{M_{1}}{(\theta x)}}$.
  Since $\theta^{-2}=s_{\alpha^{-1}}(u_{\alpha^{-1}})u_{\alpha}$ 
  (see~\cite{Virelizi}), the couple $(M,\theta_{M})$ 
  is a object in $\Rib\,(H)$.
  Conversely, let $(N,t)$ be an object in $\Rib_{\alpha}(H)$, 
  with $\alpha\in\pi$. Define the action of $\RT_{\alpha}(H)$ on $N$ via
  \begin{equation}\label{e:action-final}
   (h+k v_{\alpha})n = hn +k t^{-1}\bigl(\lidx{\alpha}{\, n}\bigr)
  \end{equation}
  for any $h,k\in H$ and $n\in N$. Let us check that the action 
  defined in~\eqref{e:action-final} is \RTalinear,
  i.e., that we provided $N$ of a structure of \RTamodule.
  For any $h_{1},k_{1},h_{2},k_{2}\in H$ and $n\in N$, we have
  $
   \bigl((h_{1}+k_{1}v_{\alpha})  \,(h_{2}+k_{2}v_{\alpha})\bigr)n 
  =\bigl(h_{1}h_{2}+k_{1}\varphi_{\alpha}(k_{2})u_{\alpha}s_{\alpha^{-1}}(u_{\alpha^{-1}})\bigr)n +
    \bigl(h_{1}k_{2}+k_{1}\varphi_{\alpha}(h_{2})\bigr)t^{-1}\bigl(\lidx{\alpha}{\, n}\bigr)
  $
  and
  $
   (h_{1}+k_{1}v_{\alpha})\bigl((h_{2}+k_{2})n\bigr) 
  =h_{1}h_{2}n+h_{1}k_{2}t^{-1}\bigl(\lidx{\alpha}{\, n}\bigr)+k_{1}t^{-1}
 \bigl(\lidx{\alpha}{\,(h_{2}n)}\bigr)
   +k_{1}t^{-1}\Bigl(\lidx{\alpha}{\,\bigl(k_{2}t^{-1}
   \bigl(\lidx{\alpha}{\, n}\bigr)\bigr)}\Bigr)$.
  By observing that
  \begin{equation*}
   k_{1}t^{-1}\bigl(\lidx{\alpha}{\,(h_{2}n)}\bigr)=
   k_{1}t^{-1}\bigl(\varphi_{\alpha}(h_{2})\lidx{\alpha}{\, n}\bigr)
   = k_{1}\varphi_{\alpha}(h_{2})t^{-1}\bigl(\lidx{\alpha}{\, n}\bigr)
  \end{equation*}
  and that
  \begin{equation*}\begin{split}
   k_{1}t^{-1}\Bigl(\lidx{\alpha}{\,\bigl(k_{2}t^{-1}\bigl(\lidx{\alpha}%
   {\, n}\bigr)\bigr)}\Bigr)
   & =k_{1}t^{-1}\biggl(\varphi_{\alpha}(k_{2})\lidx{\alpha}{\,\Bigl(t^{-1}
   \bigl(\lidx{\alpha}{\, n}\bigr)\Bigr)}\biggr)
   =k_{1}\varphi_{\alpha}(k_{2})\bigl(t^{-1}\circ\leftidx{^{\alpha}}{t}{^{-1}}\bigr)
   \bigl(\lidx{\alpha}{\, n}\bigr)\\
   & =k_{1}\varphi_{\alpha}(k_{2})t^{-2}\Bigl(\lidx{\alpha^{2}}{n}\Bigr)=
   k_{1}\varphi_{\alpha}(k_{2})u_{\alpha}s_{\alpha^{-1}}(u_{\alpha^{-1}})n\text{,}
  \end{split}\end{equation*}
  we obtain
  \begin{multline*}
  (h_{1}+k_{1}v_{\alpha})\bigl((h_{2}+k_{2})n\bigr)
  =\\ \bigl(h_{1}h_{2}+k_{1}\varphi_{\alpha}(k_{2})u_{\alpha}s_{\alpha^{-1}}(u_{\alpha^{-1}})\bigr)n +
    \bigl(h_{1}k_{2}+k_{1}\varphi_{\alpha}(h_{2})\bigr)t^{-1}\bigl(\lidx{\alpha}{\, n}\bigr)\text{,}
  \end{multline*}
  i.e., the action defined in~\eqref{e:action-final} is \RTalinear. 
  To complete the proof that $\Rib\,(H)$ and $\Rep\bigl(\RT(H)\bigr)$ 
  are isomorphic categories is now trivial.
  
  Since $\Rib\,(H)$ and $\Rep\bigl(\RT(H)\bigr)$ are isomorphic, 
  the ribbon \Tcat structure of $\Rep\bigl(\RT(H)\bigr)$
  induces a ribbon \Tcat structure on $\Rib\,(H)$.
  This is the structure described above.
  The only nontrivial point is to show that the tensor
  product induced in $\Rib\,(H)$ is the same given 
  in~\eqref{e:povera-civetta}.
  Let $(M_{1},t_{1})$ be an object in $\Rib_{\alpha}(H)$ and let $(M_{2},t_{2})$ 
  be an objects in $\Rib_{\beta}(H)$, with $\alpha,\beta\in\pi$.
  To prove that $\theta_{M_{1}\otimes M_{2}}=t_{1}\twten t_{2}$, an easy but long
  computation shows that,
  for any $m_{1}\in M_{1}$ and $m_{2}\in M_{2}$, both
  $\theta_{M_{1}\otimes M_{2}}(m_{1}\otimes m_{2})$
  and $(t_{1}\twten t_{2})(m_{1}\otimes m_{2})$ are equal to
  $\lidx{\alpha\beta}{\bigl(\theta_{\alpha}\zeta_{(\alpha).i}\xi_{(\alpha).j}m_{1}\bigr)}\otimes 
        \lidx{\alpha\beta}{\bigl(\theta_{\beta}\varphi_{\alpha^{-1}}(\xi_{(\alpha\beta\alpha^{-1}).i})\zeta_{(\beta).j}m_{2}\bigr)}$.
 
 Let us prove that $\Rib\,(H)$ and $\mathcal{N}
 \Bigl(\bigl(\Repf(H)\bigr)^{N}\Bigr)$ are isomorphic.
 If $(M,\theta_{M})$ is an object in $\Rib_{\alpha}(H)$, with $\alpha\in\pi$, then,
 by Lemma~\ref{l:final-lemma},
 for any $m\in M$ we have $\Omega_{M}\bigl(\lidx{\alpha^{-2}}{m}\bigr)=u_{\alpha}\theta^{-2}_{M}(m)$,
 so that $(M,\theta_{M})$ is an object in 
 $\mathcal{N}\Bigl(\bigl(\Repf(H)\bigr)^{N}\Bigr)$.
 Conversely, if $(M,t)$ is an object in 
 $\mathcal{N}\Bigl(\bigl(\Repf(H)\bigr)^{N}\Bigr)$,
 then, by Lemma~\ref{l:reflexivity}, and 
 Lemma~\ref{l:final-lemma} for any $m\in M$ we have
 $t\bigl(\lidx{\alpha^{-2}}{m}\bigr)=u_{\alpha}\theta^{-2}_{M}(m)$, i.e, 
 $(M,t)$ is an object in $\Rib\,(H)$. 
 The rest follows easily.
 \end{proof}
 
 Since, by Theorem~\ref{thm:-pre-final},
 $\mathcal{N}\Bigl(\bigl(\Repf(H)\bigr)^{N}\Bigr)$ 
 is isomorphic to $\Repf\bigl(\RT(H)\bigr)$, the
 balanced \Tcat $\mathcal{N}\Bigl(\bigl(\Repf(H)\bigr)^{N}\Bigr)$ 
 has also a natural structure of ribbon \Tcat.
 In particular, when $H$ is the quantum double of a finite-type
 \Tcoalg $H'$, this structure of a ribbon \Tcat is the same induced
 by the isomorphism between 
 $\Repf\Bigl(\RT\bigl(\overline{D}(H')\bigr)\Bigr)$
 and $\mathcal{D}\bigl(\Repf(H')\bigr)$, 
 so that we obtain the following corollary.
 
 \begin{cor}\label{dgfngj}
  If $H'$ is a finite-type \Talg, then 
  $\Repf\Bigl(\RT\bigl(\overline{D}(H')\bigr)\Bigr)$ 
  and $\mathcal{D}\bigl(\Repf(H')\bigr)$ are isomorphic as
  ribbon \Tcats. 
 \end{cor}
 
\bibliography{CategoricalDouble}

\providecommand{\bysame}{\leavevmode\hbox to3em{\hrulefill}\thinspace}
\providecommand{\MR}{\relax\ifhmode\unskip\space\fi MR }
\providecommand{\MRhref}[2]{%
  \href{http://www.ams.org/mathscinet-getitem?mr=#1}{#2}
}
\providecommand{\href}[2]{#2}
\begin{thebibliography}{10}

\bibitem{Bruguieres}
Alain Brugui{\`e}res, Private communication.

\bibitem{Drn}
Vladimir~G. Drinfeld, \emph{Quantum groups}, {P}roceedings of the
  {I}nternational {C}ongress of {M}athematicians \textup{(}{B}erkeley,
  1986\/\textup{)} (Providence, RI), vol.~1, Amer. Math. Soc., 1987,
  pp.~798--820. \MR{89f:17017}

\bibitem{JS0}
Andr{\'e} Joyal and Ross Street, \emph{The geometry of tensor calculus. {I}},
  Adv. in Math. \textbf{88} (1991), no.~1, 55--112. \MR{92d:18011}

\bibitem{JS-tortile}
\bysame, \emph{Tortile {Y}ang-{B}axter operators in tensor categories}, J. Pure
  Appl. Algebra \textbf{71} (1991), no.~1, 43--51. \MR{92e:18006}

\bibitem{JS}
\bysame, \emph{Braided tensor categories}, Adv. in Math. \textbf{102} (1993),
  no.~1, 20--78. \MR{94m:18008}

\bibitem{Kas}
Christian Kassel, \emph{Quantum {G}roups}, Graduate Texts in Mathematics~155,
  Springer-Verlag, New York, 1995. \MR{96e:17041}

\bibitem{KasTur}
Christian Kassel and Vladimir~G. Turaev, \emph{Double construction for monoidal
  categories}, Acta Math. \textbf{175} (1995), no.~1, 1--48. \MR{96m:18015}

\bibitem{LeTur}
Thang Le and Vladimir~G. Turaev, \emph{Quantum groups and ribbon
  ${G}$\nobreakdash-\hspace{0pt}categories}, Preprint
  {\texttt{math.QA/0103017}}, 2001.

\bibitem{ML-coherence}
Saunders Mac~Lane, \emph{Natural associativity and commutativity}, Rice Univ.
  Studies \textbf{49} (1963), no.~4, 28--46. \MR{30 \#1160}

\bibitem{ML}
\bysame, \emph{Categories for the {W}orking {M}athematician}, second ed.,
  Graduate Texts in Mathematics~5, Springer-Verlag, New York, 1998. \MR{1 712
  872}

\bibitem{Muger}
Michael M{\"u}ger, \emph{From subfactors to categories and {T}opology~{II}},
  Preprint {\texttt{math.CR/0111205}}, 2001.

\bibitem{RT}
Nicolai~Yu. Reshetikhin and Vladimir~G. Turaev, \emph{Ribbon graphs and their
  invariants derived from quantum groups}, Comm. Math. Phys. \textbf{127}
  (1990), no.~1, 1--26. \MR{91c:57016}

\bibitem{Shum}
Mei~C. Shum, \emph{Tortile tensor categories}, J. Pure Appl. Algebra
  \textbf{93} (1994), no.~1, 57--110. \MR{95a:18008}

\bibitem{Street-double}
Ross Street, \emph{The quantum double and related constructions}, J. Pure Appl.
  Algebra \textbf{132} (1998), no.~2, 195--206. \MR{99e:16054}

\bibitem{Sweedler}
Moss~E. Sweedler, \emph{Hopf {A}lgebras}, Mathematics Lecture Note Series, W.
  A. Benjamin, Inc., New York, 1969. \MR{40/5705}

\bibitem{Tur-QG}
Vladimir~G. Turaev, \emph{Quantum {I}nvariants of {K}nots and 3-{M}anifolds},
  de Gruyter Studies in Mathematics~18, Walter de Gruyter \& Co., Berlin, 1994.
  \MR{95k:57014}

\bibitem{Tur-pi}
\bysame, \emph{Homotopy field theory in dimension $2$ and group-algebras},
  Preprint {\texttt{math.QA/9910010}}, 2000.

\bibitem{Tur-CPC}
\bysame, \emph{Homotopy field theory in dimension $3$ and crossed
  group-categories}, Preprint {\texttt{math.GT/000529}}, 2000.

\bibitem{Virelizi}
Alexis Virelizier, \emph{{H}opf group-coalgebras}, Preprint
  \texttt{math.QA/0012073}. To appear in \textit{{J}. {P}ure {A}ppl.
  {A}lgebra}, 2000.

\bibitem{Virelizi2}
\bysame, \emph{Alg{\`e}bres de {H}opf gradu{\'e}es et fibr{\'e}s plats sur le
  $3$\nobreakdash-\hspace{0pt}vari{\'e}t{\'e}s \textup{(In English)}}, Ph.D.
  thesis, Universit{\'e} Louis Pasteaur (Strasbourg), 2001.

\bibitem{Zunino-CHC}
Marco Zunino, \emph{Double constructions for crossed hopf coalgebras},
  \texttt{math.QA/0212192}.

\end{thebibliography}

\end{document}